\documentclass{amsart}
\usepackage[utf8]{inputenc}
\usepackage[T1]{fontenc}
\usepackage[T2A]{fontenc}
\usepackage{enumerate}
\usepackage[table,xcdraw]{xcolor}
\usepackage{imakeidx}
\usepackage{amsfonts}
\usepackage{amssymb}
\usepackage{amsmath}
\usepackage{amsthm}
\usepackage{tikz-cd}
\usepackage{float}
\usepackage{graphicx}
\usepackage{tikz}
\usepackage{thmtools}
\usepackage{thm-restate}
\usepackage{mathrsfs}
\usepackage{enumitem}
\usepackage{mathtools}
\usepackage{multirow}
\usepackage{wrapfig}
\usepackage{verbatim}
\usepackage[all]{xy}
  \SelectTips{cm}{10}     
  \everyxy={<2.5em,0em>:} 

\usepackage{hyperref}
\usepackage[backend=biber, style=alphabetic, sorting=nyt, maxcitenames=10, minbibnames=10, maxbibnames=10]{biblatex}
\theoremstyle{definition}
\newtheorem{theorem}{Theorem}[section]
\newtheorem{definition}[theorem]{Definition}

\newtheorem*{notation*}{Notation}

\newtheorem{example}[theorem]{Example}
\newtheorem{construction}[theorem]{Construction}
\newtheorem{remark}[theorem]{Remark}
\newtheorem{corollary}[theorem]{Corollary}
\newtheorem{proposition}[theorem]{Proposition}
\newtheorem{conjecture}[theorem]{Conjecture}

\newtheorem{lemma}[theorem]{Lemma}

\newtheorem{claim}[theorem]{Claim}

\newtheorem*{acknowledgements*}{Acknowledgements}

\newtheorem*{question*}{Question I}
\newtheorem*{question'}{Question II}
\newtheorem*{question''}{Question III}

\newtheorem*{theorem*}{Theorem}

\newcommand{\ZZ}{{\mathbb Z}}
\newcommand{\NN}{{\mathbb N}}
\newcommand{\RR}{{\mathbb R}}
\newcommand{\QQ}{{\mathbb Q}}

\newcommand{\CC}{{\mathbb C}}

\newcommand{\PP}{{\mathbb P}}

\newcommand{\Affine}{{\mathbb A}}

\newcommand{\cX}{{\mathcal X}}

\newcommand{\cU}{{\mathcal U}}

\newcommand{\cO}{{\mathcal O}}

\newcommand{\cM}{{\mathcal M}}
\newcommand{\cL}{{\mathcal L}}

\newcommand{\cE}{{\mathcal E}}
\newcommand{\cD}{{\mathcal D}}

\newcommand{\cA}{{\mathcal A}}

\DeclareMathOperator{\Spec}{Spec}

\DeclareMathOperator{\Sym}{Sym}

\DeclareMathOperator{\Aut}{Aut}

\DeclareMathOperator{\Diagat}{D_{at}}

\DeclareMathOperator{\qf}{qf}

\DeclareMathOperator{\Cl}{Cl}
\DeclareMathOperator{\NS}{NS}

\DeclareMathOperator{\GVF}{GVF}

\DeclareMathOperator{\Val}{Val}

\DeclareMathOperator{\an}{an}

\DeclareMathOperator{\Tr}{Tr}
\DeclareMathOperator{\Norm}{N}

\DeclareMathOperator{\cdiv}{div} 

\DeclareMathOperator{\Adiv}{ADiv}

\DeclareMathOperator{\Apic}{APic}

\DeclareMathOperator{\LdivQ}{LDiv_{\QQ}}

\DeclareMathOperator{\ULatQ}{ULat_{\QQ}}

\DeclareMathOperator{\LPicQ}{LPic_{\QQ}}

\DeclareMathOperator{\Div}{Div} 
\DeclareMathOperator{\ord}{ord} 
\DeclareMathOperator{\height}{ht}

\DeclareMathOperator{\VF}{VF}

\DeclareMathOperator{\id}{id}
\DeclareMathOperator{\ev}{ev}
\DeclareMathOperator{\qftp}{qftp}
\DeclareMathOperator{\tp}{tp}
\DeclareMathOperator{\loc}{loc}

\DeclareMathOperator{\BPF}{BPF}

\DeclareSymbolFont{yhlargesymbols}{OMX}{yhex}{m}{n} 
\DeclareMathAccent{\yhwidehat}{\mathord}{yhlargesymbols}{"62}

\newcommand{\smallk}{k}

\newcommand{\Field}{F}

\newcommand{\adiv}{\widehat{\cdiv}}

\newcommand{\oOF}{\ov{\Omega}_{\Field}}

\newcommand{\OF}{\Omega_{\Field}}
\newcommand{\OFp}{\Omega_{\Field'}}

\newcommand{\GVFk}{\ov{\smallk(t)}}
\newcommand{\GVFQ}{\ov{\QQ}}

\newcommand{\ov}{\overline}

\newcommand{\meet}{\wedge}
\newcommand{\join}{\vee}
\newcommand{\bigmeet}{\bigwedge}
\newcommand{\bigjoin}{\bigvee}
\newcommand{\aerror}{e}

\title{Globally valued fields: foundations}
\author{Itaï Ben Yaacov}
\address{Itaï Ben Yaacov -- Universite Claude Bernard Lyon 1, CNRS, Centrale Lyon, INSA Lyon, Université Jean Monnet, ICJ UMR5208, 69622 Villeurbanne, France}
\author{Pablo Destic}
\address{Pablo Destic -- Universite Claude Bernard Lyon 1, CNRS, Centrale Lyon, INSA Lyon, Université Jean Monnet, ICJ UMR5208, 69622 Villeurbanne, France}
\author{Ehud Hrushovski}
\address{Ehud Hrushovski -- Mathematical Institute, University of Oxford, Andrew Wiles Building, Radcliffe Observatory Quarter, Woodstock Road, Oxford OX2 6GG, UK}
\author{Michał Szachniewicz}
\address{Michał Szachniewicz -- Institute for Advanced Study, 1 Einstein Drive, Princeton, New Jersey, 08540, USA}
\date{}

\begin{document}

\pagestyle{plain}

\begin{abstract}
    We present foundations for \textit{globally valued fields}, a class of enriched fields capturing some aspects of the geometry of global fields, based on the product formula. We show that this class is axiomatizable in the framework of unbounded continuous logic, opening the way to decidability
    and transfer theorems.  We provide a dictionary between various data defining such extra structure: syntactic, Arakelov or cycle-theoretic, and measure theoretic. In particular we obtain a representation theorem relating globally valued fields and adelic curves defined by Chen and Moriwaki.
\end{abstract}

\maketitle

\tableofcontents

\section{Introduction}

A global field $K$, i.e., either a number field, or a finite extension of the field $k(t)$ (for some finite base-field $k$), satisfies a \textit{product formula}, which means that there is a natural choice of absolute values $(|\cdot|_i)_{i \in I}$ on $K$ such that
\begin{equation} \label{eq:zero}
    \prod_{i \in I} |a|_i = 1 \textnormal{ for all } a \in K^{\times}.
\end{equation}
In fact, if one of the absolute values either has discrete value group with finite residue field or is Archimedean, the product formula characterises this class of fields, see the Artin-Whaples theorem \cite[Theorem 3]{artin1945axiomatic}. Geometrically, the product formula can be interpreted as the fact that every `function' $a \in K^{\times}$ has the same number of zeroes and poles, counted with multiplicities. Indeed, by taking $-\log$ of the previous formula, we get
\begin{equation} \label{eq:one}
    \sum_{i \in I} -\log|a|_i = 0,
\end{equation}
where we can think of $-\log|a|_i$ as of the order of vanishing of $a$ at $i \in I$. This point of view and the analogy between function fields and number fields are very useful, see \cite[Chapter I.6]{Hartshorne_AG}, and \cite[Chapter IV]{Basic_number_theory_Andre_Weil}.

However, despite Artin-Whaples, there are `product formulas' appearing in mathematical practice that do not come from number fields or function fields of curves over finite fields. First, note that function fields of curves over any field give such examples, as if $C$ is a smooth projective curve over any field, then non-zero rational functions on $C$ have the same number of zeroes and poles (counted with multiplicities). Even more generally, if $X$ is a one-dimensional proper regular Noetherian scheme, it may happen that rational functions on $X$ satisfy the product formula with respect to orders of vanishing at closed points of $X$. A spectacular example of this behaviour is when $X$ is the Fargues-Fontaine curve, see \cite{FF_curve_Lurie}. Also, in the world of analytic functions, Jensen's formula
\begin{equation} \label{eq:two}
    \sum_{0<|a|<r} \ord_a(f) \log \frac{r}{|a|} + \frac{1}{2\pi} \int_{0}^{2\pi} -\log |f(r e^{i\theta})| d\theta = -\log |f(0)|
\end{equation}
can be thought of as a product formula, with an error $-\log|f(0)| = O_r(1)$. This observation is a part of Vojta's dictionary between diophantine approximation and value distribution theory, see \cite{Vojta1987DiophantineAA}. Moreover, by using the above Jensen's formula for $r=1$ together with $p$-adic Jensen's formulas \cite[Remark 2.8]{khoai1989heights}, one can deduce that for $f \in \QQ(z)^{\times}$ we have
\begin{equation} \label{eq:three}
    \sum_{p} -\log \|f\|_p + \frac{1}{2\pi} \int_{0}^{2\pi} -\log |f(e^{i\theta})| d\theta + \sum_{x \in \PP_{\QQ}^1} \ord_x(f) \cdot \height(x) = 0,
\end{equation}
where $\|f\|_p$ is the $p$-adic Gauss norm of $f$ (i.e., the supremum of the $p$-adic norm of $f$ on the $p$-adic unit disc, if $f \in \QQ[z]$) and for $x \in \PP_{\QQ}^1$ the number $\height(x)$ is the Weil height of $x$. 

Equation (\ref{eq:three}) is a natural example of a product formula with a continuous measure (at least over Archimedean absolute values), contrary to the counting measure case in Equation (\ref{eq:one}). Moreover, even if we restrict to number fields, while each has a discrete space of places, when we pass to $\GVFQ$, the space becomes locally compact but not discrete. This suggests that a more general notion of `global fields' should include not only counting measures on the space of absolute values on a field. There are a few candidates for such a notion in the literature. First, Gubler in \cite{GublerMfields} defined an \textit{M-field}, which is a field $K$ together with a measure space $(M,\mu)$ and $\mu$-almost everywhere defined maps 
\[ (v \in M) \mapsto (|x|_v \in \RR_{\geqslant 0}) \]
for $x \in K$, such that $\mu$-almost everywhere $|\cdot|_v$ is an absolute value, and $\log|z|_v \in L^1(M,\mu)$ for $z \in K^{\times}$. This formulation allows errors like in Equation (\ref{eq:two}), and if we want the product formula to hold, we need to additionally assume that $\log|z|_v \in L_0^1(M,\mu)$ (integral must be zero). The consideration of such more general product formulas found applications soon after the appearance of \cite{GublerMfields}, as Moriwaki \cite{moriwaki2000arithmetic} used height functions on finitely generated fields over $\QQ$, to recover the Raynaud's theorem (Manin-Mumford conjecture). Later, Chen and Moriwaki introduced \textit{(proper) adelic curves} in \cite{Adelic_curves_1}, and Ben Yaacov and Hrushovski introduced \textit{(globally or) multiply valued fields} in \cite{GVF2}. Both of these notions are defined similarly to M-fields (satisfying the product formula). Let us also mention the work \cite{yuan2024adelic} in which Yuan and Zhang extend the theory of Moriwaki developed in \cite{moriwaki2000arithmetic, moriwaki2001canonical}.

In this paper we present the basics of the theory of globally valued fields, compare it to adelic curves, and describe how to see those as models of an unbounded continuous logic theory $\GVF_e$ (where $e \geqslant 0$). The logic interpretation of globally valued fields allows one to define various operations (e.g. ultraproducts) typical of model theory, in the context of global fields. In particular this implies some transfer principles from number fields to function fields, or in Nevanlinna theory, see Example \ref{example:ultraproduct_transfer_number_fields_function_fields} and Example \ref{example:ultraproduct_transfer_Nevalinna_function_fields} respectively. 
The connection to stability in model theory would be desirable, especially because it has been useful in other algebro-geometric contexts, like in the solution of the geometric Mordell-Lang conjecture \cite{hrushovski1996mordell}.

This paper is partially based on the unpublished notes \cite{GVF3}, however, some definitions are adjusted, and some new results are obtained. We hope that it will provide a useful reference for globally valued fields in the future. 

Let us describe the main definitions and results of this text in more detail. Recall that an adelic curve is a field $\Field$ together with a measure space $(\Omega, \cA, \nu)$ and a map $\phi: \omega \mapsto |\cdot|_\omega$ from $\Omega$ to the space of absolute values $M_{\Field}$ on $\Field$, such that for all $a \in \Field^{\times}$ the function
\[ \omega \mapsto \log|a|_\omega \]
is integrable (and measurable) on $\Omega$. Given such data, we can define functions $h:\Field^n \to \RR \cup \{-\infty\}$ given by
\begin{equation} \label{eq:four}
    h(a_1, \dots, a_n) = \int_{\Omega} \max(\log|a_1|_{\omega}, \dots, \log|a_n|_{\omega}) d\nu(\omega).
\end{equation}
Denote by $\mathbb{A}(F) := \bigcup\limits_{n\in\mathbb{N}} 
F^n$ and $\mathbb{P}(F) := \bigcup\limits_{n\in\mathbb{N}} \mathbb{P}^n(F)$. We denote by $\otimes$ the Segre product $\mathbb{A}^n(F)\times \mathbb{A}^m(F) \rightarrow \mathbb{A}^{mn}(F)$, defined by \[ (x_1,\ldots,x_n)\otimes(y_1,\ldots,y_m) := (x_1y_1,\ldots,x_1y_m,x_2y_1,\ldots,x_ny_m).\]
We also denote by $\otimes : \mathbb{P}^n(F)\times \mathbb{P}^m(F) \rightarrow \mathbb{P}^{mn+m+n}(F)$ the induced map on projective spaces (the usual Segre product).
It is easy to check that the above functions form a `height' on $\Field$, which we define below.

\begin{definition} \label{definition:heights:intro}
    A \emph{height} on $F$ is a function $h:\mathbb{A}(F)\rightarrow \mathbb{R}\cup\{-\infty\}$ satisfying the following axioms, for some \emph{Archimedean error} $\aerror\geqslant 0$.
\[
    \begin{array}{lll}
        \textnormal{Height of zero:} & \forall x\in F^n, & h(x) = -\infty \Leftrightarrow x=0 \\
        \textnormal{Height of one:} & & h(1,1) = 0 \\
        \textnormal{Invariance:} & \forall x\in F^n,\, \forall \sigma\in \Sym_n, & h(\sigma x) = h(x) \\
        \textnormal{Additivity:} & \forall x\in F^n,\, \forall y\in F^m, & h(x\otimes y) = h(x) + h(y) \\
        \textnormal{Monotonicity:} & \forall x\in F^n,\, \forall y\in F^m, & h(x) \leqslant h(x,y) \\
        \textnormal{Triangle inequality:} & \forall x,y\in F^n, & h(x+y) \leqslant h(x,y) + \aerror
    \end{array}
\]
    A height $h$ is called \emph{global} if $h\restriction_{F^\times} \equiv 0$. It then induces a function $\mathbb{P}(F) \rightarrow \mathbb{R}_{\geqslant 0}$ and we write $\height(x) := h[x:1]$ for $x \in \Field$.
\end{definition}

The main theorem of this text says that the data of (global) heights on $\Field$ is equivalent to a few other types of data, each of a different flavour. We first state it and then explain the nature of various structures appearing therein.

\begin{theorem} \label{theorem:intro:main}(Theorem \ref{theorem_equivalent_MVF_structure_definitions})
    There is an explicit bijective correspondence between the following structures on $\Field$:
    \begin{enumerate}
        \item Local terms (Definition \ref{definition_GVF_Globally_valued_fields_by_predicates}),
        \item Heights (Definition \ref{definition:heights:intro}),
        \item Positive functionals (Definition \ref{definition_positive_functional}),
        \item Families of local measures (Definition \ref{definition_family_of_local_measures}),
        \item Renormalisation classes of admissible measures (Definition \ref{definition_global_measure}),
        \item Equivalence classes of lattice valuations (Definition \ref{definition:lattice_MVF_structure}).
    \end{enumerate}
    Moreover, this correspondence respects the `global' ones, i.e., the ones satisfying the product formula. 
\end{theorem}

We define a \textit{multiply valued field} (resp. \textit{globally valued field}) to be a field with any one of these equivalent structures (resp. that are global); thus, if K is a globally valued field, we may freely talk about the Weil heights, the local terms,
the valuation into the $L^1$-lattice, or about a renormalisation class of a globalizing measure. We use the abbreviations MVF and GVF respectively.

If we require a specific presentation, we will simply say ``a GVF, presented via local terms'' (or via Weil heights etc.)  A product formula field is another synonym that has been used, especially for the local term presentation.  

Note  that ``local'' is a slightly misleading shorthand;  ``mean local terms'' would be a better description, as they arise by integration over
the various localities (and similarly for Weil heights over MVFs).

We note that while MVFs can be useful for local or ``mean-local'' work, they do not form an axiomatisable class. Any axiomatisable class containing all MVFs would allow non-zero elements whose height is $-\infty$.  We can avoid this by adding an axiom insisting on a lower bound for the heights of non-zero elements; but any such bound will actually imply that $0$ is a lower bound, leading precisely to GVFs.

The above theorem says that GVF structures on a field can be seen as coming from seemingly different worlds.
\begin{itemize}
    \item First, one can define them through heights (points (1) and (2)), which easily fits into the realm of unbounded continuous logic.
    \item Second, on a countable characteristic zero field, point (3) specialises to the data of a linear functional $\Apic(\Field) \to \RR$ which is non-negative on the effective cone, where $\Apic(\Field)$ is the colimit of the spaces of adelic divisors on arithmetic varieties with function field $\Field$. Also, if $\Field$ is finitely generated over a subfield $k$ (of any characteristic), then global functionals trivial on $k$ correspond to linear functionals $\NS(\Field/k) \to \RR$ that are non-negative on the effective cone, where $\NS(\Field/k)$ is the colimit of the Néron–Severi groups of proper varieties over $k$ with function field $\Field$. We give details on this approach in Section \ref{section_representations}.
    \item Third, GVF structures can come from measures on the space of pseudo-absolute values like in Equation (\ref{eq:three}), or from adelic curves. These correspond to points (4) and (5). Another measure theoretic way of presenting a GVF structure on a field is a global `valuation' valued in a Banach $L^1$-lattice as in point (6), see Definition \ref{definition:lattice_MVF_structure} or \cite{allcock2008banach} for an example.
\end{itemize}
We summarise this dictionary in the table below.
\begin{table}[ht]
\centering
\begin{tabular}{|c|c|c|c|}
    \hline
    Weil heights & Algebraic cycles & Measure theory & Banach lattices \\
    \hline
    Global & $\Apic(\Field) \to \RR$ & Proper adelic & $v:\Field^{\times} \to \Gamma \cup \{\infty\}$ \\
   
    heights & or $\NS(\Field/k) \to \RR$ & curves, M-fields & `valuation' valued in \\
   
    as in & positive on the & satisfying the & an $L^1$-lattice $\Gamma$, \\
    
    Definition \ref{definition:heights:intro} & effective cone & product formula & see Definition \ref{definition:lattice_MVF_structure} \\

    & \cite{szachniewicz2023existential} or \cite{GVF2, dang2022intersection} & \cite{Adelic_curves_1, GublerMfields} & and \cite{allcock2008banach}\\
    \hline
\end{tabular}
\end{table}

In particular, if the field $\Field$ is countable, it follows from the proof of Theorem \ref{theorem:intro:main} that every (global) height on $\Field$ actually comes from a (proper) adelic curve structure, as in Equation (\ref{eq:four}).
\begin{corollary}[Corollary \ref{corollary_every_MVF_structure_comes_from_an_adelic_curve}] \label{corollary:intro:representation_by_adelic_curves}
    Any GVF structure on a countable field $\Field$ has a proper adelic curve representation with $\Omega = M_{\Field}$.
\end{corollary}
This means that GVF structures on a countable field $\Field$ can be seen as equivalence classes of adelic structures on $\Field$ with respect to the equivalence relation defined by inducing the same height via Equation (\ref{eq:four}). This corollary is important for the future development of globally valued fields, as it allows one to use results e.g. from \cite{Adelic_curves_2,Adelic_curves_3,Adelic_curves_4, dolce2023generalisation, liu2024arithmetic, Wenbin_Luo_relative_Siu, A_Sedillot_diff_of_relative_volume} to study model theory of GVF structures. For example, the Hilbert-Samuel theorem in \cite{Adelic_curves_3} (see also \cite{yaacov2018estimates, yaacov2024vandermonde}) may be useful in proving amalgamation theorems over some globally valued fields. However, it is natural to treat $h$ in Definition \ref{definition:heights:intro} as a continuous logic predicate, so by `model theory' we mean the model theory in \textit{unbounded continuous logic} introduced in \cite{Ben_Yaacov_unbdd_cont_FOL}. Another reason for considering continuous logic is that for example $\GVFQ$ with Weil heights considered in discrete logic (with an additional sort for values of heights) is undecidable, as it interprets\footnote{By the formula: $v(x)=0$ for almost all non-Archimedean $v$, and $v(i-x / i+x) = 0$ for almost all Archimedean $v$.} the totally real algebraic integers (known to be undecidable by \cite{robinson1959undecidability}).

We provide a survey of the unbounded continuous logic (in a simplified setting) in Section \ref{section_continuous_logic}. In Section \ref{section_existential_closedness}, we give some background on existential closedness in this logic and state the following problem.
\begin{conjecture}(Conjecture \ref{conjecture_GVF_model_companion})
    The theory $\GVF_e$ has a model companion, for every Archimedean error $e \geqslant 0$.
\end{conjecture}
This conjecture asks whether we can axiomatize (in a first order way) existentially closed globally valued fields. Intuitively, a GVF $\Field$ is existentially closed, if all polynomial equations together with height inequalities that are consistent (solved in an extension of $\Field$), have a solution in $\Field$. Existential closedness of GVFs implies existence of small sections of certain adelic vector bundles, see e.g. Minkowskian line bundles in \cite[Section 5.2.8]{Adelic_curves_3}. This weaker condition may be equivalent to existential closedness, but in any case forms an important part and is in itself axiomatisable.

In \cite{GVF2, szachniewicz2023existential} respectively, the globally valued fields $\ov{k(t)}$ (with $\height(t) = 1$, $\textnormal{ht}|_{k} =0$) and $\GVFQ$ (equipped with Weil projective heights) have been proven to be existentially closed. Over $\GVFQ$ this implies, for example, that a finiteness statement (e.g. Bogomolov) is true in every GVF extension of $\GVFQ$ if and only if it is true in $\GVFQ$ uniformly in parameters. Conjecture \ref{conjecture_GVF_model_companion} would imply that all existentially closed globally valued fields are elementarily equivalent to one of these two kinds of GVFs. In particular, the problem of approximating the essential infimum of an adelic line bundle (up to some $\varepsilon > 0$) would be decidable.

In the course of proving Theorem \ref{theorem:intro:main} we introduce a few gadgets associated to $\Field$ and prove some results that may be of independent interest. We call a function $|\cdot|:\Field \to [0, \infty]$ a \textit{pseudo-absolute value} if $|xy|=|x|\cdot|y|$, $|x+y| \leqslant |x|+|y|$, $|0|=0, |1|=1$ (with $0 \cdot \infty$ and $\infty \cdot 0$ undefined). By an \textit{absolute value} we mean a pseudo-absolute value valued in $[0, \infty)$. We provide a Riemann-Zariski-Berkovich characterisation of the space of pseudo-absolute values $\Theta_{\Field}$ on $\Field$ with the pointwise convergence topology. Let us remark that a similar result (in the case of a field extension while fixing the absolute value on the smaller field) was independently obtained by Antoine Sédillot in \cite[Birational approach I.11.2.2]{SédillotPhD2024}.
\begin{proposition}(Proposition \ref{proposition_theta_is_a_homeo}) \label{intro:proposition_theta_intro}
    Let $\Field$ be of characteristic zero. Then there is a homeomorphism $\theta:\Theta_{\Field} \to \cM(\Field)$ between the space of pseudo-absolute values on $\Field$, and
    \[ \cM(\Field) := \varprojlim \cM(\cX), \]
    where $\cX$ ranges over the system of arithmetic varieties over $\Spec(\ZZ)$ with the data of an embedding of their function field into $\Field$, and $\cM$ is the Berkovich analytification functor.
\end{proposition}

\begin{proposition}(Remark \ref{remark_theta_is_homeo_trivial_valuation}) \label{intro:proposition_theta__relative_intro}
    In any characteristic, if either $\Field/k$ is finitely generated, or $\Field$ is countable, then there is a homeomorphism 
    \[\theta:\Theta_{\Field/k} \to \varprojlim X^{\an},\]
    where $\Theta_{\Field/k}$ is the space of pseudo-absolute values on $\Field$ that restrict to the trivial absolute value on $k$, and $X$ are normal projective varieties over $k$ together with an embedding of $k(X)$ to $F$ over $k$. Here $(-)^{\an}$ is the Berkovich analytification with respect to the trivial absolute value on $k$.
\end{proposition}

This description is useful, because using the results about Berkovich spaces over $\ZZ$ from \cite{lemanissier2020espaces} we easily get the following result (related to \cite[Lemma 3.1.1]{yuan2024adelic}).

\begin{corollary}(Corollary \ref{corollary_density_of_norms_on_a_countable_field} and Corollary \ref{corollary_density_of_norms_on_a_countable_field_relative_version})
    The set of absolute values on a countable field $\Field$ is dense in the space of pseudo-absolute values $\Theta_{\Field}$.

    Also, if $\Field/k$ is finitely generated (or both are countable), then the set of absolute values on $\Field$ that are trivial on $k$ is dense in $\Theta_{\Field/k}$.
\end{corollary}

We use this corollary to identify positive elements of the \textit{universal lattice of $\Field^{\times}$}, i.e., the divisible hull of the Grothendieck group of the lattice of finite subsets of $\Field^{\times}$ (with $+$ being the multiplication of sets and $\vee$ being the union). We denote this lattice by $\ULatQ(\Field)$, see Definition \ref{definition:universal_lattice} for more details. Let $|\cdot|$ be a pseudo-absolute value on $\Field$ and $\alpha = \bigjoin_{x\in I} \adiv(x) - \bigjoin_{y\in J} \adiv(y) \in \ULatQ(\Field)$. Using the notation $v(x) := -\log|x|$ we define $v(\alpha) = \max_{x \in I} \min_{y \in J} v(x/y)$. We call $\alpha \in \ULatQ(\Field)$ \textit{positive} if for all $v$ as above, $v(\alpha) \geqslant 0$. The following characterisation of positive elements is used to get equivalence of points (2) and (3) in Theorem \ref{theorem:intro:main}.

\begin{theorem}(Corollary \ref{corollary:tropical_divisors_positivity_conditions})
    For each $\alpha\in\ULatQ(F)$, the following conditions are equivalent:
\begin{enumerate}[label=(\roman*)]
    \item For every pseudo-absolute value $|\cdot|$ on $F$, $v(\alpha) \geqslant 0$.
    \item For every $\varepsilon\in\mathbb{R}_{>0}$, there is some $m\in\mathbb{N}$ such that $m\alpha\in\Gamma_{m\varepsilon}$.
\end{enumerate}
\end{theorem}

In this statement, $\Gamma_{\varepsilon} \subset \ULatQ(F)$ (for any $\varepsilon>0$) is the set of elements of the form
\[
    \bigmeet_{i=1}^p \adiv\left(\sum_{j=1}^q m_{i,j}a_j\right) - \bigmeet_{j=1}^q \adiv(a_j),
\]
where $a_1,\ldots,a_q\in F^\times$, $p,q\in\mathbb{N}^*$, and $(m_{i,j})_{\substack{1\leqslant i \leqslant p \\ 1 \leqslant j \leqslant q}} \in \mathbb{Z}^{pq}$ are such that for all $i\leqslant p$ we have $\sum_{j=1}^q |m_{i,j}| < 2^\varepsilon$. Moreover, if $F$ is countable, in (i) one can only consider absolute values. The significance of this theorem is that it gives a syntactic criterion for positivity (which is a semantic property).

Furthermore, we study the quotient $\LdivQ(\Field)$ of $\ULatQ(\Field)$ by the kernel ideal $I$ of the pairing $(v,\alpha) \mapsto v(\alpha)$. The following characterisations are provided:
\begin{itemize}
    \item if $\Field$ is countable of characteristic zero, then $\LdivQ(\Field)$ can be seen as a sublattice of the space of Arakelov divisors on $\Field$-models, see Corollary \ref{corollary_description_of_LDiv_Arakelov_divisors};
    \item if $\Field/k$ is a finitely generated extension, then $\LdivQ(\Field)$ divided by the ideal generated by $\LdivQ(k)$ is isomorphic to the space of b-divisors on the Riemann-Zariski space of $\Field/k$, see Corollary \ref{corollary_quotient_of_lattice_divisors_relative_case}.
\end{itemize}
The lattice structure on the space of Arakelov divisors on $\Field$-models from the first bullet has been introduced in \cite{MR1408559, ru2016birationalnevanlinnaconstantconsequences} in the language of (generalised or b-) Weil functions. We use a different description of this lattice from \cite{szachniewicz2023existential}. The second bullet has been discovered (in a slightly different language) also in \cite[Theorem 5.2]{kaveh2014note}. Interestingly, already Andr{\'e} Weil was already aware of this interpretation of (b-)divisors in the 1950s, see \cite[Theorem 13]{weil1951arithmetic}.

The \textit{positive functional} appearing in the statement of Theorem \ref{theorem:intro:main} is by definition a linear functional $l:\LdivQ(\Field) \to \RR$ which is non-negative on the positive cone $\{\alpha \geqslant 0\} \subset \LdivQ(\Field)$.

Let us elaborate on the meaning of point (5) in Theorem \ref{theorem:intro:main} (point (4) is similar). We introduce a locally compact topological space $\Omega^{\circ}_{\Field}$ which is defined as the set of \textit{pseudo-valuations} $v:\Field \to [-\infty,\infty]$, i.e., functions of the form $-c\log|\cdot|$ for $c \in \RR_{>0}$ and non $\{0,1,\infty\}$-valued pseudo-absolute values $|\cdot|$ on $\Field$. Thus, one can think of $\Omega^{\circ}_{\Field}$ as an additive version of $\cM(\Field)$ (see the difference between Equation (\ref{eq:zero}) and Equation (\ref{eq:one})). For technical reasons, we sometimes need to work with a slightly bigger space $\Omega_{\Field} \supset \Omega^{\circ}_{\Field}$, but in practice one can restrict attention to pseudo-valuations, see Definition \ref{definition_almost_compactification_of_the_space_of_pre_vasluations} and Remark \ref{remark_pre_valuations_measure_can_be_on_them}.
The following representation theorem follows from equivalence of (1) and (5) in Theorem \ref{theorem:intro:main}.
\begin{corollary} \label{corollary:intro:representation_theorem_GVFs}
    Let $\Field$ be any field equipped with a GVF structure. Then there exists a measure $\mu$ on $\Omega^{\circ}_{\Field}$ such that the heights on $\Field$ are given by the formula
    \begin{equation*}
        h(a_1, \dots, a_n) = \int_{\Omega^{\circ}_{\Field}} -\min(v(a_1),\dots,v(a_n)) d\mu(v).
    \end{equation*}
\end{corollary}
There are two differences between Corollary \ref{corollary:intro:representation_by_adelic_curves} and Corollary \ref{corollary:intro:representation_theorem_GVFs}: the latter works also for uncountable fields, and the ambient topological space $\Omega^{\circ}_{\Field}$ does not vary. In fact, the former corollary follows from the latter. In general, the flexibility of choosing $\Omega$ to be any measurable space in the definition of adelic curves is useful, as it allows one to compare different metrics over the same absolute value, see for example Harder-Narasimhan filtrations in \cite{Adelic_curves_1, Adelic_curves_3, liu2024arithmetic}.

After proving Theorem \ref{theorem:intro:main} we give some examples and focus on the relative situation when $K \subset \Field$ is an extension of globally valued fields. Using a disintegration-type result, we get a relative version of Corollary \ref{corollary:intro:representation_theorem_GVFs}, namely Corollary \ref{corollary_disintegration_of_a_GVF_measure} and Corollary \ref{corollary:GVFs_as_measure_on_Berkovich_spaces}. This is interesting in light of recent applications of families of measures on Berkovich spaces, for example in \cite{poineau2024dynamique1, poineau2024dynamique2}.

In Section \ref{section_finite_extensions} we prove uniqueness (up to scaling) of the GVF structure on classical global fields (Lemma \ref{lemma:unique_GVF_function_field_curve} and Lemma \ref{lemma:uniqueness_of_the_GVF_structure_on_number_fields}), and the following.

\begin{proposition}(Proposition \ref{proposition:uniqueness_of_invariant_Galois_extension}) \label{intro:proposition_GVF_extensions}
    Let $K$ be a GVF and consider a finite Galois field extension $K \subset \Field$. Then there is a unique GVF structure on $\Field$ extending the one on $K$ that is invariant under the action of $\Aut(\Field/K)$.
\end{proposition}

It is worth mentioning that using representations by adelic curves, the existence (but not uniqueness) follows from \cite[Section 3.4]{Adelic_curves_1}. Also, under an additional assumption on the GVF structure, this result was independently obtained by Antoine Sédillot in \cite[Propositions II.5.3.1 and II.5.4.1]{SédillotPhD2024}.

\begin{acknowledgements*}
    We would like to thank Leo Gitin, Nuno Hultberg, Konstantinos Kartas, Klaus Künnemann, Jonathan Pila, Jérôme Poineau, Thomas Scanlon and the participants of logic seminars in Oxford and Wrocław for helpful conversations and remarks. The second and fourth authors would like to thank the organizers and the participants of the 2024 Students' Conference on Non-Archimedean, Tropical and Arakelov Geometry for their interest in this project. There we learned that variants of Proposition \ref{intro:proposition_theta_intro} and Proposition \ref{intro:proposition_GVF_extensions} were independently obtained in \cite{SédillotPhD2024}. We warmly thank Antoine Sédillot for informing us about that.
    This material is based upon work supported by the National Science Foundation under Grant No. DMS-2424441, and by the IAS School of Mathematics.
\end{acknowledgements*}

\section{Valuations}

Let $\Field$ be a field. We will later need an enlargement of the space of valuations on $\Field$ to be able to consider Archimedean valuations even if $\Field$ has cardinality bigger than the continuum. That is why we introduce the following definition (partially appearing already in \cite{weil1951arithmetic}).

\begin{definition} \label{definition_pre_norms_pre_valuations}~
    \begin{itemize}
        \item A \textit{pseudo-absolute value} $|\cdot|$ on $\Field$ is a function $|\cdot|:\Field \to \RR_{>0} \cup \{0, \infty\}$ satisfying $|x+y| \leqslant |x|+|y|, |0|=0, |1|=1$, and $|xy|=|x||y|$ whenever the product is defined, i.e., when $(|x|,|y|)\notin\{(0,\infty),(\infty,0)\}$.
        \item We call a pseudo-absolute value \textit{trivial} if it attains values only from $\{0, 1, \infty\}$. Otherwise, we call it a non-trivial pseudo-absolute value. 
        \item A \textit{pseudo-valuation} $v$ on $\Field$ is a map $v:\Field \to \RR \cup \{\infty,-\infty\}$ attaining at least one value outside of the set $\{-\infty, 0, \infty\}$ and such that $v(0)=\infty, v(1)=0, v(2) \neq - \infty,  v(x+y) \geqslant \min(v(x),v(y)) + \min(v(2),0)$, and $v(xy)=v(x)+v(y)$ whenever the sum is defined, i.e., when $(v(x),v(y))\notin\{(\infty,-\infty),(-\infty,\infty)\}$. 
        \item An \textit{abstract pseudo-valuation} $v$ on $\Field$ is a map $v:\Field \to \Gamma \cup \{\infty, - \infty\}$, where $\Gamma$ is a divisible (totally) ordered abelian group, which satisfies the above axioms of a pseudo-valuation (but with operations interpreted in $\Gamma$). Later we write $\{ \pm \infty \}$ for $\{\infty, - \infty\}$.
    \end{itemize}
\end{definition}

\begin{lemma} \label{lemma_quotient_and_restricted_abstract_prevaluation}
    Let $v:\Field \to \Gamma \cup \{\infty,-\infty\}$ be an abstract pseudo-valuation on $\Field$ and let $\Delta \subset \Gamma$ be a convex subgroup, i.e., a subgroup satisfying $a,b\in\Delta \Rightarrow [a,b]\subset \Delta$, where the interval $[a,b]$ is defined as $\{ c\in\Gamma, \,:\, a\leqslant c\leqslant b\}$. Then, the quotient $\Gamma/\Delta$ is also a divisible (totally) ordered abelian group, where we declare $a+\Delta \leqslant b+\Delta$ if and only if $a-b \in \Delta$ or $ a\leqslant b$. Assume that $v(2) \not< \Delta$. Consider the functions $u(a) = v(a) \mod \Delta$ and
    \[
      w(a) =
      \begin{cases}
        \infty & \text{if $v(a) > \Delta$,} \\
        v(a) & \text{if $v(a) \in \Delta$,} \\
        -\infty & \text{if $v(a) < \Delta$.}
      \end{cases}
    \]
    Then $u,w$ satisfy the axioms of abstract pseudo-valuations on $\Field$, possibly without the axiom asserting that they attain at least one value outside of the set $\{-\infty, 0, \infty\}$.
\end{lemma}
\begin{proof}
    Since the quotient map $\Gamma \rightarrow \Gamma/\Delta$ is an increasing group morphism, the conditions $v(x+y) \geqslant \min(v(x),v(y)) + \min(v(2),0)$ and $v(xy)=v(x)+v(y)$ are conserved for $u$. The conditions $u(0)=\infty$, $u(1)=0$ and $u(2) \neq -\infty$ are immediately verified.
    
    For $w$, the conditions $w(0) = \infty$, $w(1) = 0$ and $w(2)\neq -\infty$ are immediate, and $w(2) \neq -\infty$ comes from the assumption that $v(2) \not< \Delta$. The remaining two conditions are proven by isolating the cases $v(x) < \Delta$, $v(x)\in\Delta$ and $v(x)>\Delta$, and doing the same for $y$.
\end{proof}

We drop the prefix ``pseudo'' in Definition \ref{definition_pre_norms_pre_valuations}, if the map does not attain $\pm\infty$ at non-zero elements. Bear in mind that this gives a slightly non-standard definition of a \textit{valuation} on $\Field$, as we exclude the $\{0,1\}$-valued ones from the usual definition, and include some that do not satisfy the ultrametric inequality without the correction $\min(v(2),0)$. The set of valuations on $\Field$ is denoted by $\Val_F$ and the set of absolute values is denoted by $M_{\Field}$. We topologise those spaces with the pointwise convergence topology. Note that for a pseudo-absolute value $|\cdot|$ on $F$, the set of elements $\cO := \{|x|<\infty\}$ is a valuation ring with the maximal ideal $\mathfrak{m}$ consisting of elements with $|x|=0$. The quotient $K=\cO/\mathfrak{m}$ is then a field with an absolute value. Thus we get an induced surjective map $p:F \to K \cup \{\infty\}$ where the preimage of infinity is the complement of a valuation ring. Such a map is called a place. Summarising, the following holds.

\begin{remark}\label{remark_char_of_prevaluations_and_prenorms}
    \indent
    \begin{enumerate}
        \item If $|\cdot|$ is a pseudo-absolute value on $\Field$, there exists a place~{$p:\Field\to K\cup\{\infty\}$} into a field $K$ with an absolute value $|\cdot|_K$ such that $|x| = |p(x)|_K$, with the convention that $|\infty|_K=\infty$.
        \item If $|\cdot|$ is a non-trivial pseudo-absolute value on $\Field$ and $c \in \RR_{>0}$, then $v(x) = -c\log|x|$ defines a pseudo-valuation on $\Field$. Note that here we extend the natural logarithm function to $[0,\infty]$ by continuity.
        \item If $v$ is a pseudo-valuation on $\Field$, then there is a non-trivial pseudo-absolute value $|\cdot|$ on $\Field$ and a constant $c \in \RR_{>0}$ such that $v(x) = -c\log|x|$ for all $x \in \Field$.
    \end{enumerate}
\end{remark}
\begin{proof}
    The proof of (1) is outlined above. For (2) note that by (1) we can assume that $|\cdot|$ is an absolute value. By passing to the completion, without loss of generality $\Field$ is complete with respect to $|\cdot|$. Then either $|\cdot|$ is a non-Archimedean absolute value on $\Field$ in which case the statement follows, or (by the classification of Archimedean normed complete fields) $(\Field, |\cdot|) \subset (\CC, \|\cdot\|^{t})$, where $\|\cdot\|$ is the standard Euclidean absolute value on $\CC$ and $t \in (0,1]$. Direct inspection of the statement on $(\CC, \|\cdot\|^{t})$ finishes the proof of (2). 

    To see (3), put $a = -\min(v(2),0)$. The case $a = 0$ is clear, so assume that $a>0$. We define $|x| = 2^{-v(x)/a}$ and check that $|\cdot|$ is a pseudo-absolute value. Here, we choose a normalization of $|\cdot|$ so that the condition $v(x+y) \geqslant \min(v(x),v(y)) - a$ becomes $|x+y| \leqslant 2 \max(|x|,|y|)$, which implies by induction that
    \[ \Bigl| \sum_{i=1}^{2^n} x_i \Bigr| \leqslant 2^n \max_i |x_i|. \]
    
    In particular $|2^n| \leqslant 2^n$ and for any integer $N$ we get $|N| \leqslant 2N$ because if we write $N=2^a+M$, where $M < 2^a$, then
    \[
        |N| \leqslant 2\max(|2^a|,|M|) \leqslant 2\max(2^a,2M) = \max(2^{a+1},4M) \leqslant 2N,
    \]
    where the last inequality follows from the fact that $2^a\leqslant N$ and $2M\leqslant N$.
    Replacing $N$ with $N^k$ and taking limit over $k \to \infty$ we get $|N| \leqslant N$. 
    Hence
    \[ |x+y|^{2^n} = |(x+y)^{2^n}| \leqslant 2^{n+1} \max_{0 \leqslant i \leqslant 2^n} \binom{2^n}{i} |x|^i |y|^{2^n-i} \leqslant 2^{n+1} (|x|+|y|)^{2^n}. \]
    Taking $2^n$'th roots and the limit $n \to \infty$ we get $|x+y| \leqslant |x|+|y|$ which finishes the proof.
\end{proof}

We can drop the prefix ``pseudo'' everywhere in the above remark and points two and three remain true. Note that non-trivial pseudo-absolute values and pseudo-valuations represent almost the same data, the difference being that one can multiply pseudo-valuations by an arbitrary positive number, while raising a non-trivial pseudo-absolute value to a positive power does not always yield a pseudo-absolute value.

\begin{definition} \label{definition_almost_compactification_of_the_space_of_pre_vasluations}
    Let $A$ be the set of functions $v:\Field \to [-\infty, \infty]$ satisfying
    \[ v(xy)=v(x)+v(y), v(0)=\infty, v(1)=0,  v(x+y) \geqslant \min(v(x),v(y)) + \min(v(2),0), \]
    and $v(2) \neq -\infty$. 
    Equip $[-\infty, \infty]^F$ with the product topology.
    We define the space $\Omega_{\Field}$ as the closure of $A$ without the set of $\{-\infty,0,\infty\}$-valued functions on $F$, i.e., as $\overline{A} \setminus \{-\infty,0,\infty\}^F$. 
    It is an open subset of the compact space $\overline{A} \subset [-\infty,\infty]^{\Field}$. Hence, it is a locally compact, Hausdorff topological space. 
    Moreover, we denote the subset $\Omega^{\circ}_{\Field} := \{v\in\Omega_{\Field}\,:\,v(2)\neq -\infty\} \subset \Omega_{\Field}$ which equivalently can be defined as $A \setminus \{-\infty,0,\infty\}^{\Field}$ and coincides with the set of pseudo-valuations.
\end{definition}

\begin{example} \label{example_an_element_of_Omega_that_is_not_a_pre_absolute_value}
    Consider $\Field=\QQ(x)$ and let $|\cdot|_n$ be the absolute value on $\Field$ coming from the embedding $\QQ(x) \subset \CC$ sending $x$ to $1+\frac{\pi}{n}$. Then for $v_n = -n \log |\cdot|_n$ we have $\lim_n v_n(x) = -\pi$ and $\lim_n v_n(2) = -\infty$. Hence any limit point $v$ of $(v_n)_{n \in \NN}$ is in $\Omega_\Field$, but is not a pseudo-valuation. 
\end{example}

We call an element $v \in \OF$ non-Archimedean if $v(2) \geqslant 0$ and Archimedean otherwise. Not every element $v \in \Omega_{\Field}$ is a non-trivial pseudo-valuation on $\Field$, however, the set of non-trivial pseudo-valuations is dense in $\Omega_{\Field}$. We choose to work with $\Omega_{\Field}$ instead of the set of non-trivial pseudo-valuations because of nicer topological properties, see Lemma \ref{lemma_Hausdorff_compact} and the discussion following it.

Note that if $\Field$ is countable, then
\[\Val_\Field = \bigcap_{a \in \Field^{\times}} \{ v \in \Omega_\Field : v(a) \neq \pm \infty \} \]
is a Borel subset of $\Omega_\Field$ and in that case we equip $\Val_\Field$ with the $\sigma$-algebra of its Borel subsets. Denote by $\Theta_{\Field}$ the space of pseudo-absolute values equipped with pointwise convergence topology - note that it is a compact Hausdorff space. Overall, we get the following diagram of spaces associated to $\Field$
\[\begin{tikzcd}
	{M_{\Field}} & {\Theta_{\Field}} \\
	{\Val_{\Field}} & {\Omega_{\Field}}
	\arrow[hook, from=1-1, to=1-2]
	\arrow[dashed, from=1-1, to=2-1]
	\arrow[dashed, from=1-2, to=2-2]
	\arrow[hook, from=2-1, to=2-2]
\end{tikzcd}\]
where the dashed arrows are $-\log$'s defined only on non-trivial (pseudo-)absolute values.

Our goal now is to prove Corollary \ref{corollary_density_of_norms_on_a_countable_field} and Corollary \ref{corollary_density_of_norms_on_a_countable_field_relative_version}. Since these corollaries hold under different assumptions, let us restrict our attention to characteristic zero $\Field$ until Remark \ref{remark_theta_is_homeo_trivial_valuation}. Call a normal, projective, generically smooth, integral scheme $\cX \to \Spec(\ZZ)$ an \textit{$\Field$-submodel} if the function field of $\cX$ is equipped with an embedding into $\Field$. 
The system of $\Field$-models forms a category where morphisms from $f':\Spec(\Field) \to \cX'$ to $f:\Spec(\Field) \to \cX$ are morphisms of schemes $h:\cX' \to \cX$ such that $h \circ f' = f$.
Recall that for a point $p\in\cX$, $\kappa(p)$ denotes the residue field of $\cX$ at $p$. For a scheme $\cX$ define its Berkovich spectrum as
\[ \cM(\cX) := \{ (p, |\cdot|) \textnormal{ with } p \in \cX \textnormal{ and } |\cdot| \textnormal{ is an absolute value on } \kappa(p) \}. \]
For a point $x = (p, |\cdot|)$ if $f \in \cO_{\cX}(U)$ where $p \in U \subset \cX$, then we write $|f(x)|$ for $|f(p)|$ where $f(p) \in \kappa(p)$ is the evaluation of $f$ at $p$. We equip $\cM(\cX)$ with the coarsest topology such that the projection $\pi:\cM(\cX) \to \cX$ which sends $(p,|\cdot|)$ to $p$ is continuous, and for every $f$ as above, the function $x \mapsto |f(x)|$ is continuous on $\pi^{-1}(U)$. This makes $\cM(\cX)$ a compact Hausdorff topological space, provided that $\cX$ is projective over $\ZZ$. 
For a thorough treatment of these spaces, see \cite{lemanissier2020espaces, Poineau2008LaDD}. 
We define
\[ \cM(\Field) := \varprojlim \cM(\cX) \]
to be the inverse limit over (the directed system of) all $\Field$-submodels of $\Field$, as a topological space. 
If $X$ is a proper scheme over a complete valued field, we write $X^{\an}$ for the Berkovich analytification with respect to the ambient absolute value, i.e., the subset of $\cM(X)$ where absolute values on residue fields extend the one on the base-field.

\begin{construction}
    Let $|\cdot| \in \Theta_{\Field}$. Denote by $\cO \subset \Field$ the valuation subring consisting of elements with $|x| < \infty$, by $\mathfrak{m}$ its maximal ideal, and by $k$ the residue field. Denote the induced absolute value on $k$ also by $|\cdot|$. Fix an $\Field$-submodel $\cX$. Consider the diagram
    \[\begin{tikzcd}
    	{\Spec(\Field)} & \cX \\
    	{\Spec(\cO)} & {\Spec(\ZZ)}
    	\arrow[from=1-2, to=2-2]
    	\arrow[from=1-1, to=2-1]
    	\arrow[from=2-1, to=2-2]
    	\arrow[from=1-1, to=1-2]
    	\arrow["\phi"{description}, dashed, from=2-1, to=1-2]
    \end{tikzcd}\]
    where the top horizontal arrow comes from the embedding of the function field of $\cX$ into $\Field$ and the diagonal arrow exists by projectivity of $\cX$ over $\ZZ$. Let $p = \phi(\mathfrak{m})$ and consider the induced field extension $\kappa(p) \subset k$. We can restrict $|\cdot|$ to $\kappa(p)$ to get an absolute value inducing a point $x \in \cM(\cX)$. Summarising, we constructed a map
    \[ \theta:\Theta_{\Field} \to \cM(\Field). \]
\end{construction}

\begin{proposition} \label{proposition_theta_is_a_homeo}
    The map $\theta$ is a homeomorphism.
\end{proposition}
\begin{proof}
    First let us check that for an $\Field$-submodel $\cX$ the composed map 
    \[ \Theta_{\Field} \to \cM(\Field) \to \cM(\cX) \to \cX \]
    is continuous with respect to the Zariski topology on the codomain. Let $\Spec(A) \subset \cX$ be an open affine subset and identify $A$ with its image under the induced embedding $A \hookrightarrow \Field$. If we denote by $p:\Theta_{\Field} \to \cX$ the composition above, then 
    \[ p^{-1}(\Spec(A)) = \{ |\cdot| \in \Theta_{\Field} : A \subset \{|x|<\infty\} \}. \]
    However, since $A$ is finitely generated over $\ZZ$, if we pick its generators $a_1, \dots, a_n$, then this set can be written as $\{ |\cdot| \in \Theta_{\Field} : \max_i |a_i|<\infty \}$ which is an open subset of $\Theta_{\Field}$. Thus $p$ is continuous. Moreover, if $f\in\mathcal{O}_\mathcal{X}(U)$, the map $|\cdot|_F\in \Theta_F \mapsto |f(\theta(|\cdot|_{F}))| = |f|_{F}$ is continuous since $\Theta_F$ is equipped with the topology of pointwise convergence. Hence, the continuity of $\theta$ follows, by the definition of topology on $\cM(\cX)$. Since both domain and codomain of $\theta$ are compact Hausdorff spaces, it is enough to show that $\theta$ is a bijection. The classical Riemann-Zariski space construction shows that valuation subrings of $\Field$ are in bijection with the inverse limit of the system of $\Field$-submodels, see e.g. \cite[Corollary 3.4.7]{temkin2011relative}. More precisely, if $\cO \subset \Field$ is a valuation ring, then there exists a unique family of points $x(i) \in \cX_i$, where $(\cX_i)_{i \in I}$ is the system of $\Field$-submodels, such that $\cO = \varinjlim \cO_{\cX_i, x(i)}$. Then absolute values on the residue field $k$ of $\cO$ correspond to consistent families of absolute values on $\kappa(x(i))$ for $i \in I$, which finishes the proof.
\end{proof}

\begin{proposition} \label{proposition_density_of_generic_fiber_Berkovich_analitification_of_an_arithmetic_variety}
    Let $\cX$ be an $\Field$-submodel. Then the generic fiber of $q:\cM(\cX) \to \cX$ is dense in $\cM(\cX)$. 
\end{proposition}
\begin{proof}
    First note that it is enough to prove that for an open $\cU \subset \cX$ the set $q^{-1}(\cU)$ is dense. Indeed, this is because $\cM(\cX)$ is compact and Hausdorff, hence a Baire space, and $\cX$ has a countable basis of open subsets, whose intersection is exactly the generic point of $\cX$.

    Fix open non-empty subsets $\cU \subset \cX$ and $V \subset \cM(\cX)$. We need to show that $V \cap q^{-1}(\cU) \neq \varnothing$. By \cite[Proposition 6.4.1 and Proposition 6.6.10]{lemanissier2020espaces} the map $\cM(\cX) \to \cM(\ZZ)$ is open. Hence, by the description of $\cM(\ZZ)$ e.g. in \cite[Corollaire 3.1.12]{Poineau2008LaDD}, there is $x \in V$ with its image $x|_{\QQ}$ in $\cM(\ZZ)$ being an absolute value on $\QQ$. Since the problem is invariant under raising an absolute value to a positive power, without loss of generality $x|_{\QQ}$ is either a $p$-adic absolute value on $\QQ$ (for some prime $p$), or the Euclidean one. Let $K$ be the completion of $\QQ$ with respect to this absolute value.
    
    Let $X$ be the generic fiber of $\cX \to \Spec(\ZZ)$, so it is a smooth, projective scheme over $\QQ$ with the function field $\Field$. Denote by $U$ the intersection of $\cU$ with $X$ and assume that the point $x \in V$ corresponds to an absolute value $|\cdot|$ on $\kappa(x_0)$ for $x_0 \in X$. We consider two cases:
    
    (1) $x|_{\QQ}$ is non-Archimedean. In that case $K=\QQ_p$ and we consider the Berkovich analytification $X_K^{\an}$ of $X_K$ with projection $\rho:X_K^{\an} \to X_K$. We get the pullback diagram
    \[\begin{tikzcd}
    	{X_K^{\an}} & {\cM(\cX)} \\
       	{x|_{\QQ}} & {\cM(\ZZ)}
       	\arrow["j", hook, from=1-1, to=1-2]
       	\arrow[from=1-1, to=2-1]
       	\arrow[from=1-2, to=2-2]
       	\arrow[hook, from=2-1, to=2-2]
    \end{tikzcd}\]
    As $U_K \subset X_K$ is an everywhere dense open, by \cite[Corollary 3.4.5]{Berkovich_Spectral_Theory} the subset $\rho^{-1}(U_K)$ is everywhere dense in $X_K^{\an}$. Thus, as $V$ intersects the fiber over $x|_{\QQ}$ witnessed by $x \in V$, the set $V \cap \rho^{-1}(U_K)$ is non-empty. But $\rho^{-1}(U_K) = j^{-1}(q^{-1}(\cU))$, so if $y \in \rho^{-1}(U_K)$, then $j(y) \in V \cap q^{-1}(\cU)$, which finishes the proof in this case.
        
    (2) $x|_{\QQ}$ is Archimedean. In that case $K = \RR$ and we look at the complex analytification $X_{\CC}^{\an}$ which is a smooth complex variety with the conjugation map $c$. Then the fiber of $\cM(\cX)$ over $x|_{\QQ}$ is $X_{\CC}^{\an}/c$, and it is enough to prove that non-empty Zariski open subsets of $X$ yield everywhere dense open subsets of $X_{\CC}^{\an}$. This is classical, and can be easily deduced e.g. from the resolution of singularities.
\end{proof}

We use the above proposition to get the following corollary.

\begin{corollary} \label{corollary_density_of_norms_on_a_countable_field}
    Let $\Field$ be a countable field of characteristic zero. The set $M_{\Field}$ is a dense subset of $\Theta_{\Field}$. In particular, the set $\Val_\Field$ is a dense subset of $\Omega_{\Field}$.
\end{corollary}
\begin{proof}
    Note that $M_{\Field} = \bigcap_{a \in \Field^{\times}} \{ |\cdot| \in \Theta_{\Field} : |a|<\infty \}$. By the fact that $\Theta_{\Field}$ is a Baire space, and $\Field$ is countable, it is enough to show that for $a \in \Field^{\times}$, the set $V = \{ |\cdot| \in \Theta_{\Field} : |a|<\infty \}$ is dense.

    Using the homeomorphism from Proposition \ref{proposition_theta_is_a_homeo}, it is enough to check that for an $\Field$-submodel $\cX$ and an open subset $U \subset \cM(\cX)$, if we denote by $s:\Theta_{\Field} \to \cM(\cX)$ the natural projection, $V \cap s^{-1}(U) \neq \varnothing$. By possibly replacing the $\Field$-submodel $\cX$ by one over it, we may assume that there is an open subscheme $\Spec(A) \subset \cX$ such that $a \in A \subset \Field$. By Proposition \ref{proposition_density_of_generic_fiber_Berkovich_analitification_of_an_arithmetic_variety}, there is a point $x \in U \subset \cM(\cX)$ such that its first coordinate (if we write is as a pair) is in $\Spec(A) \subset \cX$.

    Now, if $|\cdot| \in \Theta_{\Field}$ is such that $s(|\cdot|)=x$, then $|\cdot| \in V \cap s^{-1}(U)$, which gives the desired non-emptiness. Since maps in the system of Berkovich analytifications of $\Field$-submodels are surjective, such a pseudo-absolute value $|\cdot|$ exists, which finishes the proof.
\end{proof}

\begin{remark} \label{remark_theta_is_homeo_trivial_valuation}
    One can perform the above constructions also in the relative context, where we fix a base-field $k \subset \Field$ and $\Field$ is either finitely generated over $k$ of arbitrary characteristic, or countable. Indeed, let $\Theta_{\Field/k} \subset \Theta_{\Field}$ be the set of pseudo-absolute values that are $1$ on $k^{\times}$ and let $M_{\Field/k}$ be the set of absolute values on $\Field$ that are trivial on $k$. We then get (using the same proof as in Proposition \ref{proposition_theta_is_a_homeo})
    \[ \Theta_{\Field/k} \simeq \varprojlim X^{\an}, \]
    where the inverse limit is taken over normal projective varieties $X$ over $k$ together with an embedding of $k(X)$ into $\Field$ over $k$. The Berkovich analytification $X^{\an}$ is with respect to the trivial absolute value on $k$. The generic fiber of the natural projection $\pi:X^{\an} \to X$ is dense, by \cite[Section 3.5]{Berkovich_Spectral_Theory}.
    
    In the case where $\Field$ is countable, we can copy the proof of Corollary \ref{corollary_density_of_norms_on_a_countable_field} to get density of $M_{\Field / k}$ in $\Theta_{\Field / k}$. In the case where $\Field$ is finitely generated over $k$, we can write
    \[ \Theta_{\Field/k} \simeq \varprojlim_{k(X) \simeq \Field} X^{\an}, \]
    where in the limit we only take varieties over $k$ with an isomorphism of their function field with $\Field$ (over $k$). If $X' \to X$ is a birational map in the system, we get a commutative diagram
    \[\begin{tikzcd}
    	{X'^{\an}} & {X^{\an}} \\
    	{M_{\Field/k}} & {M_{\Field/k}}
    	\arrow[from=1-1, to=1-2]
    	\arrow[hook, from=2-1, to=1-1]
    	\arrow["\id"', from=2-1, to=2-2]
    	\arrow[hook, from=2-2, to=1-2]
    \end{tikzcd}\]
    where the vertical maps come from isomorphisms $k(X), k(X') \simeq \Field$. Taking inverse limits, we also get density in that case. Summarising, we get the following.
\end{remark}

\begin{corollary} \label{corollary_density_of_norms_on_a_countable_field_relative_version}
    Let $k \subset \Field$ be either a finitely generated extension, or an extension of countable fields. Then the set $M_{\Field / k}$ is a dense subset of $\Theta_{\Field / k}$.
\end{corollary}

Let us come back to the context of a general field $\Field$. Define $\ov{\Omega}_{\Field}$ to be the quotient of $\Omega_{\Field}$ by the scaling action of $\RR_{>0}$ (with the quotient topology) and call two elements of $\Omega_{\Field}$ \textit{equivalent} if they lie in a single orbit of this action. We use the symbol $\sim$ to denote this equivalence relation. For $a \in \Field^{\times}$ we define sets
\[ \Omega_{\Field}(a) := \{v \in \Omega_{\Field} : v(a) = 1\}, \]
and $\ov{\Omega}_{\Field}(a) := \pi(\Omega_{\Field}(a))$ for the projection $\pi:\Omega_{\Field} \to \ov{\Omega}_{\Field}$. Note that $\ov{\Omega}_{\Field}(a)$ is the set of classes of $v \in \Omega_{\Field}$ with $0 < v(a) < \infty$. Moreover, $\Omega_{\Field}(a)$ is compact and $\pi$ defines a bijection between $\Omega_{\Field}(a)$ and $\ov{\Omega}_{\Field}(a)$. We also have
\[ \ov{\Omega}_{\Field} = \bigcup_{a \in \Field^{\times}} \ov{\Omega}_{\Field}(a). \]

\begin{remark}
    In general $\ov{\Omega}_{\Field}$ is not Hausdorff. For example, if $\Field=\QQ(x)$, then the relation of being equivalent in $\Omega_{\Field}$ is not closed. Indeed, consider the elements $0, x, x-p^{-1}$ in $\Field$ and let $|\cdot|$ be an appropriate power of a Gauss absolute value of suitable radius (with respect to the $p$-adic absolute value on $\QQ$) so that we have
    \[ |x| = 2, |x-p^{-1}| = 2, |p^{-1}| = 1+\frac{1}{n}. \]
    Let $v_n$ be the corresponding valuation on $\Field$ and note that $v_n \sim n \cdot v_n$. By taking limits (possibly we need to pass to a subnet in $[-\infty, \infty]^{\Field}$), we see that $v_n$ converges to a pseudo-valuation $v$ where the above triangle has sides $(2,2,1)$, while $n \cdot v_n$ converges to a pseudo-valuation $w$ where the same triangle has sides $(\infty, \infty, e)$. Note that these values are not contained in $\{-\infty, 0, \infty\}$, so $v, w \in \Omega_{\Field}$. Moreover, $v$ is not equivalent to $w$.
\end{remark}
Even though the full space $\ov{\Omega}_{\Field}$ may not be Hausdorff, we have the following lemma.
\begin{lemma} \label{lemma_Hausdorff_compact}
    The sets $\ov{\Omega}_{\Field}(a)$ are Hausdorff topological subspaces, for $a \in \Field^{\times}$.
\end{lemma}
\begin{proof}
    Consider $\pi^{-1}(\ov{\Omega}_{\Field}(a)) = \{v \in \Omega_{\Field}: 0 < v(a) < \infty\}$. We show that the relation of being equivalent is closed when restricted to this subset. Let $f_b:\pi^{-1}(\ov{\Omega}_{\Field}(a)) \to [-\infty, \infty]$ be a function given by $f_b(v) = \frac{v(b)}{v(a)}$, for $b \in \Field$. Note that these are continuous functions and for $v, w \in \pi^{-1}(\ov{\Omega}_{\Field}(a))$ we have $v \sim w$ if and only if for all $b \in \Field$ the equality $f_b(v) = f_b(w)$ holds, which is a closed condition.
\end{proof}
In particular, because restricted quotients $\pi:\Omega_{\Field}(a) \to \ov{\Omega}_{\Field}(a)$ are continuous, they are homeomorphisms between compact Hausdorff spaces.

Consider an extension of fields $F \subset F'$. The following shows that the Chevalley's extension theorem holds for pseudo-absolute values.

\begin{lemma} \label{Lemma_extending_pseudo_absolute_values}
    Let $v \in \Omega_{F}$. Then, there exists $v' \in \Omega_{\Field'}$ such that $v'|_{F} = v$.
\end{lemma}
\begin{proof}
    Let $a \in \Field$ be such that $0<v(a)<\infty$. Without loss of generality, we can assume that $v(a) = 1$. Consider the restriction map $\OFp(a) \to \OF(a)$. It is a continuous map between compact Hausdorff spaces, so to prove it is surjective it is enough to prove it has dense image. Thus, we can assume that $v(2) \neq -\infty$. By Remark \ref{remark_char_of_prevaluations_and_prenorms} we can assume that $v = -\log|\cdot|$ for some non-trivial pseudo-absolute value $|\cdot|$.

    Let $\cO \subset \Field$ be the valuation ring of elements with $|x| \neq \infty$ and let $k$ and $\mathfrak{m}$ be its residue field and maximal ideal respectively. By Chevalley's theorem on extending valuations, there is a valuation ring $\cO' \subset \Field'$ with maximal ideal $\mathfrak{m}'$ and residue field $k'$ such that $\cO' \cap \Field = \cO$ and $\mathfrak{m}' \cap \cO = \mathfrak{m}$. The non-trivial pseudo-absolute value $|\cdot|$ induces an absolute value $|\cdot|$ on $k$ and we are reduced to proving that the corresponding valuation extends to a pseudo-valuation on $k'$.

    Let $K$ be a complete, algebraically closed, valued field, such that the absolute value $|\cdot|$ on $k$ comes from an embedding $k \subset K$. Such exists, as we can extend the induced absolute value on the completion of $k$ to its algebraic closure. Take an ultrapower $K^*$ of $K$ such that $k'$ embeds into $K^*$ over $k$. The ultrapower of the valuation on $K$ gives an abstract pseudo-valuation $v^*$ on $K^*$ with the value group $\Gamma$. Let $\Delta \subset \Gamma$ be the convex hull of $\RR$ and $\Delta_0 \subset \Delta$ be the convex subgroup of infinitesimal elements. We use Lemma \ref{lemma_quotient_and_restricted_abstract_prevaluation} twice: first restricting $v^*$ to $\Delta$, and next quotienting $\Delta$ by $\Delta_0$. Identifying $\Delta/\Delta_0$ with $\mathbb{R}$, we get a (non-abstract) pseudo-valuation $u$ on $K^*$ extending $v$ on $k$. By Remark \ref{remark_char_of_prevaluations_and_prenorms}, we can write $u = -\log|\cdot|_u$ for some pseudo-absolute value $|\cdot|_u$ that extends $|\cdot|$ on $k$ (without loss of generality we can take $c=1$ for the constant from Remark \ref{remark_char_of_prevaluations_and_prenorms}). By restricting $|\cdot|_u$ to $k'$, we are done.
\end{proof}

\section{Seminorms}

In this section we introduce a submultiplicative variant of (pseudo-)absolute values and show how to obtain (pseudo-)absolute values from them.

\begin{definition}
    A \emph{pseudo-seminorm} $|\cdot|$ on a ring $R$ is a map $|\cdot| : R \rightarrow \mathbb{R}_{\geqslant 0}\cup\{+\infty\}$ satisfying $|x+y|\leqslant |x|+|y|$, $|xy|\leq|x|\cdot|y|$, $|0|=0$, $|\pm 1| = 1$, with $0\cdot\infty$ undefined. A pseudo-seminorm $|\cdot|$ is called \emph{power-multiplicative} if $|x^n|=|x|^n$ for all $x\in R$ and $n\in\mathbb{N}$. It is called a pseudo-absolute value if $|xy|=|x||y|$ for all $x,y\in R$. We define a partial order on the set of pseudo-seminorms of $F$ by
    \[
        |\cdot|_1\leqslant |\cdot|_2 \Leftrightarrow \forall x\in F,\, |x|_1\leqslant |x|_2.
    \]
\end{definition}
Note that a pseudo-seminorm on a ring $R$ is simply the data of a subring $R'\subset R$ (characterized by $R' = \{ x\in R\,:\, |x|\neq \infty\}$) and a seminorm $|\cdot|_0$ on $R'$ such that $|1|_0=1$. When $R$ is a field, the notion of pseudo-absolute value coincides with the one of Definition \ref{definition_pre_norms_pre_valuations}.

\begin{lemma} \label{lemma:minimal_preseminorm_prenorm_field}
    Let $|\cdot|$ be a minimal pseudo-seminorm on a field $R$. Then $|\cdot|$ is a pseudo-absolute value.
\end{lemma}
\begin{proof}
    Let $R' = \{|x|<\infty\}$. Note that if there was a seminorm on $R'$ smaller than the restriction of $|\cdot|$ to $R'$, then one could extend it to a pseudo-seminorm on $R$ by taking the value $\infty$ at $R \setminus R'$, which would contradict minimality of $|\cdot|$. Hence, the restriction of $|\cdot|$ is a minimal seminorm on $R'$. Consider its spectral radius norm
    \[
        \rho(x) = \left\{
            \begin{array}{ll}
                \lim\limits_{n\rightarrow +\infty} |x^n|^{1/n} &\text{ if } |x|\neq \infty \\
                \infty &\text{ if } |x| = \infty.
            \end{array}
        \right.
    \]
    It is a classical result (see \cite[Corollary 1.3.3]{Berkovich_Spectral_Theory}) that the spectral radius of a seminorm is a power-multiplicative seminorm. By minimality, we get $|\cdot|=\rho$, which shows that $|\cdot|$ is power-multiplicative on $R'$.

    Note that if $|x|=0$ for $x \in R'$, then for any $y \in R'$, we have $|xy| \leqslant |x||y| = 0$, so in fact $|xy|=|x||y|$ in this case. On the other hand, take $x \in R'$ with $|x| \neq 0$. Let us define the following function on $R'$:
    \[ p(z) = \inf_{n \in \mathbb{N}} \frac{|x^n z|}{|x|^n}. \]
    By the power-multiplicativity of $|\cdot|$ on $R'$, $p(1)=1$ and it is easy to check that $p$ is subadditive, and submultiplicative on $R'$. Moreover, looking at $n=0$ one gets $p \leqslant |\cdot|$ on $R'$. By minimality, this means that $p=|\cdot|$, so $|y|=p(y) \leqslant \frac{|xy|}{|x|}$, which proves that $|xy|=|x||y|$. Thus, $|\cdot|$ is an absolute value on $R'$.

    Let $\mathfrak{p}=\{|x|=0\} \subset R'$ be the prime kernel ideal of $|\cdot|$. By the Chevalley's extension theorem, there is a valuation subring $\mathcal{O} \subset R$ with a maximal ideal $\mathfrak{m}$ satisfying $R' \cap \mathfrak{m} = \mathfrak{p}$. By Lemma \ref{Lemma_extending_pseudo_absolute_values}, there is a pseudo-absolute value $v'$ on $\mathcal{O}/\mathfrak{m}$ extending the norm $|\cdot|$ on $R'/\mathfrak{p}$. Extending this pseudo-absolute value by $\infty$ on $R \setminus \mathcal{O}$, we get that there is a pseudo-absolute value $v$ on $R$ extending the seminorm $|\cdot|$ on $R'$. By minimality (as a pseudo-seminorm on $R$) of $|\cdot|$, we get that $|\cdot| = v$, i.e., that $|\cdot|$ is a pseudo-absolute value.
\end{proof}

\begin{corollary} \label{corollary:existence_smaller_prenorm}
    Let $|\cdot|$ be a pseudo-seminorm on $F$. Then, there exists a pseudo-absolute value $\lVert\cdot\rVert$, such that $\lVert\cdot\rVert\leqslant |\cdot|$.
\end{corollary}

\begin{proof}
    It is clear that if $S$ is a totally ordered set of pseudo-seminorms of $F$, the lower bound $|\cdot|_{\inf}$ defined by
    \[
        \forall x\in F,\, |x|_{\inf} = \inf_{|\cdot|_s\in S} |x|_s.
    \]
    is again a pseudo-seminorm. Thus, we conclude using Zorn's lemma and Lemma \ref{lemma:minimal_preseminorm_prenorm_field}.
\end{proof}

\begin{lemma} \label{lemma:existence_prenorm_ring_map}
    Let $f:R\rightarrow F$ be a ring homomorphism from a ring to a field, and let $|\cdot|$ be a pseudo-seminorm on $R$ such that
    \[
        \forall x\in R,\, f(x) = 1 \Rightarrow |x|\geqslant 1.
    \]
    Then, there exists a pseudo-absolute value $\lVert \cdot \rVert$ on $F$ such that
    \[
        \forall x\in R,\, \lVert f(x)\rVert \leqslant |x|.
    \]
\end{lemma}

\begin{proof}
    Define a function $\lVert\cdot\rVert_0 : F\rightarrow \mathbb{R}_{\geqslant 0}\cup\{+\infty\}$ by
    \[
        \lVert y\rVert_0 = \inf_{x\in f^{-1}(\{y\})} |x|.
    \]
    (thus $\lVert y\rVert_0 = +\infty$ if $y$ is not in the image of $f$).
    The condition $f(x) = 1 \Rightarrow |x|\geqslant 1$ implies that $\lVert 1\rVert_0=1$, hence $\lVert \cdot \rVert_0$ is a pseudo-seminorm. By Corollary \ref{corollary:existence_smaller_prenorm}, there exists a pseudo-absolute value $\lVert\cdot\rVert$ on $F$ such that $\lVert \cdot \rVert \leqslant \lVert\cdot\rVert_0$. By construction, it satisfies the lemma's conclusion.
\end{proof}

\section{Lattices}

Important objects appearing in this paper are certain universal lattices of $\Field^{\times}$. Before we define those, we recall the definitions of a lattice monoid and a lattice group.

\begin{definition}
A \emph{lattice monoid (group)} is a commutative monoid (group) $(\Gamma,+)$ equipped with a \textit{join} operation $\join$ which satisfies the following axioms.
\begin{itemize}
    \item Idempotence: $\forall x\in \Gamma,\; x\join x = x$
    \item Commutativity: $\forall x,y\in \Gamma,\; x\join y = y\join x$
    \item Associativity: $\forall x,y,z \in \Gamma,\; (x\join y)\join z = x\join (y\join z)$
    \item Translation-invariance: $\forall x,y,z \in \Gamma,\; (x+z)\join(y+z) = (x\join y)+z$
\end{itemize}
A morphism of lattice monoids (groups) is a map from a lattice monoid (group) to another which preserves addition and join. We call a lattice group divisible, if the ambient abelian group is divisible. Every lattice monoid has a natural partial order defined by
\[ x \geqslant y \iff x \vee y = x. \]
\end{definition}

It is easy to see that this partial order is compatible with the group law of the lattice group, and that the operations $\vee$ and $+$ are increasing with respect to this order.
Lattice groups have the following property which may not hold in lattice monoids (for example in the lattice monoid of formal joins in $\Field$, see Definition~\ref{definition_formal_joins}).
\begin{lemma}\label{lemma:lattice_group_multiplication_by_n}
    Let $\Gamma$ be a lattice group. Then for any $x, y \in \Gamma$ and $n \in \NN$ we have
    \[ n(x \vee y) = (nx) \vee (ny). \]
\end{lemma}
\begin{proof}
    By writing $x = x^{+} + x^{-}$, where $x^{+} = x \vee 0, x^{-} = x \wedge 0$, and similarly for $y$, we are reduced to showing that $(nx)^{\pm} = nx^{\pm}$. This follows from the calculation:
    \[ (nx)^+ + nx^+ = (nx) \vee 0 + \bigvee_{k=0}^n (kx) = \bigvee_{k=0}^{2n} (kx) = 2nx^{+}, \]
    and similarly for $(nx)^{-}$.
\end{proof}
The forgetful functor from the category of lattice groups to the category of abelian groups has a left adjoint, which is called the universal lattice group functor. Moreover, by Lemma \ref{lemma:lattice_group_multiplication_by_n}, the forgetful functor from the category of divisible lattice groups to the category of lattice groups also has a left adjoint. By composing those we get a functor from the category of abelian groups to the category of divisible lattice groups (left adjoint to the forgetful functor), which we call the universal divisible lattice group functor.  Let us describe how to construct these functors. 

\begin{definition}~\label{definition_formal_joins}
    Let $(G,+)$ be an abelian group. A \textit{formal join in $G$} is an expression of the form $\bigjoin_{x\in I} [x]$, where $I$ is a finite subset of $G$. The set of such expressions forms a lattice monoid with operations defined as follows:
    \[
    \bigjoin_{x\in I} x \join \bigjoin_{y\in J} y = \bigjoin_{z\in I\cup J} z,
    \]
    \[
    \bigjoin_{x\in I} x + \bigjoin_{y\in J} y = \bigjoin_{z\in I + J} z.
    \]
    In other words, the formal joins form a lattice monoid isomorphic to the lattice monoid of finite subsets of $G$ with join being the union and addition being the sum of sets.

    The Grothendieck group of this lattice monoid inherits its lattice group structure. We call it the universal lattice group generated by $G$. Its divisible hull satisfies the universal property of the universal divisible lattice group of $G$. Equivalently, one could also take the Grothendieck group of the lattice monoid of formal joins in $G_{\QQ} = G \otimes \QQ$. 
\end{definition}

Now, we define the lattice group of tropical polynomials, which we use to encode global structures.

\begin{definition}
    Let $x_1,\ldots,x_n$ be symbols for abstract variables. Elements of the universal lattice group generated by $\ZZ x_1 + \ldots + \ZZ x_n$ will be called $\ZZ$-tropical polynomials in $x_1,\ldots,x_n$. Then, if $G$ is a lattice group and $g_1,\ldots,g_n\in G$, by the universal property of the universal lattice groups generated by $\ZZ x_1 + \ldots + \ZZ x_n$, the evaluation map $\ZZ x_1 + \ldots + \ZZ x_n \rightarrow G$ which sends $x_i$ to $g_i$ lifts uniquely to a map from the lattice group of $\ZZ$-tropical polynomials in $x_1,\ldots,x_n$ to $G$. If $t$ is such a tropical polynomial, its image under this map will be called the \textit{evaluation} of $t$ at $(g_1,\ldots,g_n)$, and denoted $t(g_1,\ldots,g_n)$.

    Similarly, elements of the universal divisible lattice group generated by $\ZZ x_1 + \ldots + \ZZ x_n$ will be called $\QQ$-tropical polynomials in $x_1,\ldots,x_n$. Then, if $G$ is a divisible lattice group and $g_1,\ldots,g_n\in G$, the same construction yields an evaluation map from the set of $\QQ$-tropical polynomials in $x_1,\ldots,x_n$ to $G$.

    If $x = (x_1,\ldots,x_n)$ is a tuple of variables, we will often write $t(x)$ instead of $t$ to indicate that $t$ is a tropical polynomial in $x_1,\ldots,x_n$. Notice that this is coherent with the above defined evaluation map, as the evaluation of $t$ at $x_1,\ldots,x_n$ is equal to $t$.
\end{definition}

Another encoding of GVF structures will use a different lattice group called the universal lattice associated to a field, which we now define.

\begin{definition}\label{definition:universal_lattice}
    We define the \textit{universal lattice} of $\Field$, denoted $\ULatQ(F)$, as the universal divisible lattice group generated by $(\Field^{\times},\times)$. When $a\in F^{\times}$, we denote by $\adiv(a)$ the associated element of $\ULatQ(F)$, i.e., the formal join of the singleton $\{a\}$.
\end{definition}

In other words, elements of the universal lattice of $\Field^{\times}$ can be represented by differences
\[ \bigjoin_{x\in I} \adiv(x) - \bigjoin_{y\in J} \adiv(y), \]
for finite $I, J \subset \Field_{\QQ}^{\times}$. It follows that if $t$ is a $\QQ$-tropical polynomial and $a = (a_1, \dots, a_n)$ is a tuple of elements from $\Field^{\times}$, then $t(\adiv(a)) = t(\adiv(a_1), \dots, \adiv(a_n))$ can be written in the above form in $\ULatQ(\Field)$ with $I$ and $J$ consisting of multiples of rational powers of $a_i$'s (and hence in any divisible lattice group with a homomorphism $\adiv:\Field^{\times} \to \Gamma$). More generally, using the same argument the following holds.

\begin{lemma} \label{lemma_tropical_polynomials_are_differences_of_maxima}
    Let $t(x)$ be a $\QQ$-tropical polynomial with $x=(x_1, \dots, x_n)$. Then there are terms $\alpha_i(x) = \sum_{k=1}^{n} q_{ik} \cdot x_k$ and $\beta_j(x) = \sum_{k=1}^{n} r_{jk} \cdot x_k$ (with $q_{ik}, r_{jk} \in \QQ$) over some finite index sets $i \in I, j \in J$ such that in any divisible lattice group $\Gamma$, the equality
    \[ t(x) = \max_{i \in I} \alpha_i(x) - \max_{j \in J} \beta_j(x) = \max_{i \in I} \min_{j \in J} (\alpha_i(x)-\beta_j(x)) \]
    holds, as functions on $\Gamma^n$.
\end{lemma}

Using the above observations, we define a pairing between $\Omega_{\Field}$ and $\ULatQ(\Field)$.
\begin{definition}
    We define a pairing
    \[ \ev:\Omega_{\Field} \times \ULatQ(\Field) \to [-\infty, \infty], \]
    by sending a pair $(v, t(\adiv(a)))$ to $v \bigl( t(\adiv(a)) \bigr) := t(v(a))$, where $a$ is a tuple. More precisely, this formula works in the case where $v$ is a valuation. In general, we define  
    \[ v(\alpha) = \max_{x \in I} \min_{y \in J} v(x/y), \]
    when
    \[ \alpha = \bigjoin_{x\in I} \adiv(x) - \bigjoin_{y\in J} \adiv(y) \in \ULatQ(\Field). \]
    Moreover, the same definition makes sense if $v$ is an abstract pseudo-valuation on the field $\Field$.
\end{definition}

\begin{remark}
    Let $t$ be a $\QQ$-tropical polynomial and $a$ be a tuple of elements from $\Field^{\times}$. From now on, when we write $t(v(a))$, we mean $v \bigl( t(\adiv(a)) \bigr)$, for any abstract pseudo-valuation $v$ on $\Field$ or $v \in \OF$. Note that the relation $t(v(a)) = v(b)$ is uniformly definable in any two-sorted structure $(K, \Gamma \cup \{\pm \infty\}, v)$ where $K$ is a field, $\Gamma$ is a divisible ordered abelian group and $v:K \to \Gamma \cup \{\pm\infty\}$ is an abstract pseudo-valuation.
\end{remark}

\begin{lemma} \label{lemma_nice_properties_of_evaluation}
    The pairing $\ev:\Omega_{\Field} \times \ULatQ(\Field) \to [-\infty, \infty]$ has the following properties:
    \begin{itemize}
        \item $v(\alpha \vee \beta) = \max(v(\alpha),v(\beta))$,
        \item $v(\alpha + \beta) = v(\alpha) + v(\beta)$,
    \end{itemize}
    for all $v \in \Omega_{\Field}$, $\alpha, \beta \in \ULatQ(\Field)$, and in the second equation the right-hand side is undefined if one term is $-\infty$ and the other is $\infty$. Moreover, for $\alpha \in \ULatQ(\Field)$, the corresponding map $\ev_{\alpha}:\Omega_{\Field} \to [-\infty, \infty]$ is continuous.
\end{lemma}
\begin{proof}
    Continuity follows from the continuity of $\max$ and $\min$ on $[-\infty, \infty]$. We omit the remaining verifications.
\end{proof}

\begin{definition} \label{definition_ideal_in_a_lattice}
    By an ideal $I$ in a lattice group $\Gamma$ we mean a subgroup such that for all $\alpha\in I$ and $\beta\in \Gamma$,
    $|\beta|\leqslant |\alpha|$ implies that $\beta\in I$. Then the quotient $\Gamma/I$ inherits the operations $\join$ and $\meet$, making it a lattice group.
\end{definition}

Now, we recall the definition of $L^1$-lattices, which can be found for instance in \cite{meyer2012banach}.

\begin{definition} \label{definition:L1_lattice}
    An \emph{$L^1$-lattice} is a lattice group $\Gamma$ which is also a $\mathbb{R}$-vector space, together with a norm $\lVert\cdot \rVert$ such that
    \begin{itemize}
        \item $(\Gamma,\lVert\cdot\rVert)$ is a Banach space.
        \item For all $\alpha, \beta \in \Gamma$ disjoint (i.e., such that $\alpha\meet \beta \leqslant 0\leqslant \alpha\join\beta$), $\lVert \alpha + \beta \rVert = \lVert \alpha \rVert + \lVert \beta \rVert$.
    \end{itemize}
    If $(\Gamma,\lVert\cdot\rVert)$ is an $L^1$-lattice, then we can define a continuous linear form $\int : \Gamma \rightarrow \mathbb{R}$ by
    \[
        \int\alpha := \lVert \alpha^+ \rVert - \lVert \alpha^- \rVert. 
    \]
\end{definition}

Let us mention that every $L^1$-lattice is isomorphic to one of the form $L^1(\mu)$ for some measure $\mu$ on some measure space, see \cite[Theorem 2.7.1]{meyer2012banach}.
\begin{definition} \label{definition:lattice_MVF_structure}
    Let $\Field$ be a field. A \emph{lattice valuation} on $F$ is the data of an $L^1$-lattice $\Gamma$ and a map $\underline{v}: \Field \rightarrow \Gamma \cup\{\infty\}$ satisfying:
    \[
    \underline{v}(xy)=\underline{v}(x)+\underline{v}(y), \quad \underline{v}(x+y) \geqslant \underline{v}(x) \wedge \underline{v}(y) + \underline{v}(2) \wedge 0, \quad \underline{v}(1)=0, \quad \underline{v}(x)=\infty \Leftrightarrow x=0,
    \]
    and such that $\Gamma$ is generated as an $L^1$-lattice by values of $\underline{v}$ (i.e., the divisible sublattice generated by those values is dense in $\Gamma$). Such a lattice valuation is called \emph{global} if for all $x\in \Field^\times$, $\int \underline{v}(x) = 0$. 
    
    There is a natural notion of equivalence of two lattice valuations $\underline{v}:\Field \to \Gamma \cup \{\infty\}$, $\underline{v}':\Field \to \Gamma' \cup \{\infty\}$, namely an isomorphism of $L^1$-lattices $\Gamma \to \Gamma'$ such that the following diagram commutes
    \[\begin{tikzcd}
    	& {\Gamma' \cup \{\infty\}} \\
    	\Field & {\Gamma \cup \{\infty\}}
    	\arrow["{\underline{v}'}", from=2-1, to=1-2]
    	\arrow["{\underline{v}}", from=2-1, to=2-2]
    	\arrow["\simeq"', from=2-2, to=1-2]
    \end{tikzcd}\]
\end{definition}

Note that we need not have imposed the condition of being generated by values of $\underline{v}$, but we do so in order to have a bijective correspondence in Theorem \ref{theorem_equivalent_MVF_structure_definitions}.

\section{Measures} \label{section_measures}

In this section we define the space of lattice divisors of a field $\Field$, and show a correspondence between positive functionals on this space and measures on $\OF$. When $\Field$ is countable, the measures can be assumed to live on the space of absolute values $M_\Field$.
Thus in that case, an adelic curve structure appears on $\Field$, in the sense of \cite{Adelic_curves_1}. 
However, even in the countable case, it is easier to produce measures on $\OF$ first, because of the topological properties of $\OF$.

\begin{definition} \label{definition_lattice_of_a_field}
    Consider the subset $I \subset \ULatQ(\Field)$ consisting of those $\alpha$ such that for all $v \in \Omega_\Field$ we have $\ev(v, \alpha) = 0$. It is an ideal by Lemma \ref{lemma_nice_properties_of_evaluation}. We define a divisible lattice group $\LdivQ(\Field) := \ULatQ(\Field)/I$ and call it the space of \textit{lattice divisors on $\Field$}. In particular, $\alpha \in \LdivQ(\Field)$ is \textit{positive} (written $\alpha \geqslant 0$) if and only if, for all $v \in \Omega_{\Field}$ we have $\ev(v, \alpha) \geqslant 0$.
\end{definition}

Let us consider a relative situation where $\Field$ is a field extending $k$. Using Lemma \ref{Lemma_extending_pseudo_absolute_values} we get the following corollary.

\begin{corollary} \label{corollary_embedding_of_fields_induces_embedding_of_LDivs}
    The natural morphism $\ULatQ(k) \to \ULatQ(\Field)$ induces an embedding $\LdivQ(k) \subset \LdivQ(\Field)$ respecting the non-negative cone.
\end{corollary}

\begin{remark}
    The ideal $I$ is in general not trivial. For example, if $\Field=\QQ$ one can consider an element
    \[ \alpha = \adiv(pq) \vee 0 - \adiv(p) \vee \adiv(q) \vee 0, \]
    for some primes $p \neq q$. Then, by using Ostrowski's theorem and direct inspection, for any pseudo-valuation $v$ on $\QQ$ we have $v(\alpha)=0$. However, if $\alpha = 0$ in $\ULatQ(\QQ)$, this would mean that there is a finite subset $A \subset \QQ^{\times}$ such that
    \[ A \cup pqA = A \cup pA \cup qA, \]
    but looking at the element of largest Euclidean norm in $A$, we get that this is impossible.
\end{remark}

Recall that for a locally compact Hausdorff topological space $X$, positive functionals on the space of compactly supported continuous functions on $X$ correspond to Radon measures on $X$ (this correspondence sometimes serves as a definition of a Radon measure). 

\begin{definition}
    Let $\alpha \in \LdivQ(\Field)$ be a non-zero, non-negative element. We define a compact space
    \[ \Omega_{\Field}(\alpha) = \{ v \in \Omega_{\Field} : v(\alpha) = 1\}, \]
    and its projection to $\ov{\Omega}_{\Field}$ denoted $\ov{\Omega}_{\Field}(\alpha)$. Moreover, let
    \[ \langle \alpha \rangle := \{ \beta \in \LdivQ(\Field) : (\exists n \in \NN)(|\beta| \leqslant n \cdot \alpha) \}. \]
    Note that $\langle \alpha \rangle$ is a sublattice of $\LdivQ(\Field)$. We call it the set of elements dominated by $\alpha$.
\end{definition}

\begin{construction} \label{construction:pairing_divisors_valuations_local}
    Fix $\alpha \geqslant 0$ from $\LdivQ(\Field)$ and consider the restricted evaluation product
    \[ \Omega_{\Field}(\alpha) \times \langle \alpha \rangle \to \RR. \]
    It defines a natural map
    \[ i:\langle \alpha \rangle \to C(\Omega_{\Field}(\alpha)), \]
    given by $\beta \mapsto (v \mapsto v(\beta))$. Note that by the max-closed/tropical Stone-Weierstrass theorem, the image of $i$ is a dense subset with respect to the uniform convergence topology.
\end{construction}

\begin{definition} \label{definition_positive_functional}
    For a single element $f \in \Field^{\times}$ we call $\adiv(f)$ the principal lattice divisor associated with $f$. A \textit{positive functional on $\Field$} is a $\QQ$-linear map $l:\LdivQ(\Field) \to \RR$ which takes non-negative lattice divisors to non-negative numbers. We call $l$ a \textit{global functional} (or a GVF \textit{functional}) if additionally it takes principal lattice divisors to zero. Note that it can equivalently be described as a linear map $l:\LPicQ(\Field) \to \RR$ non-negative on the effective cone. Furthermore, we say that $l$ is \textit{normalised} if $\Field$ is of characteristic zero and $l(\adiv(2)^{+}) = \log 2$.
\end{definition}

Now, the following corollary follows from applying the Riesz-Markov-Kakutani representation theorem to the pairing from Construction \ref{construction:pairing_divisors_valuations_local}.

\begin{corollary} \label{corollary_positive_functional_gives_local_measures}
    Let $l$ be a positive functional on $\Field$. Then for every non-zero $\alpha \geqslant 0$ in $\LdivQ(\Field)$, there is a unique finite Radon measure $\mu_{\alpha}$ on $\Omega_{\Field}(\alpha)$ such that for all $\beta \in \langle \alpha \rangle$ we have
    \[ l(\beta) = \int_{\Omega_{\Field}(\alpha)} v(\beta) d\mu_{\alpha}(v). \]
    We also use the notation $\mu_{\alpha}$ for the pushforward of this measure to $\ov{\Omega}_{\Field}(\alpha)$.
\end{corollary}

Hence, a positive functional $l$ on $\Field$ induces a family of measures $\mu_{\alpha}$ on $\ov{\Omega}_{\Field}(\alpha)$ for non-zero, non-negative $\alpha \in \LdivQ(\Field)$.
\begin{lemma}
    For any non-zero, non-negative $\alpha, \gamma \in \LdivQ(\Field)$ the local measures $\mu_{\alpha}, \mu_{\gamma}$ associated to a positive functional $l$ on $\Field$ satisfy
    \[ \frac{d\mu_{\alpha}}{d \mu_{\gamma}} = \frac{v(\alpha)}{v(\gamma)}, \]
    when restricted to $\ov{\Omega}_{\Field}(\alpha) \cap \ov{\Omega}_{\Field}(\gamma)$, where $\frac{v(\alpha)}{v(\gamma)}$ is a function taking the class of $v$ to $\frac{v(\alpha)}{v(\gamma)}$ (which does not depend on the representative).
\end{lemma}
\begin{proof}
    First, note that we have
    \[ l(\gamma) \geqslant l(\gamma \wedge n \alpha) = \int_{\Omega_{\Field}(\alpha)} v(\gamma \wedge n \alpha) d\mu_{\alpha}(v) \geqslant n \cdot \mu_{\alpha}(\{ v \in \Omega_{\Field}(\alpha) : v(\gamma) = \infty \}). \]
    Since $n$ is an arbitrary natural number here, we get that
    \[ \mu_{\alpha}(\{ v \in \Omega_{\Field}(\alpha) : v(\gamma) = \infty \}) = 0. \]
    Thus, for $\beta \in \langle \alpha \rangle \cap \langle \gamma \rangle$ we can write
    \[ l(\beta) = \int_{\Omega_{\Field}(\alpha)} v(\beta) d\mu_{\alpha}(v) = \int_{\ov{\Omega}_{\Field}(\alpha) \cap \ov{\Omega}_{\Field}(\gamma)} \frac{v(\beta)}{v(\alpha)} d\mu_{\alpha}(v). \]
    Symmetrically, we have
    \[ l(\beta) = \int_{\Omega_{\Field}(\gamma)} v(\beta) d\mu_{\gamma}(v) = \int_{\ov{\Omega}_{\Field}(\alpha) \cap \ov{\Omega}_{\Field}(\gamma)} \frac{v(\beta)}{v(\gamma)} d\mu_{\gamma}(v). \]
    This finishes the proof as functions sending a class of $v$ to $\frac{v(\beta)}{v(\alpha)}$, for $\beta \in \langle \alpha \rangle \cap \langle \gamma \rangle$ have dense real span in the space of continuous functions on $\ov{\Omega}_{\Field}(\alpha) \cap \ov{\Omega}_{\Field}(\gamma)$.
\end{proof}

For $a \in \Field^{\times}$ we denote by $\mu_a$ the measure on $\ov{\Omega}_{\Field}(a)$ coming from $\alpha = \adiv(a)^{+}$. Note that for this choice we have $\ov{\Omega}_{\Field}(a) = \ov{\Omega}_{\Field}(\alpha)$.

\begin{definition} \label{definition_family_of_local_measures}
    By a \textit{family of local measures on $\Field$} we mean a family of finite Radon measures $\mu_a$ on $\ov{\Omega}_{\Field}(a)$ for $a \in \Field^{\times}$, such that for any $a, b \in \Field^{\times}$ we have
    \[ \frac{d\mu_{a}}{d \mu_{b}} = \frac{v(a)}{v(b)} \]
    on $\ov{\Omega}_{\Field}(a) \cap \ov{\Omega}_{\Field}(b)$ and $\mu_a \{v \in \ov{\Omega}_{\Field}(a) : v(b) = \infty \} = 0$. We call such a family \textit{global} if in addition, for every $a \in \Field^{\times}$, we have
    \[ \mu_a(\ov{\Omega}_{\Field}(a)) = \mu_{a^{-1}}(\ov{\Omega}_{\Field}(a^{-1})). \]
\end{definition}

Now we pass to a global description of GVF functionals.

\begin{definition} \label{definition_global_measure}
    An \textit{admissible measure on $\Field$} is a Radon measure $\mu$ on $\Omega_{\Field}$ which satisfies the following conditions.
    \begin{enumerate}
        \item For every $a \in \Field^{\times}$ the function $\ev_a:v \mapsto v(a)$ is $\mu$-integrable.
        \item There is a family of Borel sections $s_a:\oOF(a) \to \pi^{-1}(\OF(a))$ of $\pi$, such that $\mu$ restricted to $\pi^{-1}(\oOF(a))$ is equal to $s_{a*}\mu_a$ for $\mu_a := \pi_*(\mu|_{\pi^{-1}(\OF(a))})$.
    \end{enumerate}
    Note that the Borel sections $s_a$ are uniquely determined, up to a set of $\mu_a$-measure zero.
    Moreover, we call $\mu$ a \textit{global measure on $\Field$} (or a GVF measure) if in addition for all $a \in \Field^{\times}$ the product formula holds:
    \[ \int_{\Omega_{\Field}} v(a) d\mu(v) = 0. \]
\end{definition}

\begin{remark} \label{remark_pre_valuations_measure_can_be_on_them}
    Note that the condition (1) in the above definition implies in particular that $\mu \{v(2)=-\infty\} =0$, hence $\mu$ can be seen as a measure on the set of pseudo-valuations $\Omega^{\circ}_{\Field}$. We choose to work with the bigger space $\Omega_{\Field}$ to ensure that the sets $\Omega_{\Field}(a)$ are compact.
\end{remark}

\begin{definition} \label{definition_renormalisation_of_an_admissible_measure_new}
    Let $\mu, \mu'$ be two admissible measures on $\Field$. We say that \textit{$\mu$ can be renormalised to $\mu'$} if for every $a \in \Field^{\times}$ and Borel sections $s_a, s'_a: \oOF(a) \to \pi^{-1}(\oOF(a))$ associated to measures $\mu, \mu'$ respectively, the measures $\mu_a, \mu'_a$ obtained as pushforwards of $\mu, \mu'$ through $\pi:\OF(a) \to \oOF(a)$ respectively, satisfy:
    \[ \frac{d\mu_a}{d\mu'_a} = \frac{s'_a}{s_a}. \]
\end{definition} 

\begin{lemma}
    Renormalisation is an equivalence relation among admissible measures on $\Field$ that respects being a global measure.
\end{lemma}
\begin{proof}
    The first part is clear. The second part follows from Lemma \ref{lemma_admissible_measures_define_MVF_structures}.
\end{proof}

We now describe how to construct an admissible measure on $\Field$ from a family of local measures.

\begin{construction} \label{construction_admissible_measure_from_local_measures}
     Fix a family of local measures $(\mu_a)_{a \in \Field^{\times}}$. Well order $\Field^{\times}$ as $(a_{\alpha})_{\alpha < \kappa}$ for some ordinal $\kappa$. We inductively build sections $s_{\alpha}:\oOF(a_{\alpha}) \to \OF$ such that for $\alpha, \beta <\kappa$ we have
    \[ s_{\alpha}|_{\oOF(a_{\alpha}) \cap \oOF(a_{\beta})} = s_{\beta}|_{\oOF(a_{\alpha}) \cap \oOF(a_{\beta})} \]
    almost everywhere with respect to $\mu_{a_{\alpha}}$ or $\mu_{a_{\beta}}$ equivalently. Assume that we have already built such $s_{\gamma}$ for $\gamma < \alpha$ and let $a = a_{\alpha}$. Pick a countable sequence of $\gamma_n < \alpha$ such that for $a_n = a_{\gamma_n}$ the set $\bigcup_n \oOF(a) \cap \oOF(a_n)$ has maximal $\mu_a$-measure. Define $s = s_{\alpha}:\oOF(a) \to \OF$, at least up to $\mu_a$-measure zero, as $s_n = s_{\gamma_n}$ on $\oOF(a) \cap \oOF(a_n)$ (this makes sense by the inductive assumption about $s_n$'s) and on the rest define it to be the inverse of the projection $\pi:\OF(a) \to \oOF(a)$. Let $\nu_{\alpha} := g_{\alpha} \cdot s_{\alpha*} \mu_{a}$, where $g_{\alpha}(v) = 1/s_{\alpha}(\pi(v))(a_{\alpha})$. Next, define a measure $\mu$ on $\OF$ to be the unique Radon measure which coincides with $\nu_{\alpha}$ on $\pi^{-1}(\oOF(a_\alpha))$. This makes sense, as every compact subset of $\OF$ is contained in a finite union of $\pi^{-1}(\oOF(a))$'s and the choice of $g_{\alpha}$'s ensures that $\nu_{\alpha}$'s agree on overlaps because of the following calculation (on $\pi^{-1}(\oOF(a_\alpha)) \cap \pi^{-1}(\oOF(a_\beta))$): 
    \[ \nu_{\alpha} = g_{\alpha} \cdot s_{\alpha*} \mu_{a} = g_{\alpha} \cdot s_{\alpha*} \biggl( \frac{v(a_{\alpha})}{v(a_{\beta})} \mu_{\beta} \biggr) = \frac{1}{s_{\alpha}(\pi(v))(a_{\alpha})} \cdot \frac{s_{\alpha}(\pi(v))(a_{\alpha})}{s_{\alpha}(\pi(v))(a_{\beta})} \cdot s_{\alpha*} \mu_{\beta} \]
    \[ \frac{1}{s_{\alpha}(\pi(v))(a_{\beta})} \cdot s_{\alpha*} \mu_{\beta} = \frac{1}{s_{\beta}(\pi(v))(a_{\beta})} \cdot s_{\beta*} \mu_{\beta} = g_{\beta} \cdot s_{\beta*} \mu_{\beta} = \nu_{\beta}. \]
    Since
    \[ \nu_{\alpha} = g_{\alpha} \cdot s_{\alpha*} \mu_{a} = s_{\alpha*} \biggl( \frac{1}{s_{\alpha}(\ov{v})(a_{\alpha})} \mu_{a} \biggr),\]
    where $\ov{v} \in \oOF(a_{\alpha})$, the measure $\mu$ satisfies the second condition in the definition of an admissible measure. Moreover, by construction for every $a = a_{\alpha} \in \Field^{\times}$ (with $s_a = s_{\alpha}, \nu_a = \nu_{\alpha}$) we have
    \[ \int_{\OF} v(a)^{+} d\mu(v) = \int_{\pi^{-1}(\oOF(a))} v(a) d\mu(v) = \int_{\pi^{-1}(\oOF(a))} v(a) d\nu_a(v) \]
    \[ = \int_{\pi^{-1}(\oOF(a))} \frac{v(a)}{s_{a}(\pi(v))(a)} ds_{a*} \mu_{a}(v) = \int_{\oOF(a)} d\mu_a = \mu_a(\oOF(a)) < \infty. \]
    Thus, $\mu$ also satisfies the first condition from the definition of an admissible measure, and a similar calculation for $v(a)^{-}$ shows that a global family of local measures induces a global measure on $\Field$.
\end{construction}

\begin{remark} \label{remark_GVF_measures_can_be_pushed_from_M_F}
    Note that in the above construction, if we start with $a_0 = 1/2$ in the characteristic zero case, then we can assume that the obtained measure $\mu$ on $\OF$ restricted to Archimedean $v \in \OF$ is concentrated on $v$'s with $v(2) = -\log 2$. In particular, if $\Field$ is a countable field of characteristic zero, then by Remark \ref{remark_char_of_prevaluations_and_prenorms} we can assume that the measure comes from a pushforward of some measure via $-\log:M_{\Field} \to \Val_{\Field} \subset \Omega_{\Field}$. More precisely, we can even assume that the measure on $M_{\Field}$ is concentrated on the union of non-Archimedean absolute values and Archimedean absolute values extending $|\cdot|_{\infty}$ on $\QQ$.
\end{remark}

\section{Positivity}
In this section we give an algebraic criterion for a lattice divisor to be positive.

\begin{definition}
    We call lattice divisors of the form $\adiv(f)$ \textit{principal}, for $f \in \Field_{\QQ}^{\times}$. The quotient of $\LdivQ(\Field)$ by the subspace of principal lattice divisors is called the space of \textit{lattice line bundles on $\Field$} and is denoted by $\LPicQ(\Field)$. The \textit{effective cone} of $\LPicQ(\Field)$ is defined as the set of classes of positive lattice divisors on $\Field$.
\end{definition}

\begin{definition} \label{definition:positive_pairs}
    Consider pairs $(\phi, t) = (\phi(x,y), t(x))$, where $x, y$ are tuples of variables, $\phi(x,y)$ is a quantifier-free ring-language formula implying that none of the variables of $x$ is zero, and $t(x)$ is a $\QQ$-tropical polynomial. We call such a pair \textit{positive} if
    \[ \VF \models (\forall a,b)(\phi(a,b) \implies t(v(a)) \geqslant 0 ). \]
    By $\VF$ we mean the common theory of triples $(K, \Gamma \cup \{\pm\infty\}, v)$, where $v$ is an abstract pseudo-valuation on a field $K$ and $\Gamma$ is the corresponding divisible ordered abelian group.
\end{definition}

\begin{remark}
    By adding dummy variables, for any tuples of variables $x,y$ we can see a $\QQ$-tropical polynomial $t(x)$ as $t(x,y)$. Thus, we can omit the second variable $y$ in the definition of positive pairs.
\end{remark}

\begin{proposition} \label{proposition_relative_positivity_of_lattice_divisors}
    Let $k$ be a subfield of $\Field$ and $\alpha \in \LdivQ(\Field)$. The following are equivalent:
    \begin{enumerate}
        \item for every $v \in \Omega_{\Field}$ that restricts to the trivial valuation on $k$, we have $v(\alpha) \leqslant 0$;
        \item there are finitely many $a_1, \dots, a_n \in k^{\times}$ such that $\alpha \leqslant \bigvee_{i=1}^n \adiv(a_i)$.
    \end{enumerate}
\end{proposition}
\begin{proof}
    The implication from (2) to (1) follows from Lemma \ref{lemma_nice_properties_of_evaluation}. To prove the other one, assume that one cannot find finitely many elements in $k$ that satisfy (2). 
    This means that for every $a_1, \dots, a_n \in k^{\times}$, there is a pseudo-valuation $v \in \Omega_{\Field}$ that satisfies $v(\alpha) > \max_i v(a_i)$. Consider the following first order theory:
    \[ T = \Diagat(\Field) \cup \{v \textnormal{ is an abstract pseudo-valuation}\} \cup \{ v(\alpha) > v(b) : b \in k^{\times} \}.  \]
    This theory is consistent by the assumption, so there is an abstract pseudo-valuation $v_0:\Field \to \Gamma \cup \{\pm \infty\}$ such that for all $b \in k^{\times}$ we have $v_0(\alpha) > v_0(b)$.
    First, assume that $v_0(\alpha) \neq \infty$.
    Let $\Delta$ be the smallest convex subgroup containing $v(\alpha)$ and let $\Delta_0 \subset \Gamma$ be the biggest convex subgroup of $\Gamma$ not containing $v(\alpha)$. Then by restricting $v_0$ to $\Delta$ and further by taking the quotient by $\Delta_0$, we get a rank one valuation $v:F \to \Delta/\Delta_0 \cup \{\pm\infty\}$ which refutes (1) once we choose an embedding $\Delta/\Delta_0 \subset \RR$.
    We are left with the case when $v_0(\alpha) = \infty$.
    Write $$\alpha = q \cdot \left(\bigvee_i \adiv(x_i) - \bigvee_j \adiv(y_j)\right)$$ for some $q \in \QQ_{>0}$ and $x_i, y_j$ from $\Field^{\times}$. Let $\cO = \{x \in \Field: v_0(x) \neq -\infty \}$ be the valuation ring of $v_0$. Choose $i_0$ such that for all and $j$ we have $v_0(x_{i_0}/y_j) =\infty$.  
    Let $u:\Field^{\times} \to \Field^{\times}/\cO^{\times} =: \Sigma$ and note that it is a non-Archimedean abstract valuation.
    Moreover, $u$ is trivial on $k^{\times}$ and $u(x_{i_0}/y_j) > 0$ for every $j$.
    Hence $\gamma = u(\alpha) > 0$.
    Define $\Delta \subset \Sigma$ as the smallest convex subgroup containing $\gamma$ and let $\Delta_0 \subset \Delta$ be the biggest convex subgroup of $\Delta$ without $\gamma$. Then by truncating $u$ to $\Delta$ and taking the quotient by $\Delta_0$ we get an abstract pseudo-valuation on $\Field$ where the value group $\Delta/\Delta_0$ can be embedded to $\RR$ and which satisfies $u(\alpha) > 0$. 
    This finishes the proof.
\end{proof}

\begin{corollary} \label{corollary_ideal_vanishing_relative_situation}
    Let $k \subset \Field$ be a finitely generated extension. Then the ideal
    \[ I_k = \{ \alpha \in \LdivQ(\Field) : \forall v \in \Val_{\Field/k} \ v(\alpha) = 0 \} \]
    is the smallest ideal containing $\LdivQ(k)$ in $\LdivQ(\Field)$.
\end{corollary}
Here, the embedding $\LdivQ(k) \subset \LdivQ(\Field)$ is the one from Corollary \ref{corollary_embedding_of_fields_induces_embedding_of_LDivs} and by $\Val_{\Field/k}$ we denote valuations on $\Field$ that are trivial on $k$. Let us also denote by $\Omega_{\Field/k}$ pseudo-valuations on $\Field$ that restrict to the trivial valuation on $k$.
\begin{proof}
    By Corollary \ref{corollary_density_of_norms_on_a_countable_field_relative_version} we can replace $\Val_{\Field/k}$ in the definition of $I_k$ by $\Omega_{\Field/k}$. Then the statement follows from Proposition \ref{proposition_relative_positivity_of_lattice_divisors}.
\end{proof}

Now, we want to describe the set of positive lattice divisors in more detail, see Corollary \ref{corollary:tropical_divisors_positivity_conditions}. Let us start with the following lemma.

\begin{lemma} \label{lemma:pure_divisor_negativity_condition}
    For every tuple $a = (a_1,\ldots,a_n) \in (F^{\times})^n$ and $\varepsilon \in \RR_{>0}$, the following are equivalent:
    \begin{enumerate}[label=(\roman*)]
        \item For every pseudo-valuation $v$ on $F$, $v\left(\bigmeet\limits_{i=1}^n\adiv(a_i)\right) \leqslant 0$ if $v$ is non-Archimedean, and $v\left(\bigmeet\limits_{i=1}^n\adiv(a_i)\right) < -\varepsilon v(2)$ if $v$ is Archimedean.
        \item There exists a family $(m_s)_{s\in\mathbb{N}^n}$ of integers with $m_s=0$ except for a finite number of indices $s$, such that 
        \[
            \sum_{s\in\mathbb{N}^n} |m_s|2^{-\varepsilon|s|} < 1,
        \]
        and
        \[
            \sum_{s\in\mathbb{N}^n} m_s a^s = 1.
        \]
    \end{enumerate}
\end{lemma}

\begin{proof}
    First assume that $\textnormal{(ii)}$ is satisfied, and that there is a $v\in\Omega_F$ such that $v\left(\bigmeet\limits_{i=1}^n\adiv(a_i)\right) > -\varepsilon v(2)^{-}$. Let $(m_s)_{s\in\mathbb{N}^n}$ be a family satisfying $\textnormal{(ii)}$. Then:
    \begin{itemize}
        \item If $v$ is non-Archimedean, we have $\min\limits_{i\leqslant n}v(a_i) > 0$, so $v(a_i) > 0$ for all $i\leqslant n$. Thus, the equality $\sum_{s\in\mathbb{N}^n} m_s a^s = 1$ implies $v(1)>0$.
        \item If $v$ is Archimedean, we may, up to multiplication by a positive constant, write $v = -\log\lVert \cdot \rVert$, where $\lVert \cdot\rVert$ is a pseudo-absolute value on $F$ with $\|2\|=2$. We then have $\min\limits_{i\leqslant n}v(a_i) > -\varepsilon v(2) = \varepsilon \log(2)$, so $\lVert a_i\rVert < 2^{-\varepsilon}$. Thus,
        \[
            1 = \lVert 1\rVert = \left\lVert \sum_{s\in\mathbb{N}^n} m_sa^s \right\rVert \leqslant \sum_{i=1}^n |m_s|2^{-\varepsilon|s|} < 1.
        \]
    \end{itemize}
    Either way, we obtain a contradiction, which shows that $\textnormal{(ii)}$ implies $\textnormal{(i)}$.

    Now, assume that $\textnormal{(ii)}$ is not satisfied. Consider the pseudo-seminorm $\lVert \cdot \rVert_{\varepsilon}$ defined on $\mathbb{Z}[X_1,\ldots,X_n]$ by
    \[
        \left\lVert \sum_{s\in\mathbb{N}^n} a_sX^s\right\rVert_{\varepsilon} = \sum_{s\in\mathbb{N}^n} |a_s| 2^{-\varepsilon|s|},
    \]
    and the ring homomorphism $f:\mathbb{Z}[X_1,\ldots,X_n]\rightarrow F$ sending $X_i$ to $a_i$ for all $i\leqslant n$. Then, $\lVert\cdot\rVert_\varepsilon$ and $f$ satisfy the conditions of Lemma \ref{lemma:existence_prenorm_ring_map}, hence there exists a pseudo-absolute value $\lVert \cdot\rVert$ of $F$ such that for all $P\in\mathbb{Z}[X_1,\ldots,X_n]$, $\lVert f(P)\rVert \leqslant \lVert P\rVert_\varepsilon$. In particular $\lVert a_i\rVert \leqslant 2^{-\varepsilon}$ for all $i\leqslant n$. 
    
    Assume that $\lVert \cdot\rVert$ is a non-trivial pseudo-absolute value. Then by Remark \ref{remark_char_of_prevaluations_and_prenorms}, $v = -\log\lVert\cdot\rVert$ is a pseudo-valuation on $F$. But $v$ satisfies $\min_{1\leqslant i\leqslant n} v(a_i) \geqslant \varepsilon \log(2)$, which contradicts (i) since, if $v$ is Archimedean, then $v(2) = -\log\lVert 2\rVert \geqslant -\log(2)$.
    
    On the other hand, if $\lVert \cdot\rVert$ is a trivial pseudo-absolute value, then $\lVert a_i\rVert\leqslant 2^{-\varepsilon} < 1$ implies $\lVert a_i\rVert = 0$ for all $1\leqslant i \leqslant n$. Let $\mathcal{O} = \{ x \in F \,:\, \lVert x \rVert < \infty\}$ be the valuation ring associated to $\lVert\cdot\rVert$, and let $v':F \rightarrow \Sigma\cup\{+\infty\}$ be the corresponding abstract valuation, where $\Sigma := F^\times/\mathcal{O}^{\times}$. Then $\beta := \min_{1\leqslant i\leqslant n} v'(a_i) > 0$ in $\Sigma$. Let $\Delta \subset \Sigma$ be the smallest convex subgroup of $\Sigma$ containing $\beta$, and $\Delta_0$ be the biggest convex subgroup not containing $\beta$. Then, by truncating $v'$ to $\Delta$ and taking the quotient by $\Delta_0$, we obtain an abstract non-Archimedean pseudo-valuation $v$ whose value group $\Delta/\Delta_0$ can be embedded into $\mathbb{R}$. Hence, $v$ is a non-Archimedean pseudo-valuation on $F$ and by construction satisfies $\min_{1\leqslant i\leqslant n} v(a_i) > 0$, which contradicts (i).
\end{proof}

\begin{definition} \label{definition:gamma_epsilon}
    For every $\varepsilon \in\mathbb{R}_{>0}$, let $\Gamma_\varepsilon \subset \ULatQ(F)$ be the set of elements of the form
    \[
        \bigmeet_{i=1}^p \adiv\left(\sum_{j=1}^q m_{i,j}a_j\right) - \bigmeet_{j=1}^q \adiv(a_j),
    \]
    where $a_1,\ldots,a_q\in F^\times$, $p,q\in\mathbb{N}^*$, and $(m_{i,j})_{\substack{1\leqslant i \leqslant p \\ 1 \leqslant j \leqslant q}} \in \mathbb{Z}^{pq}$ is such that for all $i\leqslant p$, $\sum\limits_{j=1}^q |m_{i,j}| < 2^\varepsilon$.
\end{definition}

\begin{remark} \label{remark:gamma_epsilon_sublattice}
    For all $\varepsilon>0$, $\Gamma_\varepsilon$ is closed under $\meet$.
    Indeed, if two elements of $\ULatQ(F)$ are written as $\bigmeet\limits_{1\leqslant i\leqslant m} \adiv(b_i) - \bigmeet\limits_{1\leqslant i\leqslant n} \adiv(a_i)$ and $\bigmeet\limits_{1\leqslant i\leqslant p} \adiv(b_i') - \bigmeet\limits_{1\leqslant i\leqslant q}\adiv(a_i')$ respectively, then their meet is equal to:
\[
    \left(\bigmeet_{\substack{1\leqslant i\leqslant m\\ 1\leqslant j\leqslant q}} \adiv(b_i a_j')\right) \meet \left(\bigmeet_{\substack{1\leqslant i\leqslant n\\ 1\leqslant j\leqslant p}} \adiv(b_j'a_i)\right) - \bigmeet_{\substack{1\leqslant i\leqslant n\\ 1\leqslant j\leqslant q}} \adiv(a_ia_j').
\]
Thus, if both elements of $\ULatQ(F)$ are of the form described in Definition \ref{definition:gamma_epsilon}, then their meet is also of this form.
\end{remark}

\begin{lemma} \label{lemma:divisor_gamma_epsilon_equivalence}
    For any $\alpha\in\ULatQ(F)$ and $\varepsilon\in\mathbb{R}_{>0}$, the following are equivalent:
    \begin{enumerate}[label=(\roman*)]
        \item For every pseudo-valuation $v$ on $F$, $v(\alpha) \geqslant 0$ if $v$ is non-Archimedean, and $v(\alpha) > \varepsilon v(2)$ if $v$ is Archimedean.
        \item There exists an integer $m\in\mathbb{N}$ such that $m\alpha\in \Gamma_{m\varepsilon}$.
    \end{enumerate}
\end{lemma}

\begin{proof}
$\textnormal{(ii)} \Rightarrow \textnormal{(i)}$ Assume that $m\alpha \in \Gamma_{m\varepsilon}$ for some $m\in\mathbb{N}$. Then, $m\alpha$ can be written as $\bigmeet_{i=1}^p \adiv(b_i) - \bigmeet_{j=1}^q \adiv(a_j)$
    where $a_1,\ldots,a_q,b_1,\ldots,b_p\in F^\times$, and  for all $i\leqslant p$, there is a tuple $(m_{i,1},\ldots,m_{i,q})\in \mathbb{Z}^q$ such that $b_i = \sum\limits_{j=1}^q m_{i,j}a_j$, and $\sum_{j=1}^q |m_{i,j}| < 2^{m\varepsilon}$. Then, for every pseudo-valuation $v$~ and $i\leqslant p$:
    \begin{itemize}
        \item if $v$ is non-Archimedean,  $\min\limits_{1\leqslant j \leqslant q} v(a_j/b_i) \leqslant 0$,
        \item if $v$ is Archimedean, $\min\limits_{1\leqslant j \leqslant q} v(a_j/b_i) < -m\varepsilon v(2)$.
    \end{itemize}
    Using the fact that $m\cdot v(\alpha) = v(m\alpha) = \min\limits_{1\leqslant i\leqslant p} \max\limits_{1\leqslant j\leqslant q} v(b_i/a_j)$, we conclude that $\alpha$ satisfies (i).
    
$\textnormal{(i)} \Rightarrow \textnormal{(ii)}$ Since any element of $\ULatQ(F)$ can be written as a meet of elements of the form $-\bigmeet_{i=1}^n \adiv(a_i)$, with $a=(a_1,\ldots,a_n)\in (F^\times)^n$, and by Remark \ref{remark:gamma_epsilon_sublattice}, it is enough to prove $\textnormal{(i)}\Rightarrow\textnormal{(ii)}$ for elements of this form. Let $\alpha = -\bigmeet_{i=1}^n \adiv(a_i)$ be such that $\alpha$ satisfies $\textnormal{(i)}$.

Let $\lVert\cdot\rVert_\varepsilon$ be the ring norm on $\mathbb{Z}[X_1,\ldots,X_n]$ defined by
\[
    \left\lVert \sum_{s\in\mathbb{N}^n} m_sX^s \right\rVert = \sum_{s\in\mathbb{N}^n} |m_s| 2^{-\varepsilon|s|}.
\]
Then, by Lemma \ref{lemma:pure_divisor_negativity_condition}, there exists a polynomial $P\in\mathbb{Z}[X_1,\ldots,X_n]$ with $\lVert P\rVert_\varepsilon < 1$, such that $P(a)=1$. 

\begin{claim} \label{claim:divisor_positivity_condition}
We have $\alpha\in \Gamma_{\delta}$, with $\delta = \log_2(n+ 2^{1+\varepsilon})$.
\end{claim}
\begin{proof}
We define the (finite) subset $E$ of $\mathbb{Z}[X_1,\ldots,X_n]$ consisting of polynomials $H$ of degree $\deg H\leqslant \deg P$, such that $\lVert H\rVert_\varepsilon\leqslant 1$ and $H(a)\neq 0$. Then, $\alpha$ may be written as
\[
    \bigmeet_{H\in E} \adiv(H(a)) - \bigmeet_{\substack{1\leqslant i \leqslant n\\ H\in E}} \adiv(a_i H(a)).
\]
We now show that $\alpha\in \Gamma_{\delta}$ by writing each $H\in E$ as a sum of at most $2^\delta$ terms of the form $a_i H'$ with $1\leqslant i \leqslant n$, $H'\in E$.

\begin{itemize}
    \item 
If $H$ has no constant term, then we may write $H = \sum\limits_{i=1}^n X_iH_i$ with $H_1,\ldots,H_n$ such that $\lVert H\rVert_\varepsilon = 2^{-\varepsilon}\sum\limits_{i=1}^n \lVert H_i\rVert_\varepsilon$. Then, each $H_i$ may be written as a sum of at most $\max(1,2\lVert H_i\rVert_{\varepsilon})$ elements of $E$. Indeed, if $\lVert H_i\rVert_{\varepsilon} \leqslant 1$ we are done, and otherwise we split $H$ into unitary monomials, and keep regrouping them until the average of their $\varepsilon$-norms is greater than $\frac{1}{2}$.
Hence, $H(a)$ is a sum of at most
\[
    \sum_{i=1}^n\max(1,2\lVert H_i\rVert_{\varepsilon}) \leqslant n + 2^{1+\varepsilon}\lVert H\rVert_\varepsilon \leqslant n + 2^{1+\varepsilon}
\]
terms of the form $a_iH'(a)$, for $H'\in E$.
    \item
If $H$ has a constant term, then $\lVert H\rVert_{\varepsilon}\leqslant 1$ implies that $H=\pm 1$. Therefore,
\[
    H(a) = \pm P(a),
\]
and we are reduced to the previous case.
\end{itemize}
\end{proof}
Now, by Lemma \ref{lemma:pure_divisor_negativity_condition}, $\alpha$ also satisfies $\textnormal{(i)}$ for some $0 < \varepsilon' < \varepsilon$. But, for all $m\in\mathbb{N}$, we have
\[
    m\alpha = -\bigmeet_{P} \adiv(P(a)),
\]
where the meet ranges over all homogeneous monomials of degree $m$ in $X_1,\ldots,X_n$, whose number is $\binom{n+m-1}{n-1}$. Thus, Claim \ref{claim:divisor_positivity_condition} applies, replacing $\alpha$ by $m\alpha$, $\varepsilon$ by $m\varepsilon'$ and $n$ by $\binom{n+m-1}{n-1}$. This implies that $m\alpha\in \Gamma_{\delta_m}$ with
\[
    \delta_m = \log_2\left(\binom{n+m-1}{n-1} + 2^{1+m\varepsilon'}\right)
\]
Since $\varepsilon'<\varepsilon$, for large enough $m$, we have $\delta_m < m\varepsilon$, which proves the lemma.
\end{proof}

\begin{corollary} \label{corollary:tropical_divisors_positivity_conditions}
    For each $\alpha\in\ULatQ(F)$, the following conditions are equivalent:
\begin{enumerate}[label=(\roman*)]
    \item For every pseudo-valuation $v$ on $F$, $v(\alpha) \geqslant 0$.
    \item For every $v\in\Omega_F$, $v(\alpha) \geqslant 0$.
    \item For every $\varepsilon\in\mathbb{R}_{>0}$, there is some $m\in\mathbb{N}$ such that $m\alpha\in\Gamma_{m\varepsilon}$.
\end{enumerate}
\end{corollary}

\begin{proof}
    Note that $\textnormal{(i)}\Leftrightarrow\textnormal{(ii)}$ since the set of pseudo-valuations is a dense subset of $\Omega_F$. 
    On the other hand, $\textnormal{(i)}\Leftrightarrow\textnormal{(iii)}$ is a consequence of Lemma \ref{lemma:divisor_gamma_epsilon_equivalence}. 
\end{proof}

If the ambient field is countable, positivity of an element $\alpha \in \LdivQ(\Field)$ can be checked by pairing with not necessarily all pseudo-valuations.

\begin{lemma} \label{lemma_density_of_valuations_on_a_countable_field}
    Assume that $\Field$ is countable and let $\alpha = t(\adiv(a)) \in \LdivQ(\Field)$. Then $\alpha \geqslant 0$ if and only if for all $v \in \Val_{\Field}$ we have $t(v(a)) \geqslant 0$.
\end{lemma}
\begin{proof}
    By Lemma \ref{lemma_nice_properties_of_evaluation} the evaluation function $v \mapsto t(v(a))$ is continuous. Thus the result follows from Corollary \ref{corollary_density_of_norms_on_a_countable_field} and Corollary \ref{corollary_density_of_norms_on_a_countable_field_relative_version}.
\end{proof}

\begin{corollary} \label{corollary_description_of_LDiv_Arakelov_divisors}
    Let $\Field$ be a countable field of characteristic zero. Then the divisible lattice group $\LdivQ(\Field)$ is isomorphic to the one introduced in \cite{szachniewicz2023existential}, i.e., to the divisible lattice group of arithmetic divisors in integral $\Field$-submodels, spanned by principal arithmetic divisors $\adiv(f), f \in \Field^{\times}$. 
\end{corollary}
\begin{proof}
    From Lemma \ref{lemma_density_of_valuations_on_a_countable_field} it follows that for $\alpha \in \LdivQ(\Field)$ we have $\alpha = 0$ if and only if for all $v \in \Val_{\Field}$ we have $v(\alpha) = 0$. Similarly, e.g. by \cite[Section 2.2]{szachniewicz2023existential}, the same holds for arithmetic divisors. Hence there is a canonical map from $\LdivQ(\Field)$ to its arithmetic ``concrete'' realisation and it is an isomorphism by construction.
\end{proof}

\section{Heights}

In this section we construct the correspondence between all the introduced structures. We start by defining globally valued fields via local terms.

\begin{definition}\label{definition_GVF_Globally_valued_fields_by_predicates}    
    \textit{Local terms} on $\Field$ are by definition a collection of functions $R_t:(\Field^{\times})^{n_t} \to \RR$ indexed by $\QQ$-tropical polynomials $t$ (of arity $n_t$ respectively), that are also denoted by
    \[ R_t(a) =: \int t(v(a)) dv, \]
    and which satisfy the following axioms:
    \begin{enumerate}
        \item $R_t$ are compatible with permutations of variables and dummy variables.
        \item (Linearity) $R_{t_1+t_2} = R_{t_1} + R_{t_2}, R_{q t} = q R_t$ for all $\QQ$-tropical polynomials $t_1, t_2, t$ and rational numbers $q$.
        \item (Local-global positivity) If $(\phi, t)$ is a positive pair (see Definition \ref{definition:positive_pairs}), then
        \[ (\forall a)(\phi(a) \implies  R_t(a) \geqslant 0). \]
    \end{enumerate}
    We use the name \textit{product formula field} to call a field whose local terms satisfy the following axiom:
    \begin{enumerate}
        \item[(4)] (Product formula) $\int v(x) dv = 0$ for all $x \in \Field^{\times}$.
    \end{enumerate}
    Fields with local terms (resp. product formula fields) form a category where morphisms are embeddings of fields such that the predicates $R_t$ of the bigger field restrict to the ones of the smaller one. We call local terms \textit{trivial} if all predicates $R_t$ are zero.
\end{definition}

\begin{lemma} \label{lemma_admissible_measures_define_MVF_structures}
    An admissible measure $\mu$ on $\Field$ defines local terms on $\Field$ via
    \[ R_t(a) := \int_{\Omega_{\Field}} t(v(a)) d\mu(v). \]
    Moreover, these only depend on the renormalisation class of $\mu$.
\end{lemma}
\begin{proof}
    Note that the map $v \mapsto t(v(a))$ can be bounded by a multiple of the maximum of the valuations of finitely many elements from $\Field^{\times}$, by Lemma \ref{lemma_tropical_polynomials_are_differences_of_maxima}. Hence, the integral will be finite by the first condition in the definition of admissibility.
    
    We check the local-global positivity axiom. If $(\phi,t)$ is a positive pair and $\Field \models \phi(a)$, then $v \mapsto t(v(a))$ in non-negative on $\Omega_{\Field}$, so the integral will be non-negative. We skip checking the other axioms of multiply valued field structures.

    Assume that $\mu, \mu'$ are in the same renormalisation class and keep the notation from Definition \ref{definition_renormalisation_of_an_admissible_measure_new}. Let $f$ be a function supported on $\pi^{-1}(\oOF(a))$ for some $a \in \Field^{\times}$ which is homogeneous with respect to the $\RR_{>0}$-action, i.e., $f(rv)=rf(v)$ for $r \in \RR_{>0}$. Then we have
    \[ \int_{\OF} f(v) d\mu(v) = \int_{\oOF(a)} f(s_a(\ov{v})) d\mu_a(\ov{v}) = 
    \int_{\oOF(a)} f(s_a(\ov{v})) \frac{d\mu_a}{d\mu'_a} d\mu'_a(\ov{v}) \]
    \[ =
    \int_{\oOF(a)} f(s_a(\ov{v})) \frac{s'_a}{s_a} d\mu'_a(\ov{v})
    = \int_{\oOF(a)} f(s'_a(\ov{v})) d\mu'_a(\ov{v})
    = \int_{\OF} f(v) d\mu'(v). \]
    Now let $f$ be a homogeneous (without loss of generality non-negative) function on $\OF$ that is both $\mu, \mu'$-integrable. Take $a_1, \dots, a_n \in \Field^{\times}$ such that the integrals with respect to $\mu, \mu'$ are approximated up to $\varepsilon$ by the integrals over $\bigcup_{i=1}^n \OF(a_i)$. The above calculation and the inclusion–exclusion principle shows that the integrals of $f$ with respect to $\mu, \mu'$ differ only at most $2 \varepsilon$, so they must be the same, as $\varepsilon$ was arbitrary. If we take $f = \ev_{\beta}$ for some $\beta \in \LdivQ(\Field)$, we are done.
\end{proof}

Note that the first part of the above lemma works, even if $\mu$ satisfies only the first condition of being an admissible measure on $\Field$. Before we proceed further, we need to prove a few properties of heights (see Definition \ref{definition:heights:intro}) which will allow us to pass from heights to positive functionals.

\begin{lemma} \label{lemma:height_invariance_duplication_coordinates}
    Let $h$ be a height on $F$. Then, $h$ is invariant under duplication of coordinates, addition of zero coordinates, and multiplication of coordinates by roots of unity. More precisely, let $\xi \in F^n$, $x\in F$, and let $\lambda\in F^\times$ be a root of unity. Then,
    \[
        h(\xi,x,x) = h(\xi,x,0) = h(\xi,\lambda x) = h(\xi,x).
    \]
\end{lemma}

\begin{proof}
    Let us start with duplication of coordinates. Using the fact that $h(1,1) = 0$, we have
    \[
        h(\xi,x,\xi,x) = h((\xi,x)\otimes(1,1)) = h(\xi,x).
    \]
    But, by monotonicity, $h(\xi,x) \leqslant h(\xi,x,x) \leqslant h(\xi,x,\xi,x)$, which shows that $h(\xi,x,x) = h(\xi,x)$.
    Now, we show that $h(1,0) = h(1,\lambda) = 0$.
    \begin{itemize}
        \item By duplication of coordinates, $h(1,0) = h(1,0,0,0) = h((1,0)^{\otimes 2}) = 2h(1,0)$, thus $h(1,0) = 0$.
        \item Let $m$ be the order of $\lambda$ in $F^{\times}$, and consider the tuple $\Lambda := (1,\lambda,\ldots,\lambda^{m-1})$. Since we know $h(1,1) = 0$, we may assume that $m\geqslant 2$. Then, up to permutation of coordinates, $\Lambda^{\otimes m} = (\Lambda,\ldots,\Lambda)$, where the tuple $\Lambda$ is repeated $m^{m-1}$ times. Hence, by duplication of coordinates, $mh(\Lambda) = h(\Lambda)$, thus $h(\Lambda) = 0$. But, by monotonicity, $0 = h(1) \leqslant h(1,\lambda) \leqslant h(\Lambda)$, which shows that $h(1,\lambda) = 0$.
    \end{itemize}
    Now, by monotonicity of $h$, we have
    \[
        h(\xi,x) \leqslant h(\xi,x,0) \leqslant h((\xi,x)\otimes(1,0)) = h(\xi,x) + h(1,0),
    \]
    and
    \[
        h(\xi,x) \leqslant h(\xi,x,\lambda x) \leqslant h((\xi,x)\otimes(1,\lambda)) = h(\xi,x) + h(1,\lambda).
    \]
    Therefore, $h(\xi,x,0) = h(\xi,x,\lambda x) = h(\xi,x)$.
\end{proof}

\begin{lemma}[Generalized triangle inequality] \label{lemma:generalized_triangle_inequality:section_heights}
    Let $h$ be a height on $F$ with Archimedean error $\aerror$. Then, for all $n,m\in\mathbb{N}$, and $(x_1,\ldots,x_n)\in(F^m)^n$, we have
    \[
        h(x_1+\ldots+x_n) \leqslant h(x_1,\ldots,x_n) + \log_2(n) \aerror.
    \]
\end{lemma}

\begin{proof}
    First, we show that if the lemma is true for some $n \in \mathbb{N}$, it is also true for $2n$. Indeed, we then have
    \begin{align*}
        h(x_1+\ldots+x_{2n}) &\leqslant h(x_1 + \ldots + x_n , x_{n+1} + \ldots + x_{2n}) + \aerror \\
            &= h((x_1,x_{n+1}) + (x_2,x_{n+2}) + \ldots + (x_n,x_{2n})) + \aerror \\ &\leqslant h(x_1,\ldots,x_{2n}) + (\log_2(n)+1)\aerror \\
        h(x_1 + \ldots + x_{2n}) &\leqslant h(x_1,\ldots,x_{2n}) + \log_2(2n)\aerror.
    \end{align*}
    By induction, the lemma is true whenever $n$ is a power of $2$.
    
    In the general case, monotonicity allows us to add any number of zero tuples, so
    \[
        h(x_1 + \ldots + x_n) \leqslant h(x_1 , \ldots, x_n) + \lceil \log_2(n) \rceil\aerror.
    \]
    Now, fix some $n,m\in\mathbb{N}$, and tuples $x_1,\ldots,x_n\in \Field^m$. For all $p\in\mathbb{N}$, by distributivity of the Segre product relative to coordinate-by-coordinate addition, the $p$-th Segre power of the sum $x_1+\ldots+x_n$ can be written
    \[
        (x_1+\ldots+x_n)^{\otimes p} = \sum_{1\leqslant i_1,\ldots,i_p\leqslant n} x_{i_1}\otimes \ldots \otimes x_{i_p}.
    \]
    On the other hand, the $p$-th Segre power is also distributive relative to concatenation, hence (up to permutation of coordinates), $(x_1,\ldots,x_n)^{\otimes p}$ is the concatenation of the tuples $x_{i_1}\otimes \ldots \otimes x_{i_p}$ for all $(i_1,\ldots,i_p)\in \{1,\ldots,n\}^p$. 
    Thus,
    \[
         h((x_1+\ldots+x_n)^{\otimes p}) \leqslant h((x_1,\ldots,x_n)^{\otimes p}) + \lceil \log_2(n^p)\rceil \aerror,
    \]
    i.e., $p\cdot h(x_1 + \ldots + x_n) \leqslant p \cdot h(x_1,\ldots,x_n) + \lceil p\log_2(n)\rceil\aerror$. We conclude by dividing by $p$ as $p$ goes to infinity.
\end{proof}

\begin{proposition} \label{proposition:equivalence_positive_functional_height:section_heights}
    Let $l : \LdivQ(F) \rightarrow \mathbb{R}$ be a positive functional. Define a height $h_l$ on $F$ by
    \[
        \forall(x_1,\ldots,x_n)\in F^n\setminus\{0\},\, h_l(x_1,\ldots,x_n) = -l\left( \bigmeet_{\substack{1\leqslant i\leqslant n\\ x_i\neq 0}} \adiv(x_i) \right).
    \]
    Then, $l\mapsto h_l$ is a bijection from the set of positive functionals on $F$ to the set of heights on $F$, which additionally restricts to a bijection from the set of global functionals to the set of global heights.
\end{proposition}

\begin{proof}
    It is immediate that $h_l$ satisfies the axioms in Definition \ref{definition:heights:intro}, with Archimedean error $\aerror = h_l(2,1)$, since the lattice divisor $\adiv(x+y) - \adiv(x) \meet \adiv(y) - \adiv(2) \meet 0$ is positive.
    
    Now, for a height $h$ on $F$, we may recover a group morphism $l_h : \ULatQ(F)\rightarrow \mathbb{R}$ by setting
    \[
        l_h\left(\bigmeet_{j=1}^m\adiv(b_j)-\bigmeet_{i=1}^n\adiv(a_i)\right) = h(a) - h(b),
    \]
    for all $a\in (\Field^\times)^n, b\in (\Field^\times)^m$. This shows injectivity. To prove surjectivity, we need to show that $l_h$ factors to a positive functional on $F$, i.e., that it is positive on the positive cone of $\ULatQ(F)$. For that, we use characterisation $\textnormal{(iii)}$ in Corollary \ref{corollary:tropical_divisors_positivity_conditions}. More precisely, let us show that for all $\varepsilon>0$, if $\alpha\in \Gamma_\varepsilon$, then $l_h(\alpha)\geqslant -\varepsilon \aerror$, where $\aerror$ is the Archimedean error of $h$. By Definition \ref{definition:gamma_epsilon}, we may write
    \[
        \alpha = \bigmeet_{j=1}^m \adiv(b_j) - \bigmeet_{i=1}^n \adiv(a_i),
    \]
    where $a = (a_1,\ldots,a_n) \in (F^\times)^n$, and each $b_j$ is a sum of exactly $\lfloor 2^\varepsilon\rfloor$ elements of the form $a_i$, $-a_i$ or $0$, i.e., coordinates of the tuple $a' = (((-1,1)\otimes a),0)\in F^{2n+1}$. Note that by Lemma \ref{lemma:height_invariance_duplication_coordinates}, $h(a') = h(a)$.
    
    Let $r := \lfloor 2^\varepsilon \rfloor$, and fix a bijection $\pi : \{1,\ldots,2n+1\}^r \rightarrow \{1,\ldots,(2n+1)^r\}$. For all $j\in\{1,\ldots,r\}$ consider the tuple $d_j\in F^{(2n+1)^r}$ such that, for all $(i_1,\ldots,i_r)\in\{1,\ldots,(2n+1)\}^r$,
    \[
        d_{j,\pi(i_1,\ldots,i_r)} = a'_{i_j}.
    \]
    Since $d_j$ only contains coordinates of $a'$, we have $h(d_1,\ldots,d_r) \leqslant h(a') = h(a)$. On the other hand, every $b_j$ is a sum of $r$ coordinates of $a'$, so by construction, $b$ is contained in $d_1+\ldots+d_r$. Thus, by monotonicity and Lemma \ref{lemma:generalized_triangle_inequality:section_heights},
    \begin{align*}
        h(b) &\leqslant h(d_1+\ldots+d_r) \\
                 &\leqslant h(d_1,\ldots,d_r) + \log_2(r)\aerror \\
        h(b) &\leqslant h(a) + \varepsilon \aerror.
    \end{align*}
    
    Thus, $l_h(\alpha)\geqslant - \varepsilon \aerror$. By additivity of $l_h$ and characterisation $\textnormal{(iii)}$ in Corollary \ref{corollary:tropical_divisors_positivity_conditions}, it follows that $l_h$ is positive on the positive cone of $\ULatQ(F)$, which concludes the proof. Furthermore, it is immediate that $h_l$ is global if and only if $l$ is a global functional.
\end{proof}

\begin{remark}
    The above proof also shows that for any height $h$, the Archimedean error $\aerror$ can always be taken to be $h(2,1)$.
\end{remark}

The main result of this section is the following theorem giving equivalence between the above introduced definitions.
How to pass between various equivalent structures mentioned in the theorem is explained in its proof.

\begin{theorem} \label{theorem_equivalent_MVF_structure_definitions}
    There is an explicit bijective correspondence between the following structures on $\Field$:
    \begin{enumerate}
        \item Local terms,
        \item Heights,
        \item Positive functionals,
        \item Families of local measures,
        \item Renormalisation classes of admissible measures,
        \item Equivalence classes of lattice valuations.
    \end{enumerate}
    Moreover, this correspondence respects the `global' ones. 
\end{theorem}
\begin{proof}
    We describe how to pass between structures with consecutive numbers in circle.
    \begin{itemize}
        \item To pass from (1) to (2), we define $h(a_1,\dots,a_n)$ to be $R_t(a_1,\dots,a_n)$ for $t = -\min(x_1,\dots,x_n)$. The axioms from Definition \ref{definition:heights:intro} are satisfied with Archimedean error $\aerror = h(2,1)$, since the lattice divisor $\adiv(x+y) - \adiv(x) \meet \adiv(y) - \adiv(2) \meet 0$ is positive.
        
        \item From (2) to (3) (even to get the equivalence between those two) we use Proposition \ref{proposition:equivalence_positive_functional_height:section_heights}.

        \item From (3) to (4), a positive functional on $\Field$ yields a family of local measures $(\mu_a)_{a \in \Field^{\times}}$ as in Corollary \ref{corollary_positive_functional_gives_local_measures}.

        \item The passage from (4) to (5) is the content of Construction \ref{construction_admissible_measure_from_local_measures}.

        \item Lemma \ref{lemma_admissible_measures_define_MVF_structures} provides a correspondence from (5) to (1).

        \item From (5) to (6), let $\mu$ be an admissible measure on $F$. This makes $L^1(\Omega_\Field,\mu)$ an $L^1$-lattice, and the map $\underline{v} : x \mapsto (v \mapsto v(x))$ a lattice valuation (here $\Gamma$ is the $L^1$-lattice generated by the image of $\underline{v}$).

        \item Finally, from (6) to (2), if $(\Gamma,\underline{v})$ is a lattice valuation on $\Field$, then we define a height on $\Field$ by
        \[
            h(a_1,\ldots,a_n) = -\int\underline{v}(a_1)\meet\ldots\meet \underline{v}(a_n).
        \]
        This clearly satisfies the axioms from Definition \ref{definition:heights:intro}.
    \end{itemize}
    It is obvious that going in a cycle, we get back the same structure.
\end{proof}

\begin{definition} \label{definition_GVF_ultimate_definition}
    We say that $\Field$ is equipped with a \textit{multiply valued field} structure (abbreviated MVF) if it is given with one of the equivalent structures from Theorem \ref{theorem_equivalent_MVF_structure_definitions}.
    We call an MVF $\Field$ a \textit{globally valued field} (abbreviated GVF) if the MVF structure is additionally global, or in other words it satisfies the product formula.
\end{definition}

\begin{corollary} \label{corollary_measures_on_OF_can_be_changed_to_admissible_ones_inducing_same_MVF}
    If $\mu$ is a Borel measure on $\OF$ satisfying the first axiom of being an admissible measure on $\Field$, then there is an admissible measure $\nu$ on $\Field$ inducing the same MVF structure on $\Field$. Hence, an MVF structure on $\Field$ can be seen as a Borel measure on $\OF$ such that the evaluation functions are integrable, up to the equivalence relation given by integrating homogeneous functions on $\OF$ as in Lemma \ref{lemma_admissible_measures_define_MVF_structures}.
\end{corollary}

Finally, we can compare globally valued fields to proper adelic curves, which we define below.

\begin{definition} \label{definition_adelic_curves} \cite[Section 3.1]{Adelic_curves_1}
    An \textit{adelic structure on $\Field$} is a measure space $(\Omega, \cA, \nu)$ together with a map $\phi: \omega \mapsto |\cdot|_\omega$ from $\Omega$ to $M_{\Field}$ such that for all $a \in \Field^{\times}$ the function
    \[ \omega \mapsto \log|a|_\omega \]
    is integrable (and measurable) on $\Omega$. Moreover, we call the adelic structure \textit{proper} if for all $a \in \Field^{\times}$
    \[ \int_{\Omega} \log|a|_{\omega} d\nu(\omega) = 0. \]
    Together, the data $(\Field, (\Omega, \cA, \nu), \phi)$ is called an \textit{adelic curve}.
\end{definition}

Similarly as in Lemma \ref{lemma_admissible_measures_define_MVF_structures}, adelic structures on $\Field$ induce multiply valued field structures by taking
\[ R_t(a) := \int_{\Omega} t(-\log|a|_{\omega}) d\nu(\omega). \]
Moreover, by Remark \ref{remark_GVF_measures_can_be_pushed_from_M_F}, we get the following.
\begin{corollary} \label{corollary_every_MVF_structure_comes_from_an_adelic_curve}
    If $\Field$ is countable, then every multiply valued field structure on $\Field$ is induced by an adelic structure. Similarly, in that case every GVF structure on $\Field$ is induced by a proper adelic structure. Moreover, $\Omega$ can be assumed to be $M_{\Field}$.
\end{corollary}

\section{Intersection theory} \label{section_representations}

One can construct examples of globally valued fields using intersection theory or Arakelov intersection theory. If $\cX$ is an arithmetic variety, let us denote by $\Adiv_{\RR}(\cX)$ the group of arithmetic $\RR$-divisors of $C^0$-type, see \cite{szachniewicz2023existential} for precise definitions. Assume that $\Field$ is a finitely generated extension of $\QQ$. Let $\Adiv_{\RR}(\Field)$ be the direct limit of $\Adiv_{\RR}(\cX)$'s, over the system of arithmetic varieties $\cX$ with their function field isomorphic to $\Field$. A GVF functional $l:\Adiv_{\RR}(\Field) \to \RR$ is a linear map which vanishes on principal arithmetic divisors and is non-negative on effective arithmetic divisors. With this notation we have the following.

\begin{theorem} \cite[Theorem C]{szachniewicz2023existential} \label{theorem_bijection_GVF_functionals_local_terms}
    There is a natural bijection:
    \[ \bigl\{ \textnormal{GVF functionals $\Adiv_{\RR}(\Field) \to \RR$}  \bigr\} \longleftrightarrow \bigl\{ \textnormal{GVF structures on $F$}  \bigr\} \]
    \[ l \longmapsto \bigl( R_t^l(a):=l(\ov{\cD}_t(a)) \bigr)_{t} \]
    where $\ov{\cD}_t(a)=t(\adiv(a))$ if $a$ is a tuple from $\Field^{\times}$ and $t$ is a $\QQ$-tropical polynomial.
\end{theorem}

The lattice structure on $\Adiv_{\QQ}(\Field)$ was constructed in \cite[Section 3]{szachniewicz2023existential} and $\ov{\cD}_t(a)$ is defined with respect to it. This result also follows from Theorem \ref{theorem_equivalent_MVF_structure_definitions} and Corollary \ref{corollary_description_of_LDiv_Arakelov_divisors}, together with density of lattice divisors in the space of arithmetic divisors from \cite[Section 3]{szachniewicz2023existential}.

In particular, if $\cX$ is an $\Field$-model of dimension $d$ over $\Spec(\ZZ)$ and $\ov{\cD}$ is an arithmetically nef $\RR$-divisor on $\cX$, then the map
\begin{equation}\label{eq:polarisation_functional}
    \ov{\cE} \mapsto \ov{\cD}^d \cdot \ov{\cE}
\end{equation}
defines a GVF functional $\Adiv_{\RR}(F) \rightarrow \mathbb{R}$, where the product taken is the arithmetic intersection product (see e.g. \cite{Ikoma_Concavity_of_arithmetic_volume}), possibly calculated on a blowup of $\cX$. In fact, any GVF structure can be approximated by such functionals in the following sense.

\begin{theorem}\cite{szachniewicz2023existential} \label{theorem_density_GVF_arithmetic_case}
    The set of GVF structures on $\Field$ coming from arithmetically ample arithmetic $\RR$-divisors $\ov{\cD}$ on $\Field$-models $\cX$ by the formula (\ref{eq:polarisation_functional}) is dense in the space of all GVF functionals $\Adiv_{\QQ}(\Field) \to \RR$, with respect to the pointwise convergence topology.
\end{theorem}

Let us remark that the example outlined in Equation (\ref{eq:three}) from the introduction comes from arithmetic intersection product by taking $\cX$ to be $\PP_{\ZZ}^1$ and $\ov{\cD}$ to be $\adiv(z) \vee 0$ (equivalently $\cO(1)$ with the Weil metric), where we write $z=\frac{x}{y}$ for the variable generating the function field $\QQ(z)$ of $\PP_{\ZZ}^1$.

Similarly, if $K$ is a finitely generated field over a field $k$ and we restrict our attention to GVF structures that are trivial when restricted to $k$, we get the following correspondence.

\begin{theorem} \cite{GVF2} \label{theorem_GVF_structures_trivial_on_the_base_field}
    There is a natural bijection:
    \[ \bigl\{ \textnormal{GVF functionals $\Div_{\QQ}(K/k) \to \RR$}  \bigr\} \longleftrightarrow \bigl\{ \textnormal{GVF structures on $K$ trivial on $k$}  \bigr\}, \]
    \[ l \longmapsto \bigl( R_t^l(a):=l(D_t(a)) \bigr)_{t} \]
    where $D_t(a) = t(\cdiv(a))$ for tuples $a$ from $\Field^{\times}$ and $\QQ$-tropical polynomials $t$.
\end{theorem}

Here $\Div_{\QQ}(K/k) = \varinjlim \Div_{\QQ}(X)$ is the direct limit of ($\QQ$-Cartier) divisors of normal projective varieties $X$ over $k$ together with isomorphisms $k(X) \simeq K$ over $k$. These are often called ($\QQ$-Cartier) b-divisors, see \cite{shokurov2003prelimiting}. The lattice structure on $\Div_{\QQ}(K/k)$ is defined as follows: if $D, E \in \Div(X)$ are effective Cartier divisors, then we define $D \wedge E$ by the class of the exceptional fiber of the blowup $X' \to X$ of $X$ at $D \cap E$. Equivalently, we define $D \wedge E$ as the class of $D \cap E$ on a blowup $X'$ where $D \cap E$ is a Cartier divisor. For $\QQ$-divisors, or non-effective divisors, we extend using axioms of a lattice group. The proof that it is well defined is the same as in \cite{szachniewicz2023existential}. A GVF functional $l:\Div_{\QQ}(K/k) \to \RR$ on the left hand side of the correspondence is a linear map which vanishes on principal divisors and is non-negative on the effective cone. Let us note, that there is a natural pairing $\Val_{K/k} \times \Div_{\QQ}(K/k) \to \RR$ satisfying the conclusion of Lemma \ref{lemma_nice_properties_of_evaluation}, and such that for $D \in \Div_{\QQ}(K/k)$ we have $D=0$ if and only if, for all $v \in \Val_{K/k}$ we have $v(D) = 0$.

\begin{proof}[Proof of Theorem \ref{theorem_GVF_structures_trivial_on_the_base_field}]
    By Theorem \ref{theorem_equivalent_MVF_structure_definitions} the data of a GVF structure on $K$ trivial on $k$ is equivalent to the data of a GVF functional $l:\LdivQ(K) \to \RR$ that vanishes on $\LdivQ(k)$. Let $I_k \subset \LdivQ(K)$ be the ideal generated by $\LdivQ(k)$ (in the sense of Definition \ref{definition_ideal_in_a_lattice}). It is enough to prove that as lattice groups
    \[ \LdivQ(K)/I_k \simeq \Div_{\QQ}(K/k). \]
    We construct the isomorphism in the following way. First, there is a map $f:\ULatQ(K) \to \Div_{\QQ}(K/k)$ coming from the universal property. Second, note that the ideal $I$ from Definition \ref{definition_lattice_of_a_field} is contained in the kernel of $f$. Indeed, if $\alpha \in I$, then $f(\alpha)$ pairs with any $v \in \Val_{K/k}$ to zero, so $f(\alpha) = 0$. Thus, $f$ factors through $f:\LdivQ(K) \to \Div_{\QQ}(K/k)$, and from Corollary \ref{corollary_ideal_vanishing_relative_situation} we get that the kernel is exactly the ideal generated by $\LdivQ(k)$, which finishes the proof.
\end{proof}

From the above proof we get the following.

\begin{corollary} \label{corollary_quotient_of_lattice_divisors_relative_case}
    The quotient of $\LdivQ(K)$ by the (lattice group) ideal generated by $\LdivQ(k)$ is isomorphic to $\Div_{\QQ}(K/k)$.
\end{corollary}

In fact, GVF functionals also vanish on numerically trivial divisors. We give a simple proof below.

\begin{lemma} \cite[Proposition 11.5]{GVF2}
    Let $l:\Div(K/k) \to \RR$ be a GVF functional. Then for a numerically trivial divisor $E$, the equality $l(E) = 0$ holds.
\end{lemma}
\begin{proof}
    Assume that $E$ is a divisor on $X$ with $k(X) \simeq K$ and fix an ample $\QQ$-divisor $A$ on $X$. Since $E$ is numerically trivial, $\QQ$-divisors $A \pm E$ are ample (by Kleiman's criterion). Hence, for some big enough $n$ we get that $n(A \pm E)$ are rationally equivalent to effective $\QQ$-divisors. Thus
    \[ n \cdot l(A \pm E) = l(n(A \pm E)) \geqslant 0. \]
    Dividing by $n$, we get $|l(E)| \leqslant l(A)$. Since $A$ was an arbitrary ample $\QQ$-divisor, by replacing it by $tA$ for arbitrarily small $t>0$, we get $l(E) = 0$.
\end{proof}
Hence, the left (and so also the right) hand side of the correspondence in Theorem \ref{theorem_GVF_structures_trivial_on_the_base_field} can be equivalently described as
\[ N_1^{+}(K/k) := \varprojlim N_1^{+}(X), \]
where $N_1^{+}(X)$ is the movable cone in $N_1(X)$ (i.e., the dual cone to the effective cone in $N^1(X)$). Let us point out that the cone $N_1^{+}(K/k)$ has appeared in \cite{dang2022intersection} under the name $\BPF^{d-1}(\cX)$ where $d$ is the transcendence degree of $K$ over $k$. 

Note that if $D$ is an ample divisor on a $d$-dimensional $X$ with $k(X) \simeq K$, then the formula
\[ E \mapsto D^{d-1} \cdot E \]
defines a GVF structure on $K$ and the set of such GVF structures on $K$ (trivial on $k$) is dense. 

\begin{theorem} \cite[Corollary 10.11]{GVF2} \label{theorem_density_GVF_geometric_case}
    The set of GVF structures on $K$ coming from ample divisors $D$ on varieties $X$ with $k(X) \simeq K$ (by the formula $E \mapsto D^{d-1} \cdot E$) is dense in the space of all GVF functionals $\Div_{\QQ}(K/k) \to \RR$, with respect to the pointwise convergence topology.
\end{theorem}

In a different language, this was independently proven in \cite[Theorem C]{dang2022intersection}.

\section{Disintegration} \label{section_disintegration}

Consider a field extension $K \subset L$ and assume that $L$ is equipped with an MVF structure. One can describe the restriction of this structure to $K$ using any of the equivalent notions from Theorem \ref{theorem_equivalent_MVF_structure_definitions}. For example, if $(\mu'_a)_{a \in L^{\times}}$ is a family of local measures on $L$, then for $a \in \Field$, one can look at restrictions $\Omega_L(a) \to \Omega_K(a)$ and if $\mu_a$ denotes the pushforward of $\mu'_a$, then $(\mu_a)_{a \in K^{\times}}$ is the family of local measures on $K$ corresponding to the restricted multiply valued field structure. Similarly, if $\mu'$ is a measure on $\Omega_L$ inducing the MVF structure on $L$, then for the partial map $r:\Omega_L \to \Omega_K$ restricting elements to $K$, the measure $\mu = r_{*}\mu'$ induces the restricted MVF structure on $K$. Also, one can look at the restriction of the corresponding positive functional, via the embedding from Corollary \ref{corollary_embedding_of_fields_induces_embedding_of_LDivs}. Finally, the restriction of the MVF or height functions is just the restriction of their domains to the smaller field.

Let us introduce some notation regarding the partial restriction map $r:\Omega_L \to \Omega_K$. It is defined on
\[ \Omega_{L:K} := \{ v \in \Omega_L : v|_K \not\in \{-\infty,0,\infty\}^K \}. \]
We denote the complement of this open set by $\Omega_{L/K}$. Also, for $v \in \Omega_K$ we write $\Omega_{L:v}$ for the set of $w \in \Omega_L$ with $w|_K = v$. Hence
\[ \Omega_{L:K} = \bigcup_{v \in \Omega_K} \Omega_{L:v}. \]
Let $\mu$ be an admissible measure on $L$. We can write $\mu$ as the sum of its restrictions to $\Omega_{L:K}, \Omega_{L/K}$. We denote these restrictions by $\mu_{L:K}, \mu_{L/K}$ respectively. Note that $\mu = \mu_{L:K} + \mu_{L/K}$.
\begin{lemma} \label{lemma_renormalising_relatively_to_an_extension}
    If $r_*\mu$ can be renormalised to an admissible measure $\nu$ on $K$, then $\mu$ can be renormalised to an admissible measure $\mu'$ on $L$ such that $r_* \mu' = \nu$. 
\end{lemma}
\begin{proof}
    Note that the local measures associated to $\nu$ and the pushforwards of the local measures associated to $\mu$ are the same. Hence, given sections $s_a:\ov{\Omega}_K(a) \to \Omega_K(a)$ for $a \in K^{\times}$ agreeing on overlaps, it is enough to construct sections $u_a:\ov{\Omega}_L(a) \to \Omega_L(a)$ (agreeing on overlaps) such that the diagram
    \[\begin{tikzcd}
    	{\ov{\Omega}_L(a)} & {\Omega_L(a)} \\
    	{\ov{\Omega}_K(a)} & {\Omega_K(a)}
    	\arrow["{s_a}", from=2-1, to=2-2]
    	\arrow["{u_a}", from=1-1, to=1-2]
    	\arrow["r", from=1-2, to=2-2]
    	\arrow["r"', from=1-1, to=2-1]
    \end{tikzcd}\]
    commutes. This is by copying Construction \ref{construction_admissible_measure_from_local_measures} for both $K$ and $L$. Start by taking any sections $u_a$ agreeing on overlaps (construct them as in Construction \ref{construction_admissible_measure_from_local_measures}). Next, let $h_a:\ov{\Omega}_L(a) \to \RR_{>0}$ be defined by
    \[ h_a(\ov{v}) = \frac{s_a(\ov{v}|_K)}{u_a(\ov{v})|_K}. \]
    Since both $s_a$'s and $u_a$'s agree on overlaps, we get that $h_a$'s also agree on overlaps, so by putting
    \[ u'_a := h_a \cdot u_a \]
    we are done.
\end{proof}

Thus, if $\nu$ is a global measure on $K$ inducing a GVF structure, and there is a GVF structure on $L$ extending the one on $K$, one can find a global measure $\mu$ on $L$ such that $r_*\mu = \nu$. 

Let us restrict our attention to the case where $K, L$ are countable. Note that then $\Omega_K, \Omega_L$ are locally compact, separable metric spaces (because a compact metric space is separable). Using a disintegration theorem, we get the following description of GVF extensions.

\begin{corollary} \label{corollary_disintegration_of_a_GVF_measure}
    Let $\nu$ be a global measure on $K$. Then for every GVF structure on $L$ extending the one on $K$, there is a measure $\mu_{L/K}$ on $\Omega_{L/K}$ and a family of measures $\mu_v$ on $\Omega_{L:K}$, for every $v \in \Omega_K$, such that:
    \begin{enumerate}
        \item $v \mapsto \mu_v$ is a Borel map and $\mu_v$ is a probability measure for $\nu$-almost every element of $\Omega_K$;
        \item if $\mu_{L:K} := \nu \otimes \mu_v$, that is $\mu(A) = \int_{\Omega_{K}} \mu_v(A) d\nu(v)$, then $\mu := \mu_{L:K} + \mu_{L/K}$ is a global measure on $L$ inducing the GVF structure there;
        \item $\mu_v$ is concentrated on $r^{-1}(v)$ for $\nu$-almost all $v \in \Omega_K$.
    \end{enumerate}
\end{corollary}
\begin{proof}
    First we use Lemma \ref{lemma_renormalising_relatively_to_an_extension} to get a global measure $\mu$ on $L$ with $r_*\mu_{L:K}=\nu$. Next, we order $K^{\times}$ as $(a_i)_{i<\omega}$ and for every $a_i$ we use the disintegration theorem \cite[Theorem A]{possobon2023geometric} to get the measures $\mu_{v}$ for $v \in \pi^{-1}(\ov{\Omega}_K(a_i) \setminus \bigcup_{j<i} \ov{\Omega}_K(a_j))$ satisfying the assumptions locally (where $\pi:\Omega_K \to \ov{\Omega}_K$ is the projection map). By countability of $K$, the obtained family $(\mu_v)_{v \in \Omega_{K}}$ satisfies the required conditions.
\end{proof}

\begin{definition}
    We call a GVF structure on a field $K$ \textit{discrete} if it is induced by a measure $\mu = \sum_{i \in I} \delta_{v_i}$, i.e., a sum of Dirac deltas at some family of pairwise non-equivalent valuations $\{v_i\}_{i \in I} \subset \Val_K$ such that for every $a \in K^{\times}$ the set $\{i \in I: v_i(a) \neq 0\}$ is finite and we have $\sum_i v_i(a) = 0$.
\end{definition}

Corollary \ref{corollary_disintegration_of_a_GVF_measure} gives a nice presentation of globally valued fields structures on finitely generated extensions of a countable discrete GVF.

\begin{corollary} \label{corollary:GVFs_as_measure_on_Berkovich_spaces}
    Let $K \subset L$ be a finitely generated extension of fields, where $K$ is countable and carries a discrete GVF structure given by the family $\{v_i\}_{i \in I} \subset \Val_K$. Let $X$ be a projective variety over $K$ with function field isomorphic to $L$. Denote by $X_i^{\an}$ the Berkovich analytification of $X_i = X \otimes_{K} K_{v_i}$ for $i \in I$, and by $X_0^{\an}$ the Berkovich analytification of $X$ treated as a variety over $K$ with the trivial absolute value. Then, there exists a measure $\mu_0$ on $X_0^{\an}$ and a family of probability measures $\mu_i$ on $X_i^{\an}$ such that
    \[ \mu = \mu_0 + \sum_{i\in I} \mu_i \]
    is a global measure on $L$. Moreover, all these measures are concentrated on points in the analytifications mapping to the generic point $\eta:\Spec(L) \to X$.

    On the other hand, any such family of measures with $\mu$ satisfying the product formula axiom yields a GVF structure on $L$ extending the one on $K$.
\end{corollary}

\begin{remark}
    Note that the measure $\mu_0$ above is not unique, because in $X_0^{\an}$ there are many points corresponding to equivalent pseudo-valuations on $L$. This is in contrast to $X_i^{\an}$ for $i \in I$, where any two different points represent non-equivalent pseudo-valuations on $L$, because all of them restrict to the same non-trivial valuation on $K$.
\end{remark}

\section{Finite extensions} \label{section_finite_extensions}

In this section we study finite extensions of some globally valued fields.

\begin{lemma} \cite[Lemma 12.1]{GVF2} \label{lemma:unique_GVF_function_field_curve}
    Let $C$ be a curve over a field $k$ with the function field $K=k(C)$ and fix $t \in K$ transcendental over $k$. Then, for every $r \geqslant 0$ there is a unique GVF structure on $K$ which is trivial on $k$ and satisfies $\height(t)=r$.
\end{lemma}
\begin{proof}
    This follows from the fact that $N_1^{+}(K/k) = N_1^{+}(C) = \RR_{\geqslant 0}$. For an alternative proof, see \cite{GVF2}.
\end{proof}

We write $\GVFk[r]$ to denote the unique GVF structure on $\GVFk$ with $\height(t)=r$ and $\height(a)=0$ for $a \in k^{\times}$. Note that this notation depends on the choice of $t$.

\begin{lemma} \label{lemma:uniqueness_of_the_GVF_structure_on_number_fields}
    Let $K$ be a number field. Then, for any $r \geqslant 0$, there is a unique GVF structure on $K$ satisfying $\height(2) = r \cdot \log 2$.
\end{lemma}
\begin{proof}    
    Let $S \subset \Val_K$ be a set of representatives of equivalence classes of valuations on $K$ such that for all $a \in K^{\times}$ we have $\sum_{v \in S} v(a) = 0$. By the Ostrowski's theorem any GVF structure on $K$ is given by a choice of weights $w_v \geqslant 0$ such that for all $a \in K^{\times}$ we have $\sum_{v \in S} v(a) \cdot w_v = 0$. 

    Let $W$ be the $\RR$-vector space freely generated by $S$. A choice of weights as above corresponds to a linear functional on $W$ that vanishes at the image of the map $f:K^{\times} \to W$, sending $a$ to $\sum_{v \in S} v(a) \cdot v$. Let $N = \{ \sum_{v \in S} a_v \cdot v: \sum_{v \in S} a_v = 0 \}$ be the kernel of the `summing of coordinates' functional on $W$. It is enough to show that the real vector space spanned by $f(K^{\times})$ is $N$, as then, a linear functional $l:W \to \RR$ that vanishes on $f(K^{\times})$ is the same as a linear functional on $W/N \simeq \RR$.

    Consider $w = \sum_{v \in S} a_v \cdot v \in N$. 
    As $\Cl(K) \otimes \RR = 0$, the class of $w$ modulo the span of $f(K^{\times})$ can be represented by an element $w' \in N$ having zero entries over non-Archimedean places. By the Dirichlet unit theorem, the dimension of the real span of elements from $f(K^{\times})$ with vanishing non-Archimedean entries is of dimension $r_1+r_2-1$, where $r_1, r_2$ are the numbers of real and complex embeddings (up to conjugation) of $K$. Thus, the class of $w'$ modulo $f(K^{\times})$ is also trivial, which finishes the proof.
\end{proof}

If $r = 1$, we call this unique GVF structure (on a number field or $\GVFQ$) the \textit{standard} GVF structure. We write $\GVFQ[1]$ for the standard GVF structure on $\GVFQ$, and $\GVFQ[r]$ for the one satisfying $\height(2) = r \cdot \log 2$. 

For a general globally valued field, we only have uniqueness of a GVF structure on a finite extension additionally assuming Galois-invariance. To prove it we need a few lemmas.

\begin{lemma} \label{lemma:action_of_Galois_group_transitive}
    Let $K \subset \Field$ be a finite Galois extension with the Galois group $G$. Then the restriction map $\Omega^{\circ}_{\Field} \to \Omega^{\circ}_K$ is defined everywhere, and identifies $\Omega^{\circ}_K$ with the topological quotient of $\Omega^{\circ}_{\Field}$ by $G$.
\end{lemma}
\begin{proof}
    First we prove that the restriction map is defined everywhere. Assume that $|\cdot|$ is a pseudo-absolute value on $\Field$ that is $\{0,1,\infty\}$-valued on $K$. Let $$\mathcal{O} = \{x \in K : |x| < \infty \} \subset \mathcal{O}' = \{ y \in F : |y| < \infty \}$$ 
    and denote by $\kappa \subset \kappa'$ the residue field extension. Then $|\cdot|$ induces an absolute value on $\kappa'$ that is trivial on $\kappa$. By the fact that $\kappa \subset \kappa'$ is algebraic we conclude that $|\cdot|$ is $\{0,1,\infty\}$-valued on $F$.

    Now we prove that the restriction map is a topological quotient by $G$. By local compactness it is enough to prove that $G$ acts on fibers transitively. In other words, we need to show that if $|\cdot|_1, |\cdot|_2$ are pseudo-absolute values on $\Field$ that coincide on $K$, then there exists $g \in G$ such that $|g(\cdot)|_1 = |\cdot|_2$. 

    We split the problem in two cases, depending on $|\cdot|$ being non-Archimedean or Archimedean. First, assume the former. Then let $o_i \subset \cO_i$ denote the inclusions $\{|x|_i \leqslant 1\} \subset \{|x|_i \neq \infty \}$ for $i=1,2$. Note that $\cO_i$ is uniquely determined by $o_i$, as the ordered group $\Field^{\times}/o_i^{\times}$ has a unique convex non-zero Archimedean subgroup. By \cite[Conjugation Theorem 3.2.15]{engler2005valued} there is a $g \in G$ such that $g o_1 = o_2$. By uniqueness, we get that also $g \cO_1 = \cO_2$, and since both $|\cdot|_1, |\cdot|_2$ extend a non-trivial pseudo-absolute value $|\cdot|$ on $K$, we get that $|g^{-1}(\cdot)|_1 = |\cdot|_2$.
    
    Now, we assume that $|\cdot|$ is Archimedean. Let $\cO_1 = \{|x|_1 \neq \infty \}, \cO_2 = \{|x|_2 \neq \infty \}$ and denote by $\cO$ their intersections with $K$. By \cite[Conjugation Theorem 3.2.15]{engler2005valued} there is a $g \in G$ such that $g \cO_1 = \cO_2$. Let us assume that $\cO_1 = \cO_2$, by replacing $|\cdot|_1$ with $|g^{-1}(\cdot)|_1$. We end up with the diagram
    \[\begin{tikzcd}
    	\Field & {l \cup \{\infty\}} \\
    	K & {k \cup \{\infty\}}
    	\arrow[from=2-1, to=2-2]
    	\arrow[hook, from=2-1, to=1-1]
    	\arrow[from=1-1, to=1-2]
    	\arrow[hook, from=2-2, to=1-2]
    \end{tikzcd}\]
    where $k, l$ are the residue fields of $\cO, \cO_1$ respectively. Use the notation $|\cdot|_1, |\cdot|_2$ for the two induced absolute values on $l$ that extend the absolute value $|\cdot|$ on $k$. Denote by $G_1 < G$ the subgroup fixing $\cO_1$ setwise. By \cite[Lemma 5.2.6]{engler2005valued}, using the fact that $k \subset l$ is separable as we are in characteristic zero, the natural map $G_1 \to \Aut(l/k)$ is surjective. Thus, it is enough to show that there exists $\sigma \in \Aut(l/k)$ such that $|\sigma(\cdot)|_1 = |\cdot|_2$. We can assume that absolute values $|\cdot|_1, |\cdot|_2$ come from two embeddings of $l$ into $\CC$ that restrict to the same embedding on $k$. By \cite[Theorem 5.2.7.(2)]{engler2005valued} the extension $k \subset l$ is Galois, which proves the existence of the desired $\sigma \in \Aut(l/k)$ and finishes the proof.
\end{proof}

\begin{lemma} \label{lemma:extending_measure_from_the_quotient}
    Let $G$ be a finite group acting on a locally compact topological space $X$ with a locally compact topological quotient $r:X \to Y = X/G$. Assume that $\nu$ is a Radon measure on $Y$. Then there is a unique $G$-invariant Radon measure $\mu$ on $X$ such that $r_*\mu = \nu$.
\end{lemma}
\begin{proof}
    Since $\nu$ is a Radon measure, without loss of generality we can assume that $Y$ is compact. In that case $\nu$ is given by a functional $\nu:C(Y) \to \RR$ and we prove that it has a unique $G$-invariant extension to $C(X)$. Note that the pullback $C(Y) \to C(X)$ maps $C(Y)$ bijectively onto $C(X)^G$. Also, there is a retraction $\Tr:C(X) \to C(X)^G \simeq C(Y)$ defined by $\Tr(f) := \frac{1}{|G|}\sum_{g \in G} g^*f$. Thus, a $G$-invariant functional $\mu:C(X) \to \RR$ extending $\nu$ must be of the form $\mu(f) = \mu(\Tr(f)) = \nu(\Tr(f))$. This finishes the proof.
\end{proof}

Note that the above proof works in the same way for a compact group $G$ (using the Haar measure instead of the counting measure). Now we prove the existence and uniqueness of a Galois-invariant GVF extension. Existence can be concluded using \cite[Section 3.4]{Adelic_curves_1}.

\begin{proposition} \label{proposition:uniqueness_of_invariant_Galois_extension}
    Let $K$ be a GVF and consider a finite Galois field extension $K \subset \Field$. Then there is a unique GVF structure on $\Field$ extending the one on $K$ that is invariant under the action of $G=\Aut(\Field/K)$.
\end{proposition}
\begin{proof}
    By Lemma \ref{lemma:action_of_Galois_group_transitive} the restriction $r:\Omega^{\circ}_{\Field} \to \Omega^{\circ}_K$ satisfies the assumptions of Lemma \ref{lemma:extending_measure_from_the_quotient}, so for any global measure $\nu$ on $K$ there is a unique $G$-invariant measure $\mu$ on $\Field$ such that $r_*{\mu} = \nu$. By Lemma \ref{lemma_renormalising_relatively_to_an_extension} we get the uniqueness.

    To prove existence it is enough to prove that the measure $\mu$ obtained as above satisfies the product formula. From the proof of Lemma \ref{lemma:extending_measure_from_the_quotient} we get that for $a \in \Field^{\times}$ we have:
    \[ \int_{\Omega_{\Field}} \ev_a d\mu = \int_{\Omega_K} \Tr(\ev_a) d\nu, \]
    where $\ev_a(v) = v(a)$. However, using $\Tr(\ev_a) = \frac{1}{|G|}  \ev_{\Norm(a)}$ (where $\Norm(a) \in K$ is the $\Field/K$-norm) we get the product formula, provided that we prove integrability of $\ev_a$.

    Let $S_k$ be the symmetric functions for $k=0, \dots, n$ (with $n=[\Field:K]$) so that $\sum_{k\leqslant n} S_k(a)a^k = 0$. Note that for any pseudo-absolute value $|\cdot|$ on $\Field$ we have 
    \[ |a|^n \leqslant [n] \max_{k<n}(|S_k(a)|\cdot|a|^k) = [n] |S_{k_0}(a)|\cdot|a|^{k_0}, \]
    where by an expression in square brackets like $[n]$ we mean $n$ if $|\cdot|$ is Archimedean and $1$ otherwise. By dividing and taking $-\log$ we get
    \[ (n-k_0)v(a) \geqslant v(S_{k_0}(a)) -\log[n] \geqslant \min_{k<n} v(S_{k}(a)) -\log[n], \]
    where $v = - \log|\cdot|$. If we assume that $v(a) \leqslant 0$, we get that $v(a) \geqslant (n-k_0)v(a) \geqslant \min_{k<n} v(S_{k}(a)) -\log[n]$. Hence, this gives a bound
    \[-v(a)^- \leqslant -\min_{k<n} v(S_{k}(a)) + \log[n].\]
    By Remark \ref{remark_GVF_measures_can_be_pushed_from_M_F} we can assume that the restriction of $\nu$ to Archimedean $v \in \Omega_{K}$ is concentrated on the set $\{v(2) = -\log 2\}$, so we can replace $\log[n]$ in the above bound by $-v(n)^-$. In this way we get integrability of the function $-\ev_a^{-}$ and by replacing $a$ with $a^{-1}$ the integrability of $\ev_a$ follows.
\end{proof}

Since we have uniqueness and existence for every invariant finite Galois extension, there is also a unique Galois-invariant GVF extension to the algebraic closure. We call it the \textit{symmetric extension}. Note that without assuming Galois-invariance, there can be multiple ways to extend a GVF structure to a finite field extension.

\begin{example}
    Consider any GVF structure on $\QQ(x,y)$ with $\height(x) \neq \height(y)$. If we look at the degree two extension $\QQ(x+y,xy) \subset \QQ(x,y)$, then the symmetric extension of the restriction of the starting GVF structure satisfies $\height(x)=\height(y)$, so yields a different GVF structure than the starting one.
\end{example}

\section{Continuous logic} \label{section_continuous_logic}
We outline unbounded continuous logic, at least for discretely metrized spaces. We closely follow \cite{Ben_Yaacov_unbdd_cont_FOL}, but in this simplified setting. For other introductions to continuous logic, see \cite{yaacov2008model} or \cite{Hart_cont_logic}.

\begin{definition}
    A language $\mathcal{L}$ consists of the following data. First, a tuple of function symbols $f$, constant symbols $c$ and continuous predicate symbols $P$. Second, for each function or predicate symbol, its arity is specified, and it is equipped with a \textit{gauge modulus} $\Delta_p$, i.e., a right-continuous increasing function $\Delta_p:(0,\infty) \to (0, \infty)$. At last, there are two distinguished predicates, a \textit{gauge symbol} (or \textit{height symbol}) $\height$ of arity $1$, and the equality predicate $=$, of arity $2$. For simplicity, we assume that $\height$ is also a continuous predicate symbol.
\end{definition}

\begin{example}
    The language of globally valued fields consists of binary function symbols $+,\cdot$, constants $0, 1$, and continuous predicates $(h_n)_{n \in \NN}$. We often skip $n$ when discussing $h_n$. In Example \ref{example_GVF_language} we determine the gauge moduli for these symbols.
\end{example}

\begin{definition} \label{definition_L_structure_unbdd_cont_logic}
    Let $\mathcal{L}$ be a language. An $\mathcal{L}$-structure $\mathfrak{M}$ is a set $M$ together with interpretations of function, constant and predicate symbols from $\mathcal{L}$, and an interpretation of a gauge symbol $\height:M \to \RR_{\geqslant 0}$. Moreover, if $f$ is a function symbol with gauge modulus $\Delta_f$, then
    \[ (\forall r>0)(\forall x \in M)(\height(x)<r \implies \height(f(x)) \leqslant \Delta_{f}(r)). \]
    Similarly, for a predicate symbol $P$ with gauge modulus $\Delta_{P}$ we demand that
    \[ (\forall r>0)(\forall x \in M)(\height(x)<r \implies |P(x)| \leqslant \Delta_{P}(r)). \]
    At last, for each constant symbol $c$ we have a number $\Delta_c$ such that $\height(c) \leqslant \Delta_c$. If the domain of $f$ or $P$ has more than one variable, we define $\height:M^k \to \RR_{\geqslant 0}$ by taking maximum of heights of coordinates and impose the same condition. 
\end{definition}

\begin{example} \label{example_GVF_language}
    Let $\Field$ be equipped with a GVF structure, given by a global height function $h$. We equip it with an $\mathcal{L}$-structure, where $\mathcal{L}$ is the language of globally valued fields in the following way. For $x \in \Field$ we define
    \[ \height(x) := h[1:x]. \]    
    Let $e>0$ be such that $\height(2) \leqslant e$. Then the following inequalities hold:
    \begin{itemize}
        \item $\height(xy) \leqslant \height(x)+\height(y)$,
        \item $\height(x+y) \leqslant \height(x)+\height(y)+e.$
    \end{itemize}
    To see this one can either prove it directly, or use Corollary \ref{corollary_every_MVF_structure_comes_from_an_adelic_curve} together with \cite[Proposition 1.13]{dolce2023generalisation}. These inequalities show that the function $\Delta(r) = 2r+e$ is a valid gauge modulus for addition and multiplication symbols. For constants $0,1$ we take $\Delta = 0$. For a predicate symbol $h_n$, $\Delta_n(r) = 1 + n \cdot r$ can serve as a gauge modulus for $h_n$ interpreted as the height $h_n:\Field^n \to \{-1\} \cup \RR_{\geqslant 0}$ with $h_n(0) = -1$.
\end{example}

Note that here we put $h(0)=-1$ instead of $h(0)=-\infty$, because the value $-\infty$ is not allowed in unbounded continuous logic.

\begin{definition}
    Let $\mathcal{L}$ be a language. We define $\mathcal{L}$-terms inductively as follows:
    \begin{itemize}
        \item Any variable is an $\mathcal{L}$-term.
        \item Any constant $c$ from $\mathcal{L}$ is an $\mathcal{L}$-term.
        \item If $f$ is a function symbol from $\mathcal{L}$ of arity $n$, and $t_1,\ldots,t_n$ are $\mathcal{L}$-terms, then $f(t_1,\ldots,t_n)$ is an $\mathcal{L}$-term.
    \end{itemize}
    
   We then define $\mathcal{L}$-formulas inductively as follows:
    \begin{itemize}
        \item If $w_1, w_2$ are terms, then $w_1=w_2$ is an atomic formula.
        \item If $P$ is a predicate symbol from $\mathcal{L}$ of arity $n$, and $t_1,\ldots,t_n$ are $\mathcal{L}$-terms, then $P(t_1,\ldots,t_n)$ is an atomic $\mathcal{L}$-formula.
        \item As connectives we use continuous functions $\RR^n \to \RR$. So whenever $\varphi_1(x), \dots, \varphi_n(x)$ are formulas in variables $x$ and $u:\RR^n \to \RR$ is continuous, then $u(\varphi_1(x), \dots, \varphi_n(x))$ is a formula.
        \item Let $f:\RR_{\geqslant 0} \to \RR$ be a compactly supported continuous function and $\varphi(x,y)$ be a formula in tuples of variables $x, y$. Then we can form formulas of variables $y$ of the form:
        \[ \psi(y) = \inf_x f(\height(x)) \cdot \varphi(x,y), \ \xi(y) = \sup_x f(\height(x)) \cdot \varphi(x,y). \]
    \end{itemize}
    A formula is called quantifier-free if it does not use the quantifiers $\inf$ and $\sup$. A variable $x$ is said to occur freely in a formula $\phi$ if it has at least one occurrence which is not bound by a quantifier. A formula with no free variables is called a sentence. We denote by $F_n(\mathcal{L})$ (or $F_n$ if $\cL$ is clear from the context) the set of $\mathcal{L}$-formulas with $n$ free variables $x_1,\ldots,x_n$, and by $F_n^{\qf}(\mathcal{L})$ (or $F_n^{\qf}$) the set of quantifier-free formulas in $F_n(\mathcal{L})$. For each $\mathcal{L}$-structure $\mathfrak{M}$ one can inductively define the interpretation of any formula $\varphi(x)$ evaluated on a tuple $a$ from $\mathfrak{M}$ (denoted $\varphi(a)^{\mathfrak{M}}$ or just $\varphi(a)$). In the base case of this induction, for an atomic formula $w_1=w_2$ we define the value $(w_1=w_2)^\mathfrak{M}(a)$ as $0$ if $w_1(a) = w_2(a)$ in $\mathfrak{M}$ and $1$ otherwise. Moreover, when interpreting quantifiers, we put
    \[ (\inf_x f(\height(x)) \cdot \varphi(x,b))^\mathfrak{M} := \min(0,\inf_a f(\height(a)) \cdot \varphi(a,b)^\mathfrak{M} ), \]
    and similarly with supremum (with min replaced by max). If $\mathfrak{M}$ has elements of unbounded height, this coincides with the natural interpretation. This modification is needed for Łoś's theorem to hold.
    
    Moreover, for each term $t(x)$ and each formula $\varphi(x)$ we define inductively its gauge moduli $\Delta_t, \Delta_{\varphi}$ so that
    \[ (\forall L-\textnormal{structures }\mathfrak{M})(\forall a \in M)(\forall r>0)(\height(a)<r \implies \height(t(a)) \leqslant \Delta_{t}(r)), \]
    \[ (\forall L-\textnormal{structures }\mathfrak{M})(\forall a \in M)(\forall r>0)(\height(a)<r \implies |\varphi(a)| \leqslant \Delta_{\varphi}(r)). \]
    Variables have their gauge moduli set to be identities. Furthermore, we impose:
    \begin{itemize}
        \item $\Delta_{f(t_1,\dots,t_n)} = \Delta_f \circ \max_{i=1, \dots, n} \{ \Delta_{t_i} + 1 \}$,
        \item $\Delta_{(w_1=w_2)} = \id + 1$,
        \item $\Delta_{P(t_1,\dots,t_n)} = \Delta_P \circ \max_{i=1, \dots, n} \{ \Delta_{t_i} + 1 \}$,
        \item $\Delta_{u(\varphi_1, \dots, \varphi_n)}(r) = \sup \{ |u(p)| : p \in \RR^n, |p_i| \leqslant \Delta_{\varphi_i}(r)+1 \}$,
        \item If $R_f = \inf\{q \in \RR_{>0} : f(r) = 0 \textnormal{ for } r>q\}$ and $M_f = \sup_{p \in \RR} |f(p)|$, then we put $$\Delta_{\inf_x f(\height(x)) \cdot \varphi(x,y)}(r) = \Delta_{\sup_x f(\height(x)) \cdot \varphi(x,y)}(r) = \sup_{p \in \RR} |f(p)| \cdot (\Delta_{\varphi}(\max \{r, R_f\})+1).$$
    \end{itemize}
\end{definition}

\begin{definition}
    A theory $T$ in a language $\mathcal{L}$ is a set of conditions of the form $\varphi = 0$ for sentences $\varphi$. For an $\mathcal{L}$-structure $\mathfrak{M}$ we say that $\mathfrak{M}$ is a model of $T$ (written $\mathfrak{M} \models T$) if in $\mathfrak{M}$ the conditions $\varphi^\mathfrak{M} = 0$ from $T$ hold. 
\end{definition}

We also allow conditions of the form $\varphi \in C$ where $C$ is a closed subset of $\RR$, as this gives an equivalent notion.

\begin{example} \label{example_translating_axioms_to_continuous_logic}
    Assume that $\mathcal{L}$ is a language and let $f$ be a binary function symbol, $P, Q$ be continuous binary predicates, and $c$ be a constant. Let us describe how to write the statement
    \[ (\forall x,y) \bigl( (x\neq c \vee y \neq c) \implies P(x,y) \leqslant Q(x,y) \bigr) \]
    as a theory $T$, i.e., how to find a theory $T$ such that the models of $T$ are exactly the $\mathcal{L}$-structures satisfying (in the naive sense) this statement.

    Let $f_n:\RR_{\geqslant 0} \to \RR$ be a sequence of continuous non-increasing functions such that $f_n|_{[0,n]} = 1$ and $f_n|_{[n+1,\infty)} = 0$. For all natural $n$ we consider the condition
    \[ \bigl[ \sup_{x,y} f_n(\height(x,y)) \cdot \varphi_n(x,y) \bigr] = 0, \]
    where
    \[ \varphi_n(x,y) = \min \bigl( \max\{(x=c), (y=c)\}, \max \{P(x,y)-Q(x,y),0\}  \bigr). \]
    The set $T$ of all such conditions does the job, as if there was a model $\mathfrak{M}$ of $T$ with $a, b \in M$ such that $a\neq c$ or $b \neq c$ in $\mathfrak{M}$ and $P(a,b) > Q(a,b)$, then for $n$ such that $\height(a), \height(b) < n$ the formula 
    \[ \sup_{x,y} f_n(\height(x,y)) \cdot \varphi_n(x,y) \]
    would be evaluated as at least $P(a,b) - Q(a,b)$ which is greater than zero.
\end{example}

\begin{definition}\label{definition:GVF_formal_def}
    We denote by $\GVF_e$ the theory of globally valued fields in the language of globally valued fields with $\height(2) \leqslant e$, see Example \ref{example_GVF_language}. It consists of the field axioms together with the following (for every natural $n,m \geqslant 1$):
\[
    \begin{array}{ll}
        \forall x_0, \dots, x_n & h_n(x_1, \dots, x_n) = -1 \Leftrightarrow x_1=\dots=x_n=0 \\
        & h_2(1,1) = 0 \ \& \ h_2(2,1) \leqslant e \\
        \textnormal{For all } \sigma\in \Sym_n, \, \forall x_0,\dots,x_n & h_n(x_1,\dots,x_n) = h_n(x_{\sigma(1)}, \dots, x_{\sigma(n)}) \\
        \forall^* x_1,\dots,x_n\, \forall^* y_1,\dots,y_m & h_{nm}(x_1 y_1, \dots,x_n y_m) = h_n(x_1,\dots, x_n) + h_m(y_1,\dots, y_m) \\
        \forall x_1,\dots,x_n\, \forall y & h_n(x_1,\dots,x_n) \leqslant h_{n+1}(x_1,\dots,x_n,y) \\
        \forall^* x_1,\dots,x_n \, \forall^* y_1,\dots,y_n & h_n(x_1+y_1, \dots, x_n+y_n) \leqslant h_{2n}(x_1,\dots,x_n,y_1,\dots,y_n) + \aerror \\
        \forall^* x & h_1(x) = 0 \ \& \ \height(x) = h_2(x,1)
    \end{array}
\]
where $\forall^* z_1, \dots z_n$ means that we require it for all $z_1, \dots, z_n$, unless all of them are zero.
More precisely, one translates these axioms to conditions of the form $\varphi = 0$, exactly as in Example \ref{example_translating_axioms_to_continuous_logic}. We declare $\GVF$ to be $\GVF_1$. Moreover, for a $\QQ$-tropical polynomial $t$, we introduce a predicate $R_t$, using heights $h_n$, monomials, and Lemma \ref{lemma_tropical_polynomials_are_differences_of_maxima}.

\end{definition}

\begin{construction}
    Fix a language $\mathcal{L}$. Let $I$ be a set, $\cU$ be an ultrafilter on $I$ and $\mathfrak{M}_i$ be $\mathcal{L}$-structures for $i \in I$. We define the ultraproduct $\mathfrak{M} := \prod \mathfrak{M}_i/\cU$ as the quotient $M := N_0/\sim$ where
    \[ N_0 = \{ (a_i)_{i \in I} \in \prod_{i \in I} M_i : \height(a_i) \textnormal{ is bounded} \} \]
    and $(a_i)_{i \in I} \sim (b_i)_{i \in I}$ if the set of indices $i$ where $a_i = b_i$ is in $\cU$. We equip $\mathfrak{M}$ with an $\mathcal{L}$-structure interpreting constants as classes of constants in the product and functions as classes of values of functions in the product. For a continuous predicate symbol $P$, we interpret it by taking:
    \[ P([a])^{\mathfrak{M}} := \lim_{i \to \cU} P(a(i))^{\mathfrak{M}_i}, \]
    where $[a_i]$ is the class of $(a(i))_{i \in I} \in N_0$ and $\lim_{i \to \cU}$ is the ultralimit with respect to $\cU$. These definitions make sense by Definition \ref{definition_L_structure_unbdd_cont_logic}.
\end{construction}

\begin{theorem}(Łoś's theorem)
    Let $\varphi$ be an $n$-ary $\mathcal{L}$-formula and $\mathfrak{M}_i$ be $\mathcal{L}$-structures for $i \in I$. 
    Fix an ultrafilter $\cU$ on $I$, let $\mathfrak{M} := \prod \mathfrak{M}_i/\cU$, and $[a_1], \dots, [a_n] \in \mathfrak{M}$. Then,
    \[ \varphi([a_1], \ldots, [a_n])^{\mathfrak{M}} = \lim_{i \to \cU} \varphi(a_1(i), \ldots, a_n(i))^{\mathfrak{M}_i}. \]
\end{theorem}
\begin{proof}
    The proof is a standard induction over the complexity of the formula $\varphi$. We only describe how to deal with the quantifier case. 

    Consider a formula $\varphi(x) = \inf_y f(\height(y)) \cdot \psi(x,y)$. Fix a tuple $[a] = ([a_1],\dots, [a_n])$ in $M$. First, assume that $\varphi([a]) < r$ and that the minimum in the inductive interpretation is not attained at $0$ (the other case is easy). This means that there is a tuple $[b]$ in $M$ such that $f(\height([b])) \cdot \psi([a],[b]) < r$. By continuity of $f$ and the inductive hypothesis, we get that
    \[ \lim_{i \to \cU} f(\height b(i)) \cdot \psi(a(i),b(i)) < r. \]
    Thus
    \[ \lim_{i \to \cU} \inf_y f(\height(y)) \cdot \psi(a(i),y) < r. \]
    On the other hand, assume that $\lim_{i \to \cU} \varphi(a(i)) < r$. Without loss of generality, assume that for all $i \in I$ there is $b(i)$ in $M_i$ such that $f(\height b(i)) \cdot \psi(a(i),b(i)) < r$. 
    
    We end up with two cases based on whether $\height(b(i))$ is bounded or not.
    In the bounded case, there is a tuple $[b]$ in $M$ such that $f(\height([b])) \cdot \psi([a],[b]) < r$, which proves that $\varphi([a]) < r$. In the unbounded case, since $f$ is compactly supported, we get that $0 < r$. By the definition of the interpretation of quantifiers, we are done.
\end{proof}

\begin{example} \label{example:ultraproduct_transfer_number_fields_function_fields}
    Pick a sequence of real positive numbers $(r_i)_{i<\omega}$ converging to zero. Let $\Field$ be the ultraproduct $\prod_i \GVFQ[r_i]/\cU$ with respect to an ultrafilter $\cU$ on $\omega$. Then $\Field$ is a characteristic zero, algebraically closed field, with the GVF structure trivial on $\GVFQ$. One can also consider a finite-to-zero characteristic variant, i.e., $K = \prod_p \ov{\mathbb{F}_p(t)}[1]/\cU$, which is again, a characteristic zero, algebraically closed GVF, with the GVF structure trivial on $\GVFQ$.
\end{example}

\begin{example} \label{example:ultraproduct_transfer_Nevalinna_function_fields}
    Jensen's formula (\ref{eq:two}):
    \begin{equation*}
        \sum_{0<|a|<r} \ord_a(f) \log \frac{r}{|a|} + \frac{1}{2\pi} \int_{0}^{2\pi} -\log |f(r e^{i\theta})| d\theta = -\log |f(0)|
    \end{equation*}
    can be used to define a GVF structure on a subfield of meromorphic functions $\cM$. Indeed, fix an ultrafilter $u$ on $\RR_{>0}$ that avoids finite Lebesgue measure sets. Let $\eta(r)$ be a function that converges to $\infty$ when $r \to \infty$. Put
    \[ \height_{\eta, r}(f) := \sum_{0<|a|<r} \max(0,\ord_a(f)) \frac{\log(r/|a|)}{\eta(r)} + \frac{1}{2\pi} \int_{0}^{2\pi} \max(0,-\log |f(r e^{i\theta})|) \frac{d\theta}{\eta(r)} \]
    and $\height_{\eta, u}(f) := \lim_{r \to u} \height_{\eta, r}(f)$. Then we define
    \[ \cM[\eta,u] = \{ f \in \cM : \height_{\eta, u}(f) < \infty \} \]
    and equip it with a GVF structure by defining heights $h:\PP^n(\cM[\eta,u]) \to \RR$ similarly. By existential closedness of $\ov{\CC(t)}[1]$ (see Section \ref{section_existential_closedness}), the universal theories of $\cM[\eta,u]$ and $\ov{\CC(t)}[1]$ coincide. This formalises a rudimentary part of Vojta's dictionary between number theory and value distribution theory, and sets a goal of formalizing more.
\end{example}

\begin{definition}
    A set of conditions $\Sigma = \{ \varphi_i \leqslant r_i : i \in I \}$ is \textit{approximately finitely satisfiable} if for all finite $I_0 \subset I$ and $\varepsilon > 0$ the set of conditions $\{ \varphi_i \leqslant r_i + \varepsilon : i \in I_0 \}$ is satisfiable.
\end{definition}

\begin{corollary}(Compactness) \label{corollary_compactness}
    Let $\Sigma$ be a set of conditions in variables $(x_{\alpha})_{\alpha < \kappa}$, and $(r_\alpha)_{\alpha<\kappa}\in \RR_{\geqslant 0}^{\kappa}$. Assume that $\Sigma' = \Sigma \cup \{ \height(x_{\alpha}) \leqslant r_{\alpha} : \alpha < \kappa\}$ is approximately finitely satisfiable. Then, it is satisfiable. Furthermore, assume that for each finite subset $\sigma$ of $\Sigma'$ and $\varepsilon > 0$, $\sigma$ is satisfiable up to $\varepsilon$ in some model $\mathfrak{M}_{\sigma,\varepsilon}\models T$. Then, $\Sigma'$ is satisfiable in some ultraproduct of the $\mathfrak{M}_{\sigma,\varepsilon}$.
\end{corollary}
\begin{proof}
Write $\Sigma' = \{\varphi_i \leqslant r_i : i\in I\}$, and let $P_{\omega}(I)$ be the set of finite subsets of $I$. For each $I_0\in P_{\omega}(I)$ and $\varepsilon > 0$, let $\mathfrak{M}_{I_0,\varepsilon}$ be a model in which $\{\varphi \leqslant r_i + \varepsilon : i\in I_0 \}$ is satisfiable, and let $\hat{I_0} \subset P_{\omega}(I)$ be the set of $I_1$ such that $I_0\subset I_1$. Then, the set
\[
    \{ \hat{I_0}\times (0,\varepsilon) : I_0\in P_\omega(I), \varepsilon\in (0,1)\}
\]
is stable under finite intersections, so we can extend it to an ultrafilter $\mathcal{U}$ on $P_{\omega}(I)\times (0,1)$. 

For all $I_0\in P_{\omega}(I)$ and $\varepsilon > 0$, let $\tilde{I_0}\in P_{\omega}(I)$ be the finite set of indices of the conditions $\height(x_{\alpha}) \leqslant r_\alpha$, for $\alpha$ such that $x_{\alpha}$ occurs freely in $\varphi_i$ for some $i\in I_0$. Now, let
\[
    \mathfrak{M}^* := \prod \mathfrak{M}_{I_0\cup \tilde{I_0},\varepsilon} / \mathcal{U}.
\]
Define a sequence $([a_{\alpha}])_{\alpha<\kappa}$ in $\mathfrak{M}^*$ in the following way: for all $I_0\in P_{\omega}(I)$ and $\varepsilon > 0$, set $([a_{\alpha}(I_0,\varepsilon)])_{\alpha < \kappa}$ to be a realisation in $\mathfrak{M}_{I_0\cup \tilde{I_0},\varepsilon}$ of $\{\varphi \leqslant r_i + \varepsilon : i\in I_0\cup \tilde{I_0} \}$. Since we may choose any value of $a_{\alpha}$ if $x_{\alpha}$ does not occur freely in $\varphi$, we set $a_{\alpha} = 0$ in this case. Thus, for all $\alpha < \kappa$, and for all $(I_0,\varepsilon) \in P_{\omega}(I)\times (0,1)$, we have $\height(a_{\alpha}(i)) \leqslant r_\alpha+1$, which ensures that $[a_{\alpha}]$ is indeed an element of $\mathfrak{M}^*$.
Now, let $i\in I$, and $\varepsilon > 0$, the set $\widehat{\{i\}}\times (0,\varepsilon)$ is in $\mathcal{U}$, therefore by Łoś's theorem, $\varphi([a_\alpha] : \alpha < \kappa)^{\mathfrak{M}^*}\leqslant r_i + \varepsilon$. Since this is true for all $\varepsilon > 0$, $([a_{\alpha}]_{\alpha < \kappa}$ is indeed a realisation of $\Sigma'$. 
\end{proof}

\begin{construction}
    Let $T$ be a theory in a language $\mathcal{L}$. We define the space of $n$-types $S_n(T)$ as the set of functions $f:F_n(\mathcal{L}) \to \RR$ such that there is $\mathfrak{M} \models T$ and $a \in M^n$ such that $f(\varphi) = \varphi(a)^{\mathfrak{M}}$ for all $\varphi \in F_n$. We equip it with the pointwise convergence topology. We write $S_{n, \height \leqslant r}(T)$ for the set of types with $\height(x_i) \leqslant r$ for $i=1, \dots, n$. Note that
    \[ S_n(T) = \bigcup_{r \geqslant 0} S_{n, \height \leqslant r}(T) \]
    and each $S_{n, \height \leqslant r}(T)$ is compact (by the compactness theorem). Moreover, $S_n(T)$ is locally compact, and equipped with the colimit topology with respect to the above union.

    Similarly, we define the space of quantifier-free types $S_n^{\qf}(T)$ as the set of functions $f:F_n^{\qf}(\cL) \to \RR$ which are satisfiable in some model of $T$. Note that there is a restriction map $S_n(T) \to S_n^{\qf}(T)$ (where the latter is the set of quantifier-free $n$-types).
\end{construction}

If $\mathfrak{M}$ is a model of $T$, we define $\mathcal{L}_M$ as the language obtained by adding constants $(c_a)_{a\in M}$ to $\mathcal{L}$. Then, $\mathfrak{M}$ is naturally an $\mathcal{L}_M$-structure by interpreting $c_a$ as $a$. An $\mathcal{L}_M$-formula may also be called an $\mathcal{L}$-formula with parameters in $M$. Let $T'$ be the theory in $\mathcal{L}_M$ consisting of $T$ and the atomic diagram of $\mathfrak{M}$ (i.e., all quantifier-free sentences about tuples from $M$, written using constants $(c_a)_{a \in M}$). Note that models of $T'$ correspond to inclusions $\mathfrak{M}\subset \mathfrak{N}$, where $\mathfrak{N}$ is a model of $T$.
If $T$ is clear from the context, we write $S_n(\varnothing)$ for $S_n(T)$ and if $\mathfrak{M} \models T$, we write $S_n(M) := S_n(T')$ with $T'$ as above (we abuse notation slightly, as $S_n(M)$ depends depends on the structure $\mathfrak{M}$, not only on its underlying set). Similarly, we define $S_n^{\qf}(M) := S_n^{\qf}(T')$. For $a\in M$, we define $\tp(a/M)\in S_n(M)$ and $\qftp(a/M)\in S_n^{\qf}(M)$ to be the respective functions $\varphi \mapsto \varphi(a)^{\mathfrak{M}}$.

\begin{definition}
    Let $\mathfrak{M} \subset \mathfrak{N}$ be an extension of models of $T$. A function $f : N^n \rightarrow \RR$ is called \emph{$M$-definable} if there exists a continuous function $\tilde{f} : S_n(M) \rightarrow \RR$ such that
    \[
        \forall a\in N, \, f(a) = \tilde{f}(\tp(a/M)).
    \]
    Similarly, $f$ is called \emph{quantifier-free $M$-definable} if there exists a continuous function $\tilde{f} : S_n^{\qf}(M) \rightarrow \RR$ such that
    \[
        \forall a\in N, \, f(a) = \tilde{f}(\qftp(a/M)).
    \]
    Note that in both cases $f$ can be defined by the same formula on any $\mathfrak{M} \subset \mathfrak{N}' \models T$.

    Similarly, a subset $X\subset N^n$ is called (quantifier-free) $M$-definable if its indicator function is (quantifier-free) $M$-definable.
\end{definition}

Note that we can treat $\varphi \in F_n(\mathcal{L}_M)$ as an $M$-definable function by putting $\varphi(p) := p(\varphi)$ for $p \in S_n(M)$. The induced map $F_n(\mathcal{L}_M) \to C(S_n(M))$ has a dense image in the compact-open topology. If two formulas $\varphi_1, \varphi_2 \in F_n(\mathcal{L})$ (without parameters) define the same function in all models of $T$, we say that they are $T$-equivalent. 

We could define type spaces over models of the universal part of $T$, but in the GVF context that does not yield a more general notion, as $\GVF_e$ is a universal theory, at least if we add the inverse symbol to the language.

\begin{remark} \label{remark_qftypes_on_V_correspond_to_GVF_extensions}
    Let $K$ be a GVF and let $K_0$ be the underlying field, treated as a structure in the language $\mathcal{L}_{\textnormal{rings}} = (0,1,+,\cdot)$. Then, the inclusion of the set of $\mathcal{L}_{\textnormal{rings}}$-formulas into the set of formulas over the language of globally valued fields induces a restriction map $\pi:S_n^{\qf}(K) \to S_n^{\qf}(K_0)$, mapping $\qftp(a/K)$ to the point $\qftp(a/K_0) \in S_n^{\qf}(K_0)$ corresponding to the irreducible subvariety $\loc(a/K) = V \subset \Affine_K^n$ defined by polynomials that vanish on $a$. Since $a$ generates $K(V)$ over $K$, $\qftp(a/K)$ is uniquely determined by $\qftp(a/K_0)$ and the heights of tuples from $K(V)$. Thus, there is a bijection between the fiber of $\pi$ over $V$ and the set of GVF structures on $K(V)$ extending the one on $K$.
\end{remark}

\begin{remark}
    It is worth mentioning that the only quantifier-free definable sets of bounded height in GVFs are the algebraically constructible sets. More precisely, if $K$ is a GVF and $X \subset S_{n, \height \leqslant r}^{\qf}(K)$ is a clopen subset, then it is induced by some constructible subset of $K^n$. Indeed, this is because the set of quantifier-free GVF types extending a given quantifier-free ACF type over $K$ is convex, and hence connected. Moreover, since $-\min(h(x),0)$ defines the characteristic function of $\{0\}$ in a GVF, the restriction map from quantifier-free GVF types to ACF types is continuous.
\end{remark}

\section{Existential closedness} \label{section_existential_closedness}

In this section we present a few notions from discrete logic in the unbounded continuous setup. We follow \cite[Chapter 3]{chang1990model}. Let $\mathcal{L}$ be a (unbounded continuous logic) language and let $T$ be an $\mathcal{L}$-theory. Consider an extension $\mathfrak{M} \subset \mathfrak{N}$ of models of $T$. 
\begin{definition} \label{definition_existential_closedness}
    Let $\mathfrak{M}\subset \mathfrak{N}$ be an inclusion of $\mathcal{L}$-structures. We write $\mathfrak{M} \prec_{\exists} \mathfrak{N}$ if $\mathfrak{M}$ is \textit{existentially closed in $\mathfrak{N}$}, which means that if $\varphi(x,y)$ is a quantifier-free formula in $\mathcal{L}$ and $b \in M^{|y|}$, then for all $\varepsilon >0, r>0$ the following holds:
    \[ (\exists a \in N^{|x|}: \height(a) \leqslant r)(\varphi(a,b)^{\mathfrak{N}} = 0) \implies (\exists a' \in M^{|x|}: \height(a') \leqslant r + \varepsilon)(|\varphi(a',b)^{\mathfrak{M}}| \leqslant \varepsilon). \]
    Equivalently, this means that for any formula $\varphi(x,y)$ valued in $[0,1]$ and any natural $n$, if
    \[ \mathfrak{N} \models \bigl( (\sup_{x} f_n(\height(x)) \cdot (1 - \varphi(x,b) \bigr) = 1, \]
    then
    \[ \mathfrak{M} \models \bigl( (\sup_{x} f_n(\height(x)) \cdot (1 - \varphi(x,b) \bigr) = 1, \]
    where $f_n$ are the functions from Example \ref{example_translating_axioms_to_continuous_logic}.
\end{definition}
The $\varepsilon$ shows up in the definition, because with it we get the following lemma directly generalising the discrete logic case.
\begin{lemma} \label{lemma_existential_closedness_equivalent_conditions}
    Let $\mathfrak{M}\subset \mathfrak{N}$ be an inclusion of $\mathcal{L}$-structures.
    The following conditions are equivalent:
    \begin{enumerate}
        \item $\mathfrak{M} \prec_{\exists} \mathfrak{N}$.
        \item In $S_n^{\qf}(M)$, the set of quantifier-free types realised in $\mathfrak{M}$ is dense in the set of quantifier-free types realised in $\mathfrak{N}$ for all $n\in\mathbb{N}$.
        \item For some ultrapower $\mathfrak{M}^*$ of $\mathfrak{M}$, there is an embedding $\mathfrak{N} \hookrightarrow \mathfrak{M}^*$ over $\mathfrak{M}$.
    \end{enumerate}
\end{lemma}
\begin{proof}
    $\textnormal{(3)} \Rightarrow \textnormal{(1)}$ follows directly from Łoś's theorem.

    $\textnormal{(1)} \Rightarrow \textnormal{(3)}$ Let $(a_\alpha)_{\alpha < \kappa}$ be an enumeration of $\mathfrak{N}$, $r_\alpha = \height(a_\alpha)$, and let $\Sigma$ be the set of conditions of the form $\varphi_i = 0$ in the variables $(x_\alpha)_{\alpha < \kappa}$, with parameters in $M$, which are satisfied by the $a_{\alpha}$. Then, by $\textnormal{(1)}$, $\Sigma$ is finitely approximately satisfiable in $\mathfrak{M}$, so by Corollary \ref{corollary_compactness}, $\Sigma$ is satisfiable in an ultrapower of $\mathfrak{M}$.
    
    $\textnormal{(1)} \Rightarrow \textnormal{(2)}$ Let $a\in N^n$, and let $U\subset S_n^{\qf}(M)$ be an open set containing $\qftp(a/M)$. Then, $U$ contains a basic open set which can be expressed as a finite intersection of conditions of the form $\{|\varphi_i| < \delta\}$, say for $i\leqslant m$, where all $\varphi_i$ are quantifier-free formulas with parameters in $M$. Then, Definition \ref{definition_existential_closedness} applied to the formula $\varphi = \max \{|\varphi_1-\varphi_1(a)^{\mathfrak{N}}|,\ldots,|\varphi_m-\varphi_m(a)^{\mathfrak{N}}| \}$ with $r = \height(a)$ and $\varepsilon$ small enough, gives the existence of an $a'\in M^n$ such that for all $i\leqslant m$, $|\varphi(a')^{\mathfrak{M}}| < \delta$. Hence, $\qftp(a'/M)\in U$, which proves $\textnormal{(2)}$.
    
    $\textnormal{(2)} \Rightarrow \textnormal{(1)}$ follows from the fact that the intersection of conditions $\height(x) < r+\varepsilon$ and $|\varphi(x,b)| < \varepsilon$ defines a neighborhood of $\qftp(a/M)$ in $S_n^{\qf}(M)$.  
\end{proof}

Similarly, we write $\mathfrak{M} \prec \mathfrak{N}$ and call the inclusion \textit{elementary} if for all tuples $a$ in $\mathfrak{M}$ and all formulas $\varphi(x)$ we have $\varphi(a)^{\mathfrak{M}} = \varphi(a)^{\mathfrak{N}}$. If the embedding $\mathfrak{N} \hookrightarrow \mathfrak{M}^*$ in the above lemma is elementary, then $\mathfrak{M} \prec \mathfrak{N}$ by Łoś's theorem (and also vice versa).

\begin{definition}
    An $\mathcal{L}$-structure $\mathfrak{M} \models T$ is \textit{existentially closed} (with respect to $T$) if for all extensions $\mathfrak{M} \subset \mathfrak{N} \models T$ we have $\mathfrak{M} \prec_{\exists} \mathfrak{N}$. We call $T$ \textit{model complete} if for all embeddings $\mathfrak{M} \subset \mathfrak{N}$ of models of $T$ we have $\mathfrak{M} \prec \mathfrak{N}$. 
\end{definition}

Model complete theories also have the following characterisation.

\begin{lemma} \label{lemma_model_complete_theory_equivalent_conditions}
    Let $T$ be a theory. Then, $T$ is model complete if and only if all models of $T$ are existentially closed.
\end{lemma}
\begin{proof}
    The forward direction is clear, since $\mathfrak{M}\prec \mathfrak{N}$ implies $\mathfrak{M}\prec_{\exists}\mathfrak{N}$.

    Conversely, assume all models of $T$ are existentially closed, let $\mathfrak{M}\subset \mathfrak{N}$ be an inclusion of models of $T$, and let $\varphi$ be a sentence with parameters in $M$. We show that $\varphi^{\mathfrak{N}} = \varphi^{\mathfrak{M}}$ by induction on the complexity of $\varphi$.

    If $\varphi$ is quantifier-free, then by definition $\varphi^{\mathfrak{N}} = \varphi^{\mathfrak{M}}$.
    
    If $\varphi = \inf_x f(\height(x)) \cdot \psi(x)$, with $\psi$ such that $\psi(a)^{\mathfrak{N}} = \psi(a)^{\mathfrak{M}}$ for all $a\in M^n$, then $\varphi^{\mathfrak{N}} = \min(0,\inf_{a\in N^n} f(\height(a)) \cdot \psi(a)^{\mathfrak{N}}) \leqslant \varphi^{\mathfrak{M}}$. On the other hand, by Lemma \ref{lemma_existential_closedness_equivalent_conditions}, there is an ultraproduct $\mathfrak{M}^*$ of $\mathfrak{M}$ such that $\mathfrak{N}$ embeds into $\mathfrak{M}^*$ over $\mathfrak{M}$. The same argument as above shows that $\varphi^{\mathfrak{M}^*}\leqslant \varphi^{N}$. But, by Łoś's theorem, we have $\varphi^{\mathfrak{M}^*} = \varphi^{\mathfrak{M}}$, hence $\varphi^{\mathfrak{N}} = \varphi^{\mathfrak{M}}$.

    If $\varphi$ is of the form $\sup_x f(\height(x)) \cdot \psi(x)$, the same argument applied to $-\psi(x)$ shows that $\varphi^{\mathfrak{N}} = \varphi^{\mathfrak{M}}$, which concludes the proof.
\end{proof}

\begin{definition}
    We say that an $\mathcal{L}$-theory $T^*$ is a model companion of $T$ if $T^*$ is model complete and:
    \begin{itemize}
        \item all models of $T$ admit an embedding into a model of $T^*$,
        \item all models of $T^*$ admit an embedding into a model of $T$.
    \end{itemize}
\end{definition}

If a model companion of $T$ exists, it is unique. Moreover, similarly as in the discrete logic, there is a following characterisation of existence of a model companion.

\begin{lemma}
    Assume that $T$ is \textit{inductive}, i.e., that the union of a chain of models of $T$ is also a model of $T$. Let $K$ be the class of existentially closed models of $T$. Then the following are equivalent:
    \begin{enumerate}
        \item $K$ is an elementary class,
        \item $T$ has a model companion.
    \end{enumerate}
    Moreover, if these conditions are satisfied, the model companion is the common theory of $K$.
\end{lemma}
\begin{proof}
    This is a straightforward consequence of Lemma \ref{lemma_model_complete_theory_equivalent_conditions}.
\end{proof}

It is worth mentioning that if $T$ is inductive, existentially closed models exist and every model of $T$ can be embedded in one. Note that theories $\GVF_e$ are inductive, for $e \geqslant 0$. A very important question for the development of globally valued fields is the following conjecture.

\begin{conjecture} \label{conjecture_GVF_model_companion}
    The theory $\GVF_e$ has a model companion for $e \geqslant 0$.
\end{conjecture}

The affirmative answer to this question would pave the way for the use of model theoretic methods in studying heights. Immediately, it would follow that existentially closed globally valued fields are elementarily equivalent to one of the following (up to a scaling of the GVF structure): $\GVFQ$ with $\height(2) = \log(2)$, $\ov{\QQ(t)}$ with $\height(t)=1$, $\height|_{\QQ} = 0$, and $\ov{\mathbb{F}_p(t)}$ with $\height(t)=1$ for a prime $p$. This is because $\ov{\QQ(t)}$ and $\GVFQ[1]$ are existentially closed by \cite{GVF2, szachniewicz2023existential} respectively. The proof methods in loc.cit. and the recent developments of adelic curves \cite{liu2024arithmetic, A_Sedillot_diff_of_relative_volume} suggest a positive answer.

Let us note that existential closedness of the above GVFs and Lemma \ref{lemma_existential_closedness_equivalent_conditions} imply the following result, which due to Corollary \ref{corollary_every_MVF_structure_comes_from_an_adelic_curve} may shed some light on the study of adelic curves.

\begin{corollary} \label{corollary:GVFs_embedd_into_three_sorts}
    Every globally valued field $K$ embeds in an ultrapower of:
    \begin{itemize}
        \item either $\GVFQ[r]$ if $K$ is of characteristic zero and $r = \height(2) > 0$,
        \item or $\overline{\QQ(t)}[1]$ (with $\height(2)=0$) if $K$ is of characteristic zero and $\height(2) = 0$,
        \item or $\ov{\mathbb{F}_p(t)}[1]$ if $K$ is of characteristic $p$.
    \end{itemize}
\end{corollary}

The GVF in the middle bullet could be replaced by any GVF of the form $\ov{k(t)}[1]$ with $k$ of characteristic zero, for example by $\ov{\CC(t)}[1]$. 

\printbibliography

\end{document}